\documentclass[11pt,oneside]{article}
\usepackage{supertabular, setspace,hyperref}
\usepackage[utf8]{inputenc}
\usepackage{amsmath}
\usepackage{mathtools}
\usepackage{amsfonts}
\usepackage{amssymb}
\usepackage{amsthm}
\newtheorem{theorem}{Theorem}[section]

\newtheorem{prop}[theorem]{Proposition}
\newtheorem{example}[theorem]{Example}

\newtheorem{remark}{\sc Remark}
\newtheorem{lemma}{\sc Lemma}[section]
\newtheorem{corollary}{\sc Corollary}[section]
\newtheorem{definition}{\sc Definition}[section]

\newcommand{\be}{\begin{eqnarray}}
\newcommand{\ee}{\end{eqnarray}}
\newcommand{\Be}{\begin{eqnarray*}}
	\newcommand{\Ee}{\end{eqnarray*}}
\newcommand{\bee}{\begin{equation}}
\newcommand{\eee}{\end{equation}}
\newcommand{\ba}{\begin{array}}
	\newcommand{\ea}{\end{array}}
\newcommand{\bl}{\begin{lemma}}
	\newcommand{\el}{\end{lemma}}
\newcommand{\bd}{\begin{definition}}
	\newcommand{\ed}{\end{definition}}
\newcommand{\bt}{\begin{theorem}}
	\newcommand{\et}{\end{theorem}}
\newcommand{\bp}{\begin{proof}}
	\newcommand{\ep}{\end{proof}}
\newcommand{\bi}{\begin{itemize}}
	\newcommand{\ei}{\end{itemize}}
\newcommand{\br}{\begin{remark}}
	\newcommand{\er}{\end{remark}}
\newcommand{\bc}{\begin{corollary}}
	\newcommand{\ec}{\end{corollary}}
\newcommand{\bex}{\begin{example}}
	\newcommand{\eex}{\end{example}}
\usepackage{chngcntr}

\counterwithin*{equation}{section}
\counterwithin*{equation}{section}
\begin{document}
	\date{}
	\title{\textbf{Some classes of projectively and dually flat Finsler spaces with Randers change}}
	\maketitle
	\begin{center}
		\author{\textbf{Gauree Shanker, Sarita Rani and Kirandeep Kaur}}
	\end{center}
	\begin{center}
		Department of Mathematics and Statistics\\
		School of Basic and Applied Sciences\\
		Central University of Punjab, Bathinda, Punjab-151 001, India\\
		Email: gshankar@cup.ac.in, saritas.ss92@gmail.com, kiran5iitd@yahoo.com
\end{center}
\begin{center}
	\textbf{Abstract}
\end{center}
\begin{small}
	In this paper, we consider Randers change of some special $ (\alpha, \beta)- $  metrics. First we find the fundamental metric tensor and Cartan tensor of these Randers changed  $ (\alpha, \beta)- $metrics.  Next, we establish a general formula for inverse of fundamental metric tensors of these metrics. Finally, we find the necessary and sufficient conditions under which the Randers change of these $ (\alpha, \beta)- $  metrics are projectively and  locally dually flat.
\end{small}
\section{Introduction}
    According to S. S. Chern \cite{Chern1996}, Finsler geometry is just Riemannian geometry without quadratic restriction. Now a days, Finsler geometry is an interesting and active area of research for both pure and applied reasons \cite{Antonelli}.Though there has been a lot of development in this area, still there is a huge scope of research work in Finsler geometry. The concept of $(\alpha, \beta)-$metric was introduced by M. Matsumoto \cite{M.Mat1972} in 1972. For a general Finsler metric $F$, C. Shibata \cite{C.Sh1984} introduced the notion of $\beta-$change in $F,$ i.e., $\bar{F}=f(F, \beta)$ in 1984. Recall \cite{ChernShenRFG} that a Finsler metric $F$ on an open subset $\mathcal{U}$ of $\mathbb{R}^n$ is called projectively flat if and only if all the geodesics are straight in $\mathcal{U}$. The concept of dually flatness in Riemannian geometry was given by Amari and Nagaoka in \cite{S.I.Ama H.Nag} while studying information geometry. Information geometry provides mathematical science with a new framework for analysis. Information geometry is an investigation of differential geometric structure in probability distribution. It is also applicable in statistical physics, statistical inferences etc. Z. Shen \cite{ShenRFGAIG} extended the notion of dually flatness in Finsler spaces. Since then, many authors (\cite{SAB.GS.2016PFDF}, \cite{GS.SAB.2017PFDF}, \cite{GS.RKS.RDSK.2017PFDF}) have worked on this topic.\\
    The current paper is organized as follows:\\
    In second section,  we give basic definitions and examples of some special Finsler spaces with $(\alpha, \beta)-$metrics obtained by Randers change. In section 3, we find fundamental metric tensors $\bar{g}_{ij}$ and Cartan tensors $\bar{C}_{ijk}$ for these metrics. In section 4, we find a generalized formula for the inverse $\bar{g}^{ij}$ of fundamental metric tensor $\bar{g}_{ij}.$ In sections 5 and 6 we find the necessary and sufficient conditions for Randers change of some special Finsler spaces with $(\alpha, \beta)-$metrics to be projectively and locally dually flat respectively.

\section{Preliminaries}
Though there is vast literature available for Riemann-Finsler geometry, here we give some basic definitions, examples and results required for subsequent sections. 
\begin{definition} Let $M$ be an $n$-dimensional smooth manifold, $T_xM$ the tangent space at $x\in M,$ and $$TM:=\bigsqcup_{x\in M}T_xM $$ be the tangent bundle of $M$ whose elements are denoted by $(x,y),$ where $x\in M$ and $y\in T_xM.$ \\
A Finsler structure on $M$ is a function $$F:TM\longrightarrow[0,\infty), $$ with the properties:
\begin{itemize}
	\item \textbf{Regularity:} $F$ is $ C^{\infty}$ on the slit tangent bundle $ TM\backslash \left\lbrace 0\right\rbrace .$
	\item \textbf{Positive homogeneity:} 
	$ F(x,\lambda y)=\lambda F(x,y) \; \forall \; \lambda >0. $
	\item \textbf{Strong convexity:} The $n\times n$ matrix $\left( g_{ij}\right) =\left( \left[ \dfrac{1}{2}F^2 \right]_{y^iy^j} \right)  $ is positive-definite at every point of $ TM\backslash \left\lbrace 0\right\rbrace .$
\end{itemize}
A smooth manifold $M$ together with the Finsler structure $F,$ i.e., $(M,F)$ is called Finsler space and the corresponding geometry is called Finsler geometry.
\end{definition}
Next, we recall \cite{M.Mat1972} following definition:
\begin{definition}
An $(\alpha, \beta)$ metric $F=\alpha \phi\left( s\right),\ s=\dfrac{\beta}{\alpha}$ with $\phi=\phi(s) $ a smooth positive function on some symmetric interval $(-b_0, b_0),$ is a Finsler metric $F,$ which is positively homogeneous function of $\alpha$ and $\beta$ of degree one,  where $\alpha = \sqrt{a_{ij}(x) y^i y^j} $ is a Riemannian metric and  $\beta=b_i(x) y^i$ is a $1-$form.
\end{definition} 
\begin{lemma}\cite{ChernShenRFG}
An $(\alpha, \beta)-$norm $F=\alpha \phi\left( \dfrac{\beta}{\alpha}\right) $ is said to be a Minkowski norm for any Riemannian metric $\alpha$ and $1-$form $\beta$ with $ \lVert\beta\rVert_\alpha <b_0 $ if and only if $\phi=\phi(s)$ satisfies the following conditions:
$$ \phi(s)>0,\ \phi(s)-s\phi^{'}(s)+\left( b^2-s^2\right)\phi^{''}(s)>0,\ \phi(s)-s\phi^{'}(s)>0,  $$
where $b$ is an arbitrary numbers satisfying $ \lvert s\rvert\leq b<b_0.$
\end{lemma}
There are so many classical examples  of $(\alpha, \beta)-$metrics, below we mention few of them:\\
Randers metric, Kropina metric, generalized Kropina metric, Z. Shen's square metric, Matsumoto metric, exponential metric, infinite series metric.\\
Recall \cite{C.Sh1984} that
\begin{definition}
Let $(M, F)$ be an $n-$dimensional Finsler space. Then a metric $\bar{F}=f(F, \beta)+\beta $ constructed via a $\beta-$change is called Randers change of $(\alpha, \beta)-$metric.	
\end{definition}
Next, we construct some special Finsler metrics via Randers change of $(\alpha, \beta)$ metrics. Our further studies will be based on these metrics.
\begin{enumerate}
	\item Kropina-Randers change of $\left( \alpha,\beta\right) $-metric:\\
	We know that $F=\dfrac{\alpha^2}{\beta},$ where $\alpha = \sqrt{a_{ij}(x) y^i y^j} $ is a Riemannian metric and  $\beta=b_i(x) y^i$ is a $1-$form, is Kropina metric. Applying  $\beta-$change and Randers change simultaneously to this metric, we obtain a new metric
	$ \bar{F}= \dfrac{F^{2}}{\beta} + \beta,$ which we call Kropina-Randers changed metric.
	
	\item Generalized Kropina-Randers change of $\left( \alpha,\beta\right) $-metric:\\
	We know that $F=\dfrac{\alpha^{m+1}}{\beta^m} \ \ \left( m\neq 0, -1\right) ,$ where $\alpha = \sqrt{a_{ij}(x) y^i y^j} $ is a Riemannian metric and  $\beta=b_i(x) y^i$ is a $1-$form, is generalized Kropina metric. Applying $\beta-$change and Randers change simultaneously to this metric, we obtain a new metric
	$ \bar{F}= \dfrac{F^{m+1}}{\beta^m} + \beta, \ \ m\neq 0, -1, $  which we call generalized Kropina-Randers changed metric.
	
	\item Square-Randers change of $\left( \alpha,\beta\right) $-metric:\\
	We know that $F=\dfrac{(\alpha+\beta)^{2}}{\alpha},$ where $\alpha = \sqrt{a_{ij}(x) y^i y^j} $ is a Riemannian metric and  $\beta=b_i(x) y^i$ is a $1-$form, is Z. Shen's square metric. Applying $\beta-$change and Randers change simultaneously to this metric, we obtain a new metric
	$ \bar{F}= \dfrac{(F+\beta)^{2}}{F} + \beta,$  which we call Square-Randers changed metric.
	
	\item Matsumoto-Randers change of $\left( \alpha,\beta\right) $-metric:\\
	Again $F=\dfrac{\alpha^{2}}{\alpha-\beta},$ where $\alpha = \sqrt{a_{ij}(x) y^i y^j} $ is a Riemannian metric and  $\beta=b_i(x) y^i$ is a $1-$form, is a well known Matsumoto metric. Applying $\beta-$change and Randers change simultaneously to this metric, we obtain a new metric
	$ \bar{F}= \dfrac{F^{2}}{F-\beta} + \beta,$  which we call Matsumoto-Randers changed metric.
	\item Exponential-Randers change of $\left( \alpha,\beta\right) $-metric:\\
	The metric $F=\alpha e^{ \beta/\alpha},$ where $\alpha = \sqrt{a_{ij}(x) y^i y^j} $ is a Riemannian metric and  $\beta=b_i(x) y^i$ is a $1-$form, is called exponential metric. Applying $\beta-$change and Randers change simultaneously to this metric, we obtain a new metric
	$ \bar{F}= F e^{ \beta/F} + \beta, $  which we call exponential-Randers changed metric.
	\item Randers change of infinite series  $\left( \alpha,\beta\right) $-metric:\\
	The metric $F=\dfrac{\beta^{2}}{\beta -\alpha},$ where $\alpha = \sqrt{a_{ij}(x) y^i y^j} $ is a Riemannian metric and  $\beta=b_i(x) y^i$ is a $1-$form, is called infinite series  metric. Applying $\beta-$change and Randers change simultaneously to this metric, we obtain a new metric
	$ \bar{F}= \dfrac{\beta^{2}}{\beta - F } + \beta, $  which we call infinite series-Randers changed  metric.
\end{enumerate}
\begin{definition}
Let $(M, F)$ be an $n-$dimensional Finsler space. If $ F=\sqrt{A}, $  where $ A=a_{i_1 i_2}(x) y^{i_1} y^{i_2} $ with  $ a_{i_1 i_2}(x) $ symmetric in both the indices, then $F$ is called square-root Finsler metric. 
\end{definition}
We use following  notations in the subsequent sections:
\begin{equation*}
\begin{split}
\dfrac{\partial \bar{F}}{\partial x^i}=\bar{F}_{x^i},\
\dfrac{\partial \bar{F}}{\partial y^i}=\bar{F}_{y^i},\
\dfrac{\partial A}{\partial x^i}=A_{x^i},\ 
\dfrac{\partial A}{\partial y^i}=A_{i},\ 
A_{x^i} y^i=A_0,\
A_{x^iy^j} y^i=A_{0j},\ \\
\dfrac{\partial \beta}{\partial x^i}=\beta_{x^i},\
\dfrac{\partial \beta}{\partial y^i}=b_{i}\  \text{or}\ \beta_i,\ 
\beta_{x^i} y^i=\beta_0,\
\beta_{x^iy^j} y^i=\beta_{0j}\  \text{etc.}
\end{split}
\end{equation*}
\section{Fundamental Metric Tensors and Cartan Tensors}
\begin{definition}
Let $(M,F)$ be an $n-$dimensional Finsler manifold. The function 
$$ \mathnormal{g}_{ij}:=\left( \dfrac{1}{2} F^2\right)_{y^i y^j}=FF_{y^i y^j}+F_{y^i}F_{y^j}={h}_{ij}+{\ell}_i \ell_j $$ is called fundamental tensor of the metric $F$.
\end{definition}
\begin{definition}
Let $(M,F)$ be an $n-$dimensional Finsler manifold. Then its Cartan tensor is defined as 
$$ C_{ijk}(y)=\dfrac{1}{2} \dfrac{\partial\mathnormal{g}_{ij}}{\partial y^k}=\dfrac{1}{4}\left(F^2 \right)_{y^i y^j y^k}, $$
which is symmetric in all the three indices $ i, j, k.$
\end{definition}
Next, we find fundamental metric tensor and Cartan tensor for all the Randers changed $(\alpha, \beta)-$metrics constructed in the previous section.\\
 Kropina-Randers changed metric is 
 \begin{equation}{\label{KR1.1}}
 \bar{F}= \dfrac{F^{2}}{\beta} + \beta.
 \end{equation}
Differentiating (\ref{KR1.1}) w.r.t. ${y^{i}},$ we get
\begin{equation}{\label{KR1.2}}
\bar{F}_{y^{i}} = \dfrac{2F}{\beta}F_{y^{i}}- \dfrac{F^2-\beta^2}{\beta^2} b_{i}.
\end{equation}
Differentiation of (\ref{KR1.2}) further w.r.t. ${y^{j}}$ gives
\begin{equation}{\label{KR1.3}}
 \bar{F}_{y^{i}y^{j}} =\dfrac{2}{\beta}\ g_{ij} -  \dfrac{2F}{\beta^2} \left( b_{i}F_{y^{j}}+b_{j}F_{y^{i}}\right)  + \dfrac{2 F^2}{\beta^3}\ b_{i} b_{j}.
\end{equation}
Now,
\begin{align*}
\bar{g}_{ij}&=\bar{F} \bar{F}_{y^{i}y^{j}} + \bar{F}_{y^{i}} \bar{F}_{y^{j}}\\
&=\left\lbrace \dfrac{F^{2}}{\beta} + \beta\right\rbrace \left\lbrace\dfrac{2}{\beta}\ g_{ij} -  \dfrac{2F}{\beta^2} \left( b_{i}F_{y^{j}}+b_{j}F_{y^{i}}\right)  + \dfrac{2 F^2}{\beta^3}\ b_{i} b_{j}\right\rbrace \\
&\ \ \ \  +\left\lbrace \dfrac{2F}{\beta}F_{y^{i}}- \dfrac{F^2-\beta^2}{\beta^2} b_{i}\right\rbrace \left\lbrace \dfrac{2F}{\beta}F_{y^{j}}- \dfrac{F^2-\beta^2}{\beta^2} b_{j}\right\rbrace.
\end{align*}
Simplifying, we get
\begin{equation}{\label{KR1.4}}
\bar{g}_{ij} = \dfrac{2(F^2+\beta^2)}{\beta^2}\ g_{ij} -  \dfrac{4F^3}{\beta^3} \left( b_{i}F_{y^{j}}+b_{j}F_{y^{i}}\right)  + \dfrac{4 F^2}{\beta^2} F_{y^{i}} F_{y^{j}}+\left( \dfrac{3F^4}{\beta^4}+1 \right)  b_{i} b_{j}.
\end{equation}
Hence, we have following:
\begin{prop}
 Let $(M, \bar{F})$ be an $n-$dimensional Finsler space with $\bar{F}= \dfrac{F^{2}}{\beta} + \beta $ as a Kropina-Randers changed metric. Then its fundamental metric tensor is given by equation (\ref{KR1.4}).
\end{prop}
Next, we find Cartan tensor for Kropina-Randers changed  metric.\\
By definition, we have
$$ 2C_{ijk}(y)=\dfrac{\partial\mathnormal{g}_{ij}}{\partial y^k}.$$
From the equation (\ref{KR1.4}), we get 
\begin{align*}
2 \bar{C}_{ijk}
=&\dfrac{\partial}{\partial y^k}\left( \bar{g}_{ij} \right) \\
=& \dfrac{2(F^2+\beta^2)}{\beta^2}\ (2C_{ijk})+ \dfrac{4\left(\beta y_{_k} - F^2 b_k\right)} {\beta^3}\left(h_{ij}+\dfrac{y_i y_j}{F^2} \right)
-\dfrac{4F^3}{\beta^3} \left( \dfrac{b_i h_{jk}+b_j h_{ik}}{F}\right)\\
&-\dfrac{12\left( \beta F y_{_k}-F^3 b_k \right) }{\beta^4}\left(\dfrac{b_i y_j+b_j y_i}{F} \right)
+ \dfrac{4 F^2}{\beta^2} \left(\dfrac{y_i h_{jk}+y_j h_{ik}}{F^2} \right)+\dfrac{8\left(\beta y_{_k}-F^2 b_k \right) }{\beta^3} \dfrac{y_i y_j}{F^2}\\
&+\dfrac{12\left( \beta F^2 y_{_k}-F^4 b_k \right) }{\beta^5} b_ib_j
\end{align*}
\begin{align*}
=& \dfrac{4(F^2+\beta^2)}{\beta^2} C_{ijk} +h_{ij}\left(\dfrac{-4F^2}{\beta^3}b_k+\dfrac{4}{\beta^2} y_{_k} \right) +h_{jk}\left(\dfrac{-4F^2}{\beta^3}b_i+\dfrac{4}{\beta^2} y_i \right) +h_{ki}\left(\dfrac{-4F^2}{\beta^3}b_j+\dfrac{4}{\beta^2} y_j \right)\\
&-\dfrac{12F^4}{\beta^5}\biggl\{ b_ib_jb_k-\dfrac{\beta}{F^2}\left( b_ib_ky_j+b_ib_jy_{_k}+b_jb_ky_i \right)-\dfrac{\beta^3}{F^6}y_iy_jy_{_k}
+\dfrac{\beta^2}{F^4}\left( y_iy_{_k}b_j+y_iy_jb_k+y_jy_{_k}b_i \right) \biggr\}
\end{align*}
\begin{align*}
=&\dfrac{4(F^2+\beta^2)}{\beta^2} C_{ijk}
-\dfrac{4F^2}{\beta^3} \sum_{\text{cyclic sum}} h_{ij}\left(b_k-\dfrac{\beta}{F^2} y_{_k} \right)-\dfrac{12F^4}{\beta^5} \prod_{\text{cyclic product}} \left(b_i-\dfrac{\beta}{F^2} y_i \right).
\end{align*}
After simplification, we get
\begin{equation}{\label{KR1.5}}
\bar{C}_{ijk}= \dfrac{2(F^2+\beta^2)}{\beta^2} C_{ijk}
-\dfrac{2F^2}{\beta^3} \left(  h_{ij}m_k+h_{jk}m_i+h_{ki}m_j\right) -\dfrac{6F^4}{\beta^5}\ m_im_jm_k.
\end{equation}
The above discussion leads to the following proposition. 
\begin{prop}
Let $(M, \bar{F})$ be an $n-$dimensional Finsler space, with $\bar{F}= \dfrac{F^{2}}{\beta} + \beta $ as a Kropina-Randers changed Finsler metric. Then its Cartan tensor is given by equation (\ref{KR1.5}).
\end{prop}
Next, we find fundamental metric tensor for generalized Kropina-Randers changed  metric
\begin{equation}{\label{KR2.1}}
 \bar{F}= \dfrac{F^{m+1}}{\beta^m} + \beta, \ \ m\neq 0, -1.
\end{equation} 
Differentiating (\ref{KR2.1}) w.r.t. ${y^{i}},$ we get
\begin{equation}{\label{KR2.2}}
\bar{F}_{y^{i}} =(m+1) \dfrac{F^m}{\beta^m}F_{y^{i}} + \left(1-m  \dfrac{F^{m+1}}{\beta^{m+1}}\right)  b_{i}.
\end{equation}
Differentiation of (\ref{KR2.2}) further w.r.t. ${y^{j}}$ gives
\begin{equation}{\label{KR2.3}}
\bar{F}_{y^{i}y^{j}} =(m+1) \dfrac{F^m}{\beta^m} F_{y^{i} y^{j}} +m(m+1)  \dfrac{F^{m-1}}{\beta^{m}}F_{y^{i}}F_{y^{j}}-m(m+1)  \dfrac{F^{m}}{\beta^{m+1}} \left( b_i F_{y^{j}} +b_jF_{y^{i}}\right) +m(m+1)  \dfrac{F^{m+1}}{\beta^{m+2}} b_i b_j.
\end{equation}
Now,
\begin{align*}
\bar{g}_{ij}&=\bar{F} \bar{F}_{y^{i}y^{j}} + \bar{F}_{y^{i}} \bar{F}_{y^{j}}\\
=&\left\lbrace \dfrac{F^{m+1}}{\beta^m} + \beta\right\rbrace \biggl\{ (m+1) \dfrac{F^m}{\beta^m} F_{y^{i} y^{j}} +m(m+1)  \dfrac{F^{m-1}}{\beta^{m}}F_{y^{i}}F_{y^{j}}-m(m+1)  \dfrac{F^{m}}{\beta^{m+1}} \left( b_i F_{y^{j}} +b_jF_{y^{i}}\right) {}\\
& +m(m+1)  \dfrac{F^{m+1}}{\beta^{m+2}} b_i b_j\biggr\}\\
&+  \left\lbrace (m+1) \dfrac{F^m}{\beta^m}F_{y^{i}} + \left(1-m  \dfrac{F^{m+1}}{\beta^{m+1}}\right)  b_{i}\right\rbrace \biggl\{ (m+1) \dfrac{F^m}{\beta^m}F_{y^{j}}+ \left(1-m  \dfrac{F^{m+1}}{\beta^{m+1}}\right)  b_{j}\biggr\} 
\end{align*}
After simplification, we get
\begin{equation}{\label{KR2.4}}
\begin{split}
\bar{g}_{ij}
=&(m+1) \left( \dfrac{F^{2m}}{\beta^{2m}}+ \dfrac{F^{m-1}}{\beta^{m-1}} \right)  g_{ij} -(m+1)\left( 2m \dfrac{F^{2m+1}}{\beta^{2m+1}}+(m-1) \dfrac{F^{m}}{\beta^{m}} \right) \left(b_{i}F_{y^{j}}+b_{j}F_{y^{i}}\right) \\
+& (m+1)\left( 2m \dfrac{F^{2m}}{\beta^{2m}}+(m-1) \dfrac{F^{m-1}}{\beta^{m-1}} \right)  F_{y^{i}} F_{y^{j}}
+\left( 1+m(2m+1) \dfrac{F^{2m+2}}{\beta^{2m+2}}+m(m-1) \dfrac{F^{m+1}}{\beta^{m+1}} \right)b_{i}b_{j}.
\end{split}
\end{equation}
Hence, we have the following proposition.
\begin{prop}
Let $(M, \bar{F})$ be an $n-$dimensional Finsler space with $\bar{F}= \dfrac{F^{m+1}}{\beta^m} + \beta \ \ \left( m\neq 0, -1\right) $ as generalized Kropina-Randers changed metric. Then its fundamental metric tensor is given by equation (\ref{KR2.4}).
\end{prop}
Next, we find Cartan tensor for generalized Kropina-Randers changed  metric.\\
By definition, we have
$$2C_{ijk}(y)=\dfrac{\partial\mathnormal{g}_{ij}}{\partial y^k}.$$
From the equation (\ref{KR2.4}), we get \\
\begin{align*}
2 \bar{C}_{ijk}
=&\dfrac{\partial}{\partial y^k}\left( \bar{g}_{ij} \right) \\
=& 2(m+1) \left( \dfrac{F^{2m}}{\beta^{2m}}+ \dfrac{F^{m-1}}{\beta^{m-1}} \right)\ C_{ijk}\\
&+(m+1)\biggl\{ 2m \dfrac{\beta F^{2m-2} y_{_k}-F^{2m} b_k}{\beta^{2m+1}}
+(m-1) \dfrac{\beta F^{m-3} y_{_k}-F^{m-1} b_k}{\beta^{m}} \biggr\}\left(h_{ij}+\dfrac{y_i y_j}{F^2} \right)\\ 
&-(m+1)\left( 2m \dfrac{F^{2m+1}}{\beta^{2m+1}}+(m-1) \dfrac{F^{m}}{\beta^{m}} \right) \left( \dfrac{b_i h_{jk}+b_j h_{ik}}{F}\right)\\
&-m(m+1)\biggl\{ 2(2m+1) \dfrac{\beta F^{2m-1} y_{_k}-F^{2m+1} b_k}{\beta^{2m+2}}
+(m-1) \dfrac{\beta F^{m-2} y_{_k}-F^{m} b_k}{\beta^{m+1}} \biggr\} \left(\dfrac{b_i y_j+b_j y_i}{F} \right)\\
&+ (m+1)\left( 2m \dfrac{F^{2m}}{\beta^{2m}}+(m-1) \dfrac{F^{m-1}}{\beta^{m-1}} \right) \left(\dfrac{y_i h_{jk}+y_j h_{ik}}{F^2} \right)\\
&+(m+1)\biggl\{ 4m^2 \dfrac{\beta F^{2m-2} y_{_k}-F^{2m} b_k}{\beta^{2m+1}}+(m-1)^2 \dfrac{\beta F^{m-3} y_{_k}-F^{m-1} b_k}{\beta^{m}}\biggr\} \dfrac{y_i y_j}{F^2}\\
&+m(m+1)\biggl\{ 2(2m+1) \dfrac{\beta F^{2m} y_{_k}-F^{2m+2} b_k}{\beta^{2m+3}}
+(m-1) \dfrac{\beta F^{m-1} y_{_k}-F^{m+1} b_k}{\beta^{m+2}} \biggr\} b_ib_j
\end{align*}
\begin{align*}
=&2(m+1) \left( \dfrac{F^{2m}}{\beta^{2m}}+ \dfrac{F^{m-1}}{\beta^{m-1}} \right)\ C_{ijk} -(m+1) \dfrac{F^{m-1}}{\beta^{m}} \left(2m\dfrac{F^{m+1}}{\beta^{m+1}}+(m-1)\right)
\sum_{\text{cyclic sum}} h_{ij}\left(b_k-\dfrac{\beta}{F^2} y_{_k} \right)\\
&-m\dfrac{F^{m+1}}{\beta^{m+2}} \left( (2m+1)(2m+2) \dfrac{F^{m+1}}{\beta^{m+1}}+(m^2-1) \right) \prod_{\text{cyclic product}} \left(b_i-\dfrac{\beta}{F^2} y_i \right).
\end{align*}
After simplification, we get 
\begin{equation}{\label{KR2.5}}
\begin{split}
\bar{C}_{ijk}
=&(m+1) \left( \dfrac{F^{2m}}{\beta^{2m}}+ \dfrac{F^{m-1}}{\beta^{m-1}} \right) C_{ijk}\\
&-\dfrac{(m+1)}{2} \dfrac{F^{m-1}}{\beta^{m}} \left(2m\dfrac{F^{m+1}}{\beta^{m+1}}+(m-1)\right)\left(  h_{ij}m_k+h_{jk}m_i+h_{ki}m_j\right)\\
&-\dfrac{m}{2}\dfrac{F^{m+1}}{\beta^{m+2}} \left( (2m+1)(2m+2) \dfrac{F^{m+1}}{\beta^{m+1}}+(m^2-1) \right)m_im_jm_k.
\end{split}
\end{equation}
The above discussion leads to the following proposition.
\begin{prop}
	Let $(M, \bar{F})$ be an $n-$dimensional Finsler space with $\bar{F}= \dfrac{F^{m+1}}{\beta^m} + \beta \  \left( m\neq 0, -1\right) $ as generalized Kropina-Randers changed metric. Then its Cartan tensor is given by equation (\ref{KR2.5}).
\end{prop}
Next, we find fundamental metric tensor for square-Randers changed  metric
\begin{equation}{\label{KR3.1}}
\bar{F}= \dfrac{(F+\beta)^{2}}{F} + \beta.
\end{equation}
Differentiating (\ref{KR3.1}) w.r.t. ${y^{i}},$ we get
\begin{equation}{\label{KR3.2}}
\bar{F}_{y^{i}} =\left( 1-\dfrac{\beta^2}{F^2}\right)  F_{y^{i}}+ \left( \dfrac{2 \beta}{F} +3\right)  b_{i}.
\end{equation}
Differentiation of (\ref{KR3.2}) further w.r.t. ${y^{j}}$ gives
\begin{equation}{\label{KR3.3}}
 \bar{F}_{y^{i}y^{j}} =\left( 1-\dfrac{\beta^2}{F^2}\right)  F_{y^i y^j} -  \dfrac{2 \beta}{F^2} \left( b_{i}F_{y^{j}}+b_{j}F_{y^{i}}\right)  + \dfrac{2 \beta^2}{F^3} F_{y^{i}}F_{y^{j}}+ \dfrac{2}{F} b_{i} b_{j}.
\end{equation}
Now, 
\begin{align*}
\bar{g}_{ij}&=\bar{F} \bar{F}_{y^{i}y^{j}} + \bar{F}_{y^{i}} \bar{F}_{y^{j}}\\
&=\left\{ F+ \dfrac{\beta^2}{F} +3 \beta \right\} \left\{ \left( 1-\dfrac{\beta^2}{F^2}\right)  F_{y^i y^j} -  \dfrac{2 \beta}{F^2} \left( b_{i}F_{y^{j}}+b_{j}F_{y^{i}}\right)  + \dfrac{2 \beta^2}{F^3} F_{y^{i}}F_{y^{j}}+ \dfrac{2}{F} b_{i} b_{j} \right\}\\
&\ \ \ \ \ \ \ + \left\{ \left( 1-\dfrac{\beta^2}{F^2}\right)  F_{y^{i}}+ \left( \dfrac{2 \beta}{F} +3\right)  b_{i} \right\} \left\{ \left( 1-\dfrac{\beta^2}{F^2}\right)  F_{y^{j}}+ \left( \dfrac{2 \beta}{F} +3\right)  b_{j} \right\}.
\end{align*}
Simplifying, we get
\begin{equation}{\label{KR3.4}}
\begin{split}
\bar{g}_{ij}
&= \left( 1+ \dfrac{3 \beta}{F} -\dfrac{3 \beta^3}{F^3} -\dfrac{\beta^4}{F^4} \right) g_{ij} + \left(3-\dfrac{4 \beta^3}{F^3}-\dfrac{9 \beta^2}{F^2} \right)  \left( b_{i}F_{y^{j}}+b_{j}F_{y^{i}}\right)\\
&\ \  \  + \left( -\dfrac{3 \beta}{F}+\dfrac{9 \beta^3}{F^3} +\dfrac{4 \beta^4}{F^4} \right)  F_{y^{i}} F_{y^{j}}+\left( 11+\dfrac{18 \beta}{F} + \dfrac{6 \beta^2}{F^2} \right)  b_{i} b_{j}.
\end{split}
\end{equation}
Hence, we have following:
\begin{prop}
	Let $(M, \bar{F})$ be an $n-$dimensional Finsler space with $\bar{F}= \dfrac{(F+\beta)^{2}}{F} + \beta$ as a Square-Randers changed metric. Then its fundamental metric tensor is given by equation (\ref{KR3.4}).
\end{prop}
Next, we find Cartan tensor for square-Randers changed  metric.\\
By definition, we have 
$$2C_{ijk}(y)=\dfrac{\partial\mathnormal{g}_{ij}}{\partial y^k}.$$
From the equation (\ref{KR3.4}), we get \\
\begin{align*}
2 \bar{C}_{ijk}
=&\dfrac{\partial}{\partial y^k}\left( \bar{g}_{ij} \right) \\
=&2\left( 1+ \dfrac{3 \beta}{F} -\dfrac{3 \beta^3}{F^3} -\dfrac{\beta^4}{F^4} \right) C_{ijk}\\
&+\biggl\{ \dfrac{3(F^2 b_k-\beta y_k)}{F^3}-\dfrac{9(\beta^2 F^2 b_k-\beta^3 y_k)}{F^5}-\dfrac{4(\beta^3 F^2 b_k-\beta^4 y_k)}{F^6} \biggr\} \left(h_{ij}+\dfrac{y_i y_j}{F^2} \right)\\ 
&-\left(3-\dfrac{4 \beta^3}{F^3}-\dfrac{9 \beta^2}{F^2} \right) \left( \dfrac{b_i h_{jk}+b_j h_{ik}}{F}\right)\\
&-\biggl\{ -\dfrac{12(\beta^2 F^2 b_k-\beta^3 y_k)}{F^5}-\dfrac{18(\beta F^2 b_k-\beta^2 y_k)}{F^4} \biggr\} \left(\dfrac{b_i y_j+b_j y_i}{F} \right)\\
&+\left( -\dfrac{3 \beta}{F}+\dfrac{9 \beta^3}{F^3} +\dfrac{4 \beta^4}{F^4} \right) \left(\dfrac{y_i h_{jk}+y_j h_{ik}}{F^2} \right)\\
&+\biggl\{ \dfrac{16(\beta^3 F^2 b_k-\beta^4 y_k)}{F^6}+\dfrac{27(\beta^2 F^2 b_k-\beta^3 y_k)}{F^5}-\dfrac{3( F^2 b_k-\beta y_k)}{F^3} \biggr\} \dfrac{y_i y_j}{F^2}\\
&+\biggl\{ \dfrac{18( F^2 b_k-\beta y_k)}{F^3}+\dfrac{12(\beta F^2 b_k-\beta2 y_k)}{F^4} \biggr\} b_ib_j\\
\end{align*}
\begin{align*}
=&2 \left( 1+ \dfrac{3 \beta}{F} -\dfrac{3 \beta^3}{F^3} -\dfrac{\beta^4}{F^4} \right)\ C_{ijk}
+\left( \dfrac{3 }{F}-\dfrac{9 \beta^2}{F^3} -\dfrac{4 \beta^3}{F^4} \right) \sum_{\text{cyclic sum}} h_{ij}\left(b_k-\dfrac{\beta}{F^2} y_{_k} \right)\\
&+\left( \dfrac{18}{F}+\dfrac{12 \beta}{F^2} \right)  \prod_{\text{cyclic product}} \left(b_i-\dfrac{\beta}{F^2} y_i \right).
\end{align*}
After simplification, we get
\begin{equation}{\label{KR3.5}}
\begin{split}
\bar{C}_{ijk}
=& \left( 1+ \dfrac{3 \beta}{F} -\dfrac{3 \beta^3}{F^3} -\dfrac{\beta^4}{F^4} \right)\ C_{ijk}
+\dfrac{1}{2}\left( \dfrac{3 }{F}-\dfrac{9 \beta^2}{F^3} -\dfrac{4 \beta^3}{F^4} \right) \left(  h_{ij}m_k+h_{jk}m_i+h_{ki}m_j\right)\\
&\ \ +\left( \dfrac{9}{F}+\dfrac{6 \beta}{F^2} \right) m_im_jm_k.
\end{split}
\end{equation}
The above discussion leads to the following proposition.
\begin{prop}
	Let $(M, \bar{F})$ be an $n-$dimensional Finsler space with $\bar{F}= \dfrac{(F+\beta)^{2}}{F} + \beta$ as a Square-Randers changed metric. Then its Cartan tensor is given by equation (\ref{KR3.5}).
\end{prop}
Next, we find fundamental metric tensor for Matsumoto-Randers changed  metric
\begin{equation}{\label{KR4.1}}
\bar{F}= \dfrac{F^{2}}{F-\beta} + \beta.
\end{equation}
Differentiating (\ref{KR4.1}) w.r.t. ${y^{i}},$ we get
\begin{equation}{\label{KR4.2}}
\bar{F}_{y^{i}} = \dfrac{F^2-2 \beta F}{(F-\beta)^2}F_{y^{i}}+ \left(  \dfrac{F^2}{(F-\beta)^2}+1\right)  b_{i}.
\end{equation}
Differentiation of (\ref{KR4.2}) further w.r.t. ${y^{j}}$ gives
\begin{equation}{\label{KR4.3}}
 \bar{F}_{y^i y^j} =\dfrac{1}{(F-\beta)^3}\biggl\{ (F^3-3 \beta F^2+2 \beta^2 F) F_{y^i y^j}-2 \beta F \left( b_{i}F_{y^{j}}+b_{j}F_{y^{i}}\right) +2 \beta^2 F_{y^i} F_{y^j} +2 F^2  b_{i} b_{j}\biggr\}.
\end{equation}
Now, 
\begin{align*}
	\bar{g}_{ij}
	&=\bar{F} \bar{F}_{y^{i}y^{j}} + \bar{F}_{y^{i}} \bar{F}_{y^{j}}\\
	&= \dfrac{F^2+\beta F-\beta^2}{(F-\beta)^4}\biggl\{ (F^3-3 \beta F^2+2 \beta^2 F) F_{y^i y^j}-2 \beta F \left( b_{i}F_{y^{j}}+b_{j}F_{y^{i}}\right) +2 \beta^2 F_{y^i} F_{y^j} +2 F^2  b_{i} b_{j}\biggr\}\\
	&\ \ \ \ +\left\{ \dfrac{F^2-2 \beta F}{(F-\beta)^2}F_{y^{i}}+ \left(  \dfrac{F^2}{(F-\beta)^2}+1\right)  b_{i}\right\}\left\{ \dfrac{F^2-2 \beta F}{(F-\beta)^2}F_{y^{j}}+ \left(  \dfrac{F^2}{(F-\beta)^2}+1\right)  b_{j} \right\}.
\end{align*}
After simplification, we get 
\begin{equation}{\label{KR4.4}}
\begin{split}
\bar{g}_{ij}
&= \dfrac{1}{(F-\beta)^4}\biggl\{ (F^4-2\beta F^3-2\beta^2 F^2+5 \beta^3 F-2\beta^4) g_{ij}+
\left( 2F^4-8\beta F^3+3 \beta^2 F^2 \right)  \left( b_{i}F_{y^{j}}+b_{j}F_{y^{i}}\right)\\
& \ \ \ \ \ \ \ \ \ \ \ \ \ \ \   +(-2 \beta F^3+8 \beta^2 F^2-3 \beta^3 F) F_{y^i y^j} + (6F^4-6\beta F^3+6 \beta^2 F^2-4 \beta^3 F+ \beta^4)  b_{i} b_{j}\biggr\}.
\end{split}
\end{equation}
Hence, we have following
\begin{prop}
	Let $(M, \bar{F})$ be an $n-$dimensional Finsler space with $\bar{F}= \dfrac{F^{2}}{F-\beta} + \beta$ as a Matsumoto-Randers changed metric. Then its fundamental metric tensor is given by equation (\ref{KR4.4}).
\end{prop}
Next, we find Cartan tensor for Matsumoto-Randers changed  metric.\\
By definition, we have 
$$2C_{ijk}(y)=\dfrac{\partial\mathnormal{g}_{ij}}{\partial y^k}.$$
From the equation (\ref{KR4.4}), we get \\
\begin{align*}
2 \bar{C}_{ijk}
=&\dfrac{\partial}{\partial y^k}\left( \bar{g}_{ij} \right) \\
=&\dfrac{2(F^4-2\beta F^3-2\beta^2 F^2+5 \beta^3 F-2\beta^4)}{(F-\beta)^4} \ C_{ijk}\\
&+ \biggl\{ \dfrac{-\beta\left( 2F^2-8\beta F+3\beta^2\right) } {F(F-\beta)^4}y_{_K} +\dfrac{ 2F^3-8\beta F^2+3\beta^2F} {(F-\beta)^4}b_k  \biggr\} \left(h_{ij}+\dfrac{y_i y_j}{F^2} \right)\\ 
&+\left( \dfrac{ 2F^4-8\beta F^3+3\beta^2F^2} {(F-\beta)^4} \right) \left( \dfrac{b_i h_{jk}+b_j h_{ik}}{F}\right)\\
&+\biggl\{ \dfrac{ 18\beta^2 F-6\beta^3 } {(F-\beta)^5}y_{_K} +\dfrac{ -18\beta F^3+6\beta^2 F^2} {(F-\beta)^5}b_k \biggr\} \left(\dfrac{b_i y_j+b_j y_i}{F} \right)\\
&+\left( \dfrac{ 8\beta^2 F^2-2\beta F^3-3\beta^3 F} {(F-\beta)^4} \right) \left(\dfrac{y_i h_{jk}+y_j h_{ik}}{F^2} \right)\\
&+\biggl[ \biggl\{\dfrac{-6\beta F^2+16\beta^2F-3\beta^3}{(F-\beta)^4} -\dfrac{4\left( -2\beta F^2+8\beta^2F-3\beta^3\right) }{(F-\beta)^5}\biggr\}y_{_k}\\
&\ \ \ +\biggl\{\dfrac{-2F^3+16\beta F^2-9\beta^2F}{(F-\beta)^4} +\dfrac{4\left( -2\beta F^3+8\beta^2F^2-3\beta^3F\right) }{(F-\beta)^5}\biggr\}b_k\biggr] \dfrac{y_i y_j}{F^2}\\
&+\biggl\{ \dfrac{ -18\beta F^3+6\beta^2F^2 } {F(F-\beta)^5}y_{_k} +\dfrac{ 18F^4-6\beta F^3} {(F-\beta)^5}b_k \biggr\} b_ib_j
\end{align*}
\begin{align*}
=&\dfrac{2(F^4-2\beta F^3-2\beta^2 F^2+5 \beta^3 F-2\beta^4)}{(F-\beta)^4} \ C_{ijk}\\
&+\dfrac{ 2F^3-8\beta F^2+3\beta^2F} {(F-\beta)^4} \sum_{\text{cyclic sum}} h_{ij}\left(b_k-\dfrac{\beta}{F^2} y_{_k} \right)
+\dfrac{ 18F^4-6\beta F^3} {(F-\beta)^5} \prod_{\text{cyclic product}} \left(b_i-\dfrac{\beta}{F^2} y_i \right)
\end{align*}
Simplifying, we get
\begin{equation}{\label{KR4.5}}
\begin{split}
\bar{C}_{ijk}
=&\dfrac{F^4-2\beta F^3-2\beta^2 F^2+5 \beta^3 F-2\beta^4}{(F-\beta)^4} \ C_{ijk}\\
&+\dfrac{ 2F^3-8\beta F^2+3\beta^2F} {2(F-\beta)^4} \left(  h_{ij}m_k+h_{jk}m_i+h_{ki}m_j\right)
+\dfrac{ 9F^4-3\beta F^3} {(F-\beta)^5} m_im_jm_k.
\end{split}
\end{equation}
The above discussion leads to the following proposition.
\begin{prop}
	Let $(M, \bar{F})$ be an $n-$dimensional Finsler space with $\bar{F}= \dfrac{F^{2}}{F-\beta} + \beta$ as a Matsumoto-Randers changed metric. Then its Cartan tensor is given by equation (\ref{KR4.5}).
\end{prop}
Next, we find fundamental metric tensor for exponential-Randers changed  metric
\begin{equation}{\label{KR6.1}}
\bar{F}= F e^{ \beta/F} + \beta.
\end{equation}
Differentiating (\ref{KR6.1}) w.r.t. ${y^{i}},$ we get
\begin{equation}{\label{KR6.2}}
\bar{F}_{y^{i}} = \left( 1- \dfrac{\beta}{F}\right) e^{ \beta/F}F_{y^{i}} + \left( 1 +  e^{ \beta/F} \right) b_{i}.
\end{equation}
Differentiation of  (\ref{KR6.2}) further w.r.t. ${y^{j}}$ gives
\begin{equation}{\label{KR6.3}}
\bar{F}_{y^{i}y^{j}} = e^{ \beta/F} \biggl\{ \left( 1 -\dfrac{\beta}{F} \right)F_{y^{i}y^{j}} - \dfrac{\beta}{F^2}\left( b_{i}F_{y^{j}} + b_{j}F_{y^{i}}\right) + \dfrac{\beta^2}{F^3} F_{y^{i}} F_{y^{j}} + \dfrac{1}{F} b_{i} b_{j} \biggr\}.
\end{equation}
Now,
\begin{align*} 
\bar{g}_{ij}
&=\bar{F} \bar{F}_{y^{i}y^{j}} + \bar{F}_{y^{i}} \bar{F}_{y^{j}}\\
&= \left( F e^{ \beta/F} + \beta \right) e^{ \beta/F} \biggl\{ \left( 1 -\dfrac{\beta}{F} \right)F_{y^{i}y^{j}} - \dfrac{\beta}{F^2}\left( b_{i}F_{y^{j}} + b_{j}F_{y^{i}}\right) + \dfrac{\beta^2}{F^3} F_{y^{i}} F_{y^{j}} + \dfrac{1}{F} b_{i} b_{j} \biggr\} \\
& \ \ \ \ + \left\lbrace \left( 1- \dfrac{\beta}{F}\right) e^{ \beta/F}F_{y^{i}} + \left( 1 +  e^{ \beta/F} \right) b_{i}\right\rbrace \left\lbrace \left( 1- \dfrac{\beta}{F}\right) e^{ \beta/F}F_{y^{j}} + \left( 1 +  e^{ \beta/F} \right) b_{j}\right\rbrace.
\end{align*} 
Simplifying, we get
\begin{equation}{\label{KR6.4}}
\begin{split}
\bar{g}_{ij}
&= \left(e^{ \beta/F} +  \dfrac{\beta}{F}\right) e^{ \beta/F} \left( 1- \dfrac{\beta}{F}\right)  g_{ij} - e^{ \beta/F} \left\lbrace -1+ \dfrac{\beta}{F} + \dfrac{\beta^2}{F^2} + \left( -1+ \dfrac{2\beta}{F} \right) e^{ \beta/F} \right\rbrace \left( b_{i}F_{y^{j}} + b_{j}F_{y^{i}}\right) \\
& \ \ \ \ + \dfrac{\beta}{F}e^{ \beta/F} \left\lbrace-1+ \dfrac{\beta}{F} + \dfrac{\beta^2}{F^2} + \left( -1+ \dfrac{2\beta}{F}\right)  e^{ \beta/F}  \right\rbrace F_{y^{i}} F_{y^{j}} + \left\lbrace 1+ 2e^{ 2\beta/F} + \left( 2+ \dfrac{\beta}{F} \right) e^{ \beta/F}  \right\rbrace b_{i} b_{j}.
\end{split}
\end{equation}
Hence, we have following:
\begin{prop}
	Let $(M, \bar{F})$ be an $n-$dimensional Finsler space with $ \bar{F}= F e^{ \beta/F} + \beta$ as an exponential-Randers changed metric. Then its fundamental metric tensor is given by equation (\ref{KR6.4}).
\end{prop}
Next, we find Cartan tensor of exponential-Randers changed  metric.\\
By definition, we have 
$$2 C_{ijk}(y)=\dfrac{\partial\mathnormal{g}_{ij}}{\partial y^k}.$$
From the equation (\ref{KR6.4}), we get
\begin{align*}
2\bar{C}_{ijk}&=  \dfrac{\partial \bar{g}_{ij}}{\partial y^k} \\
&= 2\left(e^{ \beta/F} +  \dfrac{\beta}{F}\right) e^{ \beta/F} \left( 1- \dfrac{\beta}{F}\right) C_{ijk} \\
& \ \ \ + \dfrac{e^{ \beta/F}}{F} \left\lbrace \left( 1+ e^{ \beta/F} \right) - \dfrac{\beta}{F} \left( 1 + \dfrac{\beta}{F} + 2 e^{ \beta/F}\right)  \right\rbrace  \left(b_k-\dfrac{\beta}{F^2} y_{_k} \right)\left( h_{ij} + \dfrac{y_{i}y_{j}}{F^2} \right) \\
& \ \ \ -\dfrac{\beta e^{ \beta/F} }{F} \left\lbrace   1- \dfrac{\beta}{F} - \dfrac{\beta^2}{F^2} + \left( 1- \dfrac{2 \beta}{F}\right)e^{ \beta/F}  \right\rbrace \left( \dfrac{h_{ik} y_j + h_{jk} y_i}{F^2}\right) \\
& \ \ \  + \left\lbrace -1+ \dfrac{\beta}{F} + \dfrac{4\beta^2}{F^2} + \dfrac{\beta^3}{F^3} + \left( -1 + \dfrac{2 \beta}{F}+  \dfrac{4 \beta^2}{F^2}\right)e^{ \beta/F}   \right\rbrace \dfrac{y_i y_j e^{ \beta/F}}{F^3} \left(b_k-\dfrac{\beta}{F^2} y_{_k} \right) \\
& \ \ \ + \dfrac{e^{ \beta/F}}{F} \left\lbrace   1- \dfrac{\beta}{F} - \dfrac{\beta^2}{F^2} + \left( 1- \dfrac{2 \beta}{F}\right)e^{ \beta/F}  \right\rbrace \left( b_i h_{jk} + b_j h_{ik}\right) \\
& \ \ \ - \dfrac{\beta e^{ \beta/F}}{F^3}\left\lbrace 3 +\dfrac{\beta}{F} + 4 e^{ \beta/F} \right\rbrace \left( b_i y_j + b_j y_i\right) \left(b_k-\dfrac{\beta}{F^2} y_{_k} \right)\\
&\ \ \ + \dfrac{e^{ \beta/F}}{F} \left\lbrace 3 +\dfrac{\beta}{F} + 4 e^{ \beta/F} \right\rbrace \left(b_k-\dfrac{\beta}{F^2} y_{_k} \right) b_i b_j
\end{align*}
\begin{align*}
=& 2\left(e^{ \beta/F} +  \dfrac{\beta}{F}\right) e^{ \beta/F} \left( 1- \dfrac{\beta}{F}\right) C_{ijk} \\
&+ \dfrac{e^{ \beta/F} }{F} \left\lbrace   1- \dfrac{\beta}{F} - \dfrac{\beta^2}{F^2} + \left( 1- \dfrac{2 \beta}{F}\right)e^{ \beta/F}  \right\rbrace \sum_{\text{cyclic sum}} h_{ij}\left(b_k-\dfrac{\beta}{F^2} y_{_k} \right) \\
& + \dfrac{e^{ \beta/F}}{F} \left\lbrace 3 +\dfrac{\beta}{F} + 4 e^{ \beta/F} \right\rbrace \prod_{\text{cyclic product}} \left(b_i-\dfrac{\beta}{F^2} y_i \right).
\end{align*}
After simplifying, we get 
\begin{equation}{\label{KR6.5}}
\begin{split}
\bar{C}_{ijk}=& \left(e^{ \beta/F} +  \dfrac{\beta}{F}\right) e^{ \beta/F} \left( 1- \dfrac{\beta}{F}\right) C_{ijk}\\
&+ \dfrac{e^{ \beta/F} }{2F} \left\lbrace   1- \dfrac{\beta}{F} - \dfrac{\beta^2}{F^2} + \left( 1- \dfrac{2 \beta}{F}\right)e^{ \beta/F}  \right\rbrace \left(  h_{ij}m_k+h_{jk}m_i+h_{ki}m_j\right)\\
& + \dfrac{e^{ \beta/F}}{2F} \left\lbrace 3 +\dfrac{\beta}{F} + 4 e^{ \beta/F} \right\rbrace m_im_jm_k.
\end{split}
\end{equation}
The above discussion leads to the following proposition.
\begin{prop}
	Let $(M, \bar{F})$ be an $n-$dimensional Finsler space with $ \bar{F}= F e^{ \beta/F} + \beta$ as an exponential-Randers changed metric. Then its Cartan tensor is given by equation (\ref{KR6.5}).
\end{prop}
Next, we find fundamental metric tensor for Randers change of infinite series  metric
\begin{equation}{\label{KR7.1}}
\bar{F}= \dfrac{\beta^{2}}{\beta - F } + \beta.
\end{equation}
Differentiating (\ref{KR7.1}) w.r.t. ${y^{i}},$ we get
\begin{equation}{\label{KR7.2}}
\bar{F}_{y^{i}} = \dfrac{\left( \beta^{2} - 2 \beta F\right) b_{i} + \beta^{2} F_{y^{i}} }{\left(\beta - F \right)^{2} } + b_{i}.
\end{equation}
Differentiation of (\ref{KR7.2}) further w.r.t. ${y^{j}}$ gives
\begin{equation}{\label{KR7.3}}
 \bar{F}_{y^{i}y^{j}} = \dfrac{\beta^{2}\left(\beta- F \right)F_{y^{i}y^{j}} - 2\beta F \left( b_{i}F_{y^{j}} + b_{j}F_{y^{i}}\right) + 2 \beta^{2} F_{y^{i}} F_{y^{j}} + 2 F^{2} b_{i} b_{j}  }{\left(\beta - F \right)^{3} }.
\end{equation}
Now,
\begin{align*} 
\bar{g}_{ij}
&=\bar{F} \bar{F}_{y^{i}y^{j}} + \bar{F}_{y^{i}} \bar{F}_{y^{j}}\\
&=\left\lbrace \dfrac{\beta^{2}}{\beta - F } + \beta \right\rbrace  \left\lbrace \dfrac{ \beta^{2}\left(\beta- F \right)F_{y^{i}y^{j}} - 2\beta F \left( b_{i}F_{y^{j}} + b_{j}F_{y^{i}}\right) + 2 \beta^{2} F_{y^{i}} F_{y^{j}} + 2 F^{2} b_{i} b_{j}}{\left(\beta - F \right)^{3}}\right\rbrace \\
&\ \ \ \ +\left\lbrace \dfrac{ \left( 2\beta^{2} - 4 \beta F+ F^{2}\right) b_{i} + \beta^{2} F_{y^{i}} }{\left(\beta - F \right)^{2}} \right\rbrace   \left\lbrace \dfrac{ \left( 2\beta^{2} - 4 \beta F+ F^{2}\right) b_{j} + \beta^{2} F_{y^{j}}  }{\left(\beta - F \right)^{2}} \right\rbrace.
\end{align*} 
Simplifying, we get
\begin{equation}{\label{KR7.4}}
\begin{split}
\bar{g}_{ij}
&=\dfrac{1}{\left( \beta-F\right) ^ 4} \biggl\{ \dfrac{\beta^{3} \left(2 \beta - F \right) \left( \beta-F\right) } {F}   g_{ij} + \left( 2 \beta^{4} - 8 \beta^{3} F + 3 \beta^{2} F^{2}\right) \left( b_{i}F_{y^{j}} + b_{j}F_{y^{i}}\right) \\
& \ \ \ \ \ \ \ \ \ \ \ \ \ \ \ \ \  +  \left( -2\beta^{5}/F + 8 \beta^{4} -3 \beta^{3}F \right) F_{y^{i}} F_{y^{j}} +  \left( 4 \beta^{4} -16 \beta^{3} F +24 \beta^{2} F^{2}- 10 \beta F^{3} + F^{4}\right) b_{i} b_{j} \biggr\}
\end{split}
\end{equation}
Hence, we have following:
\begin{prop}
	Let $(M, \bar{F})$ be an $n-$dimensional Finsler space with $ \bar{F}= \dfrac{\beta^{2}}{\beta - F } + \beta$ as an infinite series-Randers changed metric. Then its fundamental metric tensor is given by equation (\ref{KR7.4}).
\end{prop}
Next, we find Cartan tensor of Randers change of infinite series  metric.\\
By definition, we have 
$$2 C_{ijk}(y)=\dfrac{\partial\mathnormal{g}_{ij}}{\partial y^k}.$$
From the equation (\ref{KR7.4}), we get 
\begin{align*}
2\bar{C}_{ijk}&=  \dfrac{\partial \bar{g}_{ij}}{\partial y^k} \\
&= \dfrac{2\beta^3 \left( 2 \beta - F\right) }{F\left( \beta - F\right)^3 }C_{ijk}
+ \dfrac{2 \beta^4-8\beta^3 F + 3 \beta^2 F^2}{F\left( \beta- F\right)^4 } \left(b_k-\dfrac{\beta}{F^2} y_{_k} \right)\left( h_{ij} + \dfrac{y_{i}y_{j}}{F^2} \right) \\
& \ \ \ + \dfrac{\beta^3 \left(-2 \beta^2+ 8\beta F- 3 F^2 \right) }{F^3 \left( \beta- F\right)^4 }\left( h_{ik} y_j + h_{jk} y_i\right) \\
&\ \ \ + \dfrac{\beta^2 \left( -2 \beta^3 + 10 \beta^2 F - 29 \beta F^2 + 9 F^3\right)}{F^3 \left( \beta- F\right)^5 }  y_i y_j \left(b_k-\dfrac{\beta}{F^2} y_{_k} \right)\\
& \ \ \ + \dfrac{2 \beta^4-8\beta^3 F + 3 \beta^2F^2}{F\left( \beta- F\right)^4}\left( b_i h_{jk} + b_j h_{ik}\right)  - \dfrac{6 \beta F \left( F- 3 \beta\right) }{\left( \beta- F\right)^5 }\left( b_i y_j + b_j y_i\right) \left(b_k-\dfrac{\beta}{F^2} y_{_k} \right)
\end{align*}
\begin{align*}
& \ \ \ + \dfrac{6 F^3 \left( F-3 \beta\right) }{\left( \beta- F\right)^5} \left(b_k-\dfrac{\beta}{F^2} y_{_k} \right) b_i b_j \\
&= \dfrac{2\beta^3 \left( 2 \beta - F\right) }{F\left( \beta - F\right)^3 }C_{ijk} + \dfrac{2 \beta^4 -8 \beta^3 F + 3\beta^2 F^2}{F\left( \beta- F\right)^4}\sum_{\text{cyclic sum}} h_{ij}\left(b_k-\dfrac{\beta}{F^2} y_{_k} \right)\\
& \ \ \ + \dfrac{6F^3\left( F-3 \beta\right)}{\left( \beta- F\right)^5}\prod_{\text{cyclic product}} \left(b_i-\dfrac{\beta}{F^2} y_i \right).
\end{align*}
After simplification, we get
\begin{equation}{\label{KR7.5}}
\bar{C}_{ijk}= \dfrac{\beta^3 \left( 2 \beta - F\right) }{F\left( \beta - F\right)^3 }C_{ijk} + \dfrac{2 \beta^4 -8 \beta^3 F + 3\beta^2 F^2}{2F\left( \beta- F\right)^4}
\left(  h_{ij}m_k+h_{jk}m_i+h_{ki}m_j\right)
+ \dfrac{3F^3\left( F-3 \beta\right)}{\left( \beta- F\right)^5} m_im_jm_k.
\end{equation}
Hence, we have the following proposition.
\begin{prop}
	Let $(M, \bar{F})$ be an $n-$dimensional Finsler space with $ \bar{F}= F e^{ \beta/F} + \beta$ as an infinite series-Randers changed metric. Then its Cartan tensor is given by equation (\ref{KR7.5}).
\end{prop}
\section{General formula for inverse $\bar{g}^{ij}$ of fundamental metric tensor $\bar{g}_{ij}$:}
Let us put
\begin{equation*}
 \dfrac{2(F^2+\beta^2)}{\beta^2}= \rho_{_{0}},\ 
 - \dfrac{4F^3}{\beta^3}= \rho_{_{1}},\ 
 \dfrac{4 F^2}{\beta^2}= \rho_{_{2}},\ 
 \left( \dfrac{3F^4}{\beta^4}+1 \right)=\rho_{_{3}} 
\end{equation*}
in the equation (\ref{KR1.4}).\\
Then the fundamental metric tensor $\bar{g}_{ij}$ for Kropina-Randers changed Finsler metric $ \bar{F}=\dfrac{F^2}{\beta}+ \beta$  takes the following form 
\begin{align*} 
\bar{g}_{ij}&= \rho_{_{0}} g_{ij} + \rho_{_{1}}\left(b_{i} \dfrac{y_{j}}{F} + b_{j} \dfrac{y_{i}}{F} \right) + \rho_{_{2}}  \dfrac{y_{i}}{F} \dfrac{y_{j}}{F} + \rho_{_{3}} b_{i} b_{j} \\
&= \rho_{_{0}} \left\lbrace g_{ij} + \dfrac{\rho_{_{1}}}{\rho_{_{0}}} \left( b_{i} + \dfrac{y_{i}}{F} \right) \left( b_{j} + \dfrac{y_{j}}{F} \right) + \left(\dfrac{\rho_{_{2}} -\rho_{_{1}}}{\rho_{_{0}}} \right)\dfrac{y_{i}}{F}\dfrac{y_{j}}{F} +  \left(\dfrac{\rho_{_{3}} -\rho_{_{1}}}{\rho_{_{0}}} \right) b_{i} b_{j}\right\rbrace, 
\end{align*} 
i.e., 
\begin{equation}{\label{INV1}}
\bar{g}_{ij}
= \rho_{_{0}} \left\lbrace g_{ij} + \lambda \left( b_{i} + \dfrac{y_{i}}{F} \right) \left( b_{j} + \dfrac{y_{j}}{F} \right) + \mu \dfrac{y_{i}}{F}\dfrac{y_{j}}{F} + \nu \  b_{i} b_{j}\right\rbrace,
\end{equation}
\begin{equation}{\label{lmn}}
\text{where\ } \lambda = \dfrac{\rho_{_{1}}}{\rho_{_{0}}}, \ \mu = \dfrac{\rho_{_{2}} -\rho_{_{1}}}{\rho_{_{0}}} \  \text{and} \  \nu = \dfrac{\rho_{_{3}} -\rho_{_{1}}}{\rho_{_{0}}}.
\end{equation}
Similarly, for all other  metrics constructed in section two, the fundamental metric tensors obtained  in equations (\ref{KR2.4}), (\ref{KR3.4}), (\ref{KR4.4}), (\ref{KR6.4}) and (\ref{KR7.4}) respectively,  can be written in the form of equation (\ref{INV1}).\\
Next, we find the inverse metric tensor $\bar{g}^{ij}$ of fundamental metric tensor $\bar{g}_{ij}$ for the metric $\bar{F}.$\\
Let
\begin{equation*}{\label{INV3}}
\begin{split}
m_{ij} &= g_{ij} +\lambda \left( b_{i} + \dfrac{y_{i}}{F} \right) \left( b_{j} + \dfrac{y_{j}}{F} \right) \\
&=  g_{ij} + \lambda c_{i} c_{j}, 
\end{split}
\end{equation*}
where $  c_{i}=   b_{i} + \dfrac{y_{i}}{F}.  $ \\
Define  
\begin{align*}  
c^2 &= g^{ij} c_{i} c_{j} \\
&= g^{ij} \left( b_{i} + \dfrac{y_{i}}{F} \right) \left( b_{j} + \dfrac{y_{j}}{F} \right) \\
&= b^2 + \dfrac{2\beta}{F} + 1.
\end{align*}
Next,
\begin{align*} 
\text{det} \  m_{ij} &= \left( 1+ \lambda c^{2} \right) \text{det} \ g_{ij} \\
&= X  \text{det} \ g_{ij},
\end{align*}
\begin{equation}{\label{valueX}}
\text{where\ } X=  1+ \lambda c^{2}  = 1+ \lambda \left(  b^2 + \dfrac{2\beta}{F} + 1 \right).
\end{equation}
Also,
$$ m^{ij}= g^{ij} - \dfrac{\lambda}{1 + \lambda c^2} c^{i} c^{j} = g^{ij} - \dfrac{\lambda}{X} A^{ij}, $$
where
\begin{align*}
A^{ij}&=  c^{i} c^{j} \\
&= g^{ih}c_{h} g^{jk} c_{k}\\
&= g^{ih}  \left( b_{h} + \dfrac{y_{h}}{F} \right)     g^{jk} \left( b_{k} + \dfrac{y_{k}}{F} \right) \\
&= \left( b_{i} + \dfrac{y_{i}}{F} \right) \left( b_{j} + \dfrac{y_{j}}{F} \right).
\end{align*}
Then equation (\ref{INV1}) becomes
\begin{equation}{\label{INV5}}
\bar{g}_{ij}
= \rho_{_{0}} \left\lbrace m_{ij} + \mu \dfrac{y_{i}}{F}\dfrac{y_{j}}{F} + \nu \  b_{i} b_{j}\right\rbrace.
\end{equation}
Further, assume that
\begin{equation}{\label{INV4}}
n_{ij}= m_{ij} + \mu d_{i} d_{j},\  \text{where}   \ d_{i}= \dfrac{y_{i}}{F}.
\end{equation}
Define  
\begin{align*}  
d^2 &= m^{ij} d_{i} d_{j} \\
&= \left\lbrace g^{ij}- \dfrac{\lambda}{X} \left( b^{i}b^{j} + b^{i}\dfrac{y^{j}}{F}+ b^{j}\dfrac{y^{i}}{F} + \dfrac{ y^{i}y^{j}}{F^2}\right) \right\rbrace  \dfrac{y^{i}}{F}\dfrac{y^{j}}{F} \\
&= 1- \dfrac{\lambda}{X}\left( 1+ \dfrac{\beta}{F}\right)^2.
\end{align*}
Next,
\begin{align*} 
\text{det} \  n_{ij} &= \left( 1+ \mu d^{2} \right) \text{det} \ m_{ij} \\
&= YX  \text{det} \ g_{ij},
\end{align*}
\begin{equation}{\label{valueY}}
\text{where\ } Y=  1+ \mu d^{2}  = 1+ \mu \left\lbrace 1- \dfrac{\lambda}{X}\left( 1+ \dfrac{\beta}{F}\right)^2 \right\rbrace. 
\end{equation}
Also,
\begin{equation}{\label{INV8}}
\begin{split}
n^{ij}=& m^{ij} - \dfrac{\mu}{1 + \mu d^2} d^{i} d^{j} \\
&= g^{ij} - \dfrac{\lambda}{X} A^{ij} - \dfrac{\mu}{Y} B^{ij}\\
&=g^{ij} - \dfrac{\lambda}{X} \left( b_{i} + \dfrac{y_{i}}{F} \right) \left( b_{j} + \dfrac{y_{j}}{F} \right) - \dfrac{\mu}{Y} B^{ij}
\end{split}
\end{equation}
where
\begin{equation}{\label{INV9}}
\begin{split}
B^{ij}&=  d^{i} d^{j} \\
&= m^{ih}d_{h} m^{jk} d_{k}\\
&= \dfrac{1}{F^2}\left\lbrace 1- \dfrac{\lambda}{X}\left( 1+ \dfrac{\beta}{F}\right) \right\rbrace ^2 y^{i} y^{j} + \dfrac{1}{F}\left\lbrace - \dfrac{\lambda}{X}\left( 1+ \dfrac{\beta}{F}\right) + \dfrac{\lambda^2}{X^2}\left( 1+ \dfrac{\beta}{F}\right)^2\right\rbrace \left( y^{i} b^{j} + y^{j} b^{i} \right) \\
& \ \ \ \ +  \dfrac{\lambda^2}{X^2}\left( 1+ \dfrac{\beta}{F}\right)^2 b^{i} b^{j}.
\end{split}
\end{equation}
From the  equations (\ref{INV5}) and (\ref{INV4}), we get
\begin{equation}{\label{INV6}}
\bar{g}_{ij}
= \rho_{_{0}} \left\lbrace n_{ij} + \nu \  b_{i} b_{j}\right\rbrace,
\end{equation}
Define
\begin{equation}{\label{valuebtilda}}
\begin{split}
\tilde{b}^2 &= n^{ij} b_{i} b_{j} \\
&= b^2 - \dfrac{\lambda}{X}\left( b^2 + \dfrac{\beta}{F}\right)^2 - \dfrac{\mu}{Y}\dfrac{\beta^2}{F^2}\left\lbrace 1- \dfrac{\lambda}{X}\left( 1+ \dfrac{\beta}{F}\right) \right\rbrace ^2 \\
& \ \ \ - 2 b^2 \dfrac{\beta}{F} \dfrac{\mu}{Y}\left\lbrace - \dfrac{\lambda}{X}\left( 1+ \dfrac{\beta}{F}\right) + \dfrac{\lambda^2}{X^2}\left( 1+ \dfrac{\beta}{F}\right)^2\right\rbrace - b^4\dfrac{\mu}{Y}\dfrac{\lambda^2}{X^2} \left( 1+ \dfrac{\beta}{F}\right)^2.
\end{split}
\end{equation}
From equation (\ref{INV6}), we have 
$$ \text{det} \  \bar{g}_{ij} = \rho_{_{0}}^n \left( 1+ \nu \tilde{b}^{2} \right) \text{det} \ n_{ij},$$
and the inverse $\bar{g}^{ij}$ of $\bar{g}_{ij}$ is given by 
\begin{align}{\label{INV7}}
\bar{g}^{ij} &= \dfrac{1}{\rho_{_{0}}} \left( n^{ij} - \dfrac{\nu}{1 + \nu \tilde{b}^2} b^{i} b^{j}\right).
\end{align}
From the equations (\ref{INV8}), (\ref{INV9}) and (\ref{INV7}), we get
\begin{equation}{\label{INV2}}
\begin{split}
\bar{g}^{ij}
&= \dfrac{1}{\rho_{_{0}}}\biggl\{ g^{ij} + \dfrac{1}{F}\left[- \dfrac{\lambda}{X} -\dfrac{\mu}{Y} \left( - \dfrac{\lambda}{X}\left( 1+ \dfrac{\beta}{F}\right) + \dfrac{\lambda^2}{X^2}\left( 1+ \dfrac{\beta}{F}\right)^2\right) \right] \left(y^i b^j+ y^j b^i \right) \\
& \ \ \ \ + \dfrac{1}{F^2} \left[- \dfrac{\lambda}{X} - \dfrac{\mu}{Y} \left(1- \dfrac{\lambda}{X}\left( 1+ \dfrac{\beta}{F}\right) \right) ^2   \right] y^{i} y^{j} + \left[- \dfrac{\lambda}{X} -\dfrac{\mu}{Y}\dfrac{\lambda^2}{X^2}\left( 1+ \dfrac{\beta}{F}\right)^2 -  \dfrac{\nu}{1 + \nu \tilde{b}^2} \right]b^{i} b^{j}  \biggr\},
\end{split}
\end{equation}
where the values of $\lambda,\  \mu$ and $\nu$ are obtained by equation (\ref{lmn}),  $X$ by equation (\ref{valueX}), $Y$ by equation (\ref{valueY}) and $\tilde{b}$ by equation (\ref{valuebtilda}).\\
The above discussion leads to the following theorem.
\begin{theorem}
Let $(M, F)$ be an $n-$dimensional Finsler space. Let $\bar{F}=f(F, \beta)$ be a Finsler metric obtained by Randers $\beta-$change of a Finsler metric $F=f(\alpha, \beta),$ where $\alpha = \sqrt{a_{ij}(x) y^i y^j} $ is a Riemannian metric and  $\beta=b_i(x) y^i$ is a $1-$form. If $\bar{g}_{ij}$ is the fundamental metric tensor of $\bar{F}$ given by (\ref{INV1}), then its inverse fundamental metric tensor, $\bar{g}^{ij}$ is given by (\ref{INV2}).
\end{theorem}
\section{Projective flatness of Finsler metrics}
The notion of projective flatness is one of the important topics in differential geometry.
The real starting point of the investigations of projectively flat metrics is Hilbert’s
fourth problem \cite{Hilbert4prob}. During the International Congress of Mathematicians,  held in Paris(1900), Hilbert asked about the spaces in which the shortest curves between any pair of points are straight lines. The first answer was
given by Hilbert’s student G. Hamel in 1903. In \cite{GHAM1903}, Hamel found the necessary
and sufficient conditions in order that a space satisfying a 
system of axioms, which is a modification of Hilbert’s system of axioms for Euclidean
geometry, be projectively flat.
After Hamel, many authors (see \cite{SAB.GS.2016PFDF}, \cite{GS.SAB.2015PF},...,\cite{GS.RKS.RDSK.2017PFDF}, \cite{RY.GS.2013PF}) have worked on this topic. Here we find necessary and sufficient conditions for Randers change of some special $(\alpha, \beta)-$metrics to be projectively flat.
For this, first we discuss some definitions and results related to projective flatness. 
\begin{definition}
Two Finsler metrics $F$ and $\bar{F}$ on a manifold $M$ are called projectively equivalent if they have same geodesics as point sets, \\
i.e., for any geodesic $\bar{\sigma}(\bar{t}) $ of $\bar{F},$ there is a geodesic $\sigma(t):=\bar{\sigma}(\bar{t}(t))$ of $F,$ where $\bar{t}=\bar{t}(t)$ is oriented re-parametrization, and vice-versa. 
\end{definition}
Next, recall (\cite{ChernShenRFG}, \cite{GHAM1903}) the following theorem:
\begin{theorem}{\label{thm5.1}}
Let $F$ and $\bar{F}$ be two Finsler metrics on a manifold $M.$ Then $F$ is projectively equivalent to $\bar{F}$ if and only if 
$$ {F}_{x^ky^\ell}y^k-{F}_{x^\ell}=0. $$ Here the spray coefficients are related by
$G^i=\bar{G}^i+Py^i,$ where $ P=\dfrac{F_{x^k}y^k}{2F}.$
\end{theorem}
\begin{definition}
Let $\bar{F}$ be a standard Euclidean norm on $\mathbb{R}^n.$ Then spray coefficients vanish, i.e., $$ \bar{G}^i=0.$$
\end{definition}
Since in $\mathbb{R}^n,$ geodesics are straight lines, we have
\begin{definition}
For a Finsler metric $F$ on an open subset $\mathcal{U}$ of $\mathbb{R}^n,$ the geodesics of $F$ are straight lines if and only if the spray coefficients satisfy 
$$G^i=Py^i. $$
\end{definition}
\begin{definition}{\label{def5.4}}\cite{ChernShenRFG}
A Finsler metric $F$ on an open subset $\mathcal{U}$ of $\mathbb{R}^n$ is called projectively flat if and only if all the geodesics are straight in $\mathcal{U},$ and a Finsler metric $F$ on a manifold $M$ is called locally projectively flat, if at any point, there is a local co-ordinate system $ (x^i) $ in which $F$ is projectively flat.
\end{definition}
Therefore, by theorem (\ref{thm5.1}) and definition (\ref{def5.4}), we have
\begin{theorem}
A Finsler metric $F$ on an open subset $\mathcal{U}$ of $\mathbb{R}^n$ is projectively flat if and only of it satisfies the following system of differential equations
$$ {F}_{x^ky^\ell}y^k-{F}_{x^\ell}=0. $$
Here $G^i=Py^i,$ where $G^i$ is spray coefficient and the scalar function $ P=\dfrac{F_{x^k}y^k}{2F} $ is called projective factor of $F.$
\end{theorem}
In further calculations, we assume that $F^2=A,$ i.e., $F$ is a square root Finsler metric.\\
First we find necessary and sufficient conditions for Kropina-Randers changed Finsler metric
\begin{equation*}
\bar{F}= \dfrac{F^{2}}{\beta} + \beta
\end{equation*}
to be projectively flat.\\
Let us put $F^2=A$ in $\bar{F},$ then 
\begin{equation}{\label{51.1}}
\bar{F}= \dfrac{A}{\beta} + \beta.
\end{equation}
Differentiating (\ref{51.1}) w.r.t. $x^k,$ we get
\begin{equation}{\label{51.2}}
\bar{F}_{x^k}=\dfrac{A_{x^k}}{\beta}-\dfrac{A}{\beta^2} \beta_{x^k}+\beta_{x^k}.
\end{equation}
Differentiation of  (\ref{51.2}) further  w.r.t. $y^\ell$ gives
\begin{equation}{\label{51.3}}
\bar{F}_{x^ky^\ell}=\dfrac{A_{x^ky^\ell}}{\beta}-\dfrac{A_{x^k}}{\beta^2} \beta_\ell -\dfrac{A}{\beta^2} \beta_{x^k y^\ell} -\dfrac{A_\ell}{\beta^2} \beta_{x^k}+
\dfrac{2A}{\beta^3} \beta_{x^k}\beta_{\ell}+\beta_{x^k y^\ell}.
\end{equation}
Contracting (\ref{51.3}) with $y^k$, we get
\begin{align*}
\bar{F}_{x^ky^\ell}y^k
&=\dfrac{A_{0 \ell}}{\beta}-\dfrac{A_0 \beta_\ell}{\beta^2}  -\dfrac{A \beta_{0 \ell}}{\beta^2}  -\dfrac{A_\ell \beta_{0}}{\beta^2} +
\dfrac{2A  \beta_{0}\beta_{\ell}}{\beta^3}+\beta_{0 \ell}\\
&=\dfrac{1}{\beta^3}\left\{\beta^2 A_{0 \ell}-\beta A_0 \beta_\ell  -\beta A \beta_{0 \ell} -\beta A_\ell \beta_{0}+ 2A \beta_{0}\beta_{\ell}+\beta^3 \beta_{0 \ell} \right\}.
\end{align*}
From the equation (\ref{51.2}), we have
$$ \bar{F}_{x^\ell}=\dfrac{A_{x^\ell}}{\beta}-\dfrac{A}{\beta^2} \beta_{x^\ell}+\beta_{x^\ell}
=\dfrac{1}{\beta^3} \left\{ \beta^2 {A_{x^\ell}}-{A}\beta \beta_{x^\ell}+\beta^3 \beta_{x^\ell} \right\}.$$
We know that $ \bar{F}$ is projectively flat if and only if $ \bar{F}_{x^ky^\ell}y^k-\bar{F}_{x^\ell}=0, $\\
i.e., 
 $$ \dfrac{1}{\beta^3} \left\{\beta^2 A_{0 \ell}-\beta A_0 \beta_\ell  -\beta A \beta_{0 \ell} -\beta A_\ell \beta_{0}+ 2A \beta_{0}\beta_{\ell}+\beta^3 \beta_{0 \ell} \right\}-\dfrac{1}{\beta^3} \left\{ \beta^2 {A_{x^\ell}}-{A}\beta \beta_{x^\ell}+\beta^3 \beta_{x^\ell} \right\}=0, $$
i.e.,
 $$A\left\{2\beta_{0}\beta_{\ell} +\beta\left(  \beta_{x^\ell}-\beta_{0 \ell}\right)\right\} +\left\{ \beta^3\left(  \beta_{0 \ell} - \beta_{x^\ell}\right) +  \beta^2\left(  A_{0 \ell}-A_{x^\ell}\right)  -\beta\left( A_0 \beta_\ell + A_\ell \beta_{0}\right)\right\}  =0. $$
 From the above equation, we conclude that $ \bar{F}$ is projectively flat if and only if following two equations are satisfied.
 \begin{equation}{\label{PF1.1}}
 2\beta_{0}\beta_{\ell} +\beta\left(  \beta_{x^\ell}-\beta_{0 \ell}\right)=0
 \end{equation}
 \begin{equation}{\label{PF1.2}}
 \beta^3\left(  \beta_{0 \ell} - \beta_{x^\ell}\right) +  \beta^2\left(  A_{0 \ell}-A_{x^\ell}\right)  -\beta\left( A_0 \beta_\ell + A_\ell \beta_{0}\right)=0
 \end{equation}
 Above discussion  leads to the following theorem.
\begin{theorem}
Let $(M, \bar{F})$ be an $n-$dimensional Finsler space with $\bar{F}= \dfrac{F^{2}}{\beta} + \beta$ as a Kropina-Randers changed metric. Then $\bar{F}$ is projectively flat if and only if equations (\ref{PF1.1}) and (\ref{PF1.2}) are satisfied.
 \end{theorem}
Next, we find necessary and sufficient conditions for generalized Kropina-Randers changed  Finsler metric 
$$ \bar{F}= \dfrac{F^{m+1}}{\beta^m} + \beta \ \left(  m\neq 0, -1\right) $$
to be projectively flat.\\ 
Let us put $ F^2=A$ in $\bar{F},$ then
\begin{equation}{\label{52.1}}
\bar{F}= \dfrac{A^{(m+1)/2}}{\beta^m} + \beta. 
\end{equation}
Differentiating (\ref{52.1}) w.r.t. $x^k,$ we get
\begin{equation}{\label{52.2}}
\bar{F}_{x^k}=\dfrac{m+1}{2\beta^m} A^{(m-1)/2} A_{x^k}-\dfrac{m}{\beta^{m+1}} A^{(m+1)/2} \beta_{x^k}+ \beta_{x^k}. 
\end{equation}
Differentiation of (\ref{52.2})  further w.r.t. $y^\ell$ gives
\begin{equation}{\label{52.3}}
\begin{split}
\bar{F}_{x^k y^\ell}=&\dfrac{m+1}{2\beta^m} A^{(m-1)/2} A_{x^k y^\ell} +\dfrac{m^2-1}{4\beta^m} A^{(m-3)/2} A_{x^k} A_{\ell}-\dfrac{m(m+1)}{2\beta^{m+1}} A^{(m-1)/2} A_{x^k} \beta_{\ell}\\
& -\dfrac{m}{\beta^{m+1}} A^{(m+1)/2} \beta_{x^k y^\ell} -\dfrac{m(m+1)}{2 \beta^{m+1}} A^{(m-1)/2} \beta_{x^k}A_\ell+\dfrac{m(m+1)}{\beta^{m+2}}\beta_\ell A^{(m+1)/2} \beta_{x^k}+\beta_{x^k y^\ell}.
\end{split}
\end{equation}
Contracting (\ref{52.3}) with $y^k$, we get
\begin{equation*}
\begin{split}
 \bar{F}_{x^k y^\ell}y^k&=\dfrac{m+1}{2\beta^m} A^{(m-1)/2} A_{0 \ell} +\dfrac{m^2-1}{4\beta^m} A^{(m-3)/2} A_{0} A_{\ell}-\dfrac{m(m+1)}{2\beta^{m+1}} A^{(m-1)/2} A_{0} \beta_{\ell}\\& -\dfrac{m}{\beta^{m+1}} A^{(m+1)/2} \beta_{0 \ell} -\dfrac{m(m+1)}{2 \beta^{m+1}} A^{(m-1)/2} \beta_{0}A_\ell+\dfrac{m(m+1)}{\beta^{m+2}}\beta_\ell A^{(m+1)/2} \beta_{0}+\beta_{0 \ell}\\
&=\dfrac{m+1}{2\beta^{m+2}}A^{(m-3)/2}\bigg\{ \beta^2 A A_{0 \ell} +\dfrac{m-1}{2}\beta^2 A_{0} A_{\ell}-m\beta A A_{0} \beta_{\ell} -\dfrac{2m}{m+1}\beta A^{2} \beta_{0 \ell}\\
&\ \ \ \ \ \ \ \ \ \ \  \ \ \ \ \ \ \ \ \ \ \ \ \ \ \ \ -m \beta A \beta_{0}A_\ell+2m \beta_\ell A^{2} \beta_{0}+\dfrac{2}{m+1}\beta^{m+2}A^{(-m+3)/2} \beta_{0 \ell} \bigg\}. 
\end{split}
\end{equation*}
From the equation (\ref{52.2}), we get
\begin{equation*}
\begin{split}
\bar{F}_{x^\ell}&=\dfrac{m+1}{2\beta^m} A^{(m-1)/2} A_{x^\ell}-\dfrac{m}{\beta^{m+1}} A^{(m+1)/2} \beta_{x^\ell}+ \beta_{x^\ell}\\
&=\dfrac{m+1}{2\beta^{m+2}}A^{(m-3)/2}\bigg\{ \beta^2 AA_{x^\ell}-\dfrac{2m}{m+1}\beta A^{2} \beta_{x^\ell}+\dfrac{2}{m+1}\beta^{m+2}A^{(-m+3)/2} \beta_{x^\ell} \bigg\}.
\end{split}
\end{equation*}
We know that $ \bar{F}$ is projectively flat if and only if $ \bar{F}_{x^ky^\ell}y^k-\bar{F}_{x^\ell}=0, $\\
i.e., 
\begin{align*}
\dfrac{m+1}{2\beta^{m+2}}&A^{(m-3)/2}\biggl[ \biggl\{ \beta^2 A A_{0 \ell} +\dfrac{m-1}{2}\beta^2 A_{0} A_{\ell}-m\beta A A_{0} \beta_{\ell} -\dfrac{2m}{m+1}\beta A^{2} \beta_{0 \ell}-m\beta A \beta_{0}A_\ell+2m \beta_\ell A^{2} \beta_{0}\\  +&\dfrac{2}{m+1}\beta^{m+2}A^{(-m+3)/2} \beta_{0 \ell} \biggr\}-\biggl\{ \beta^2 AA_{x^\ell}-\dfrac{2m}{m+1}\beta A^{2} \beta_{x^\ell}+\dfrac{2}{m+1}\beta^{m+2}A^{(-m+3)/2} \beta_{x^\ell} \biggr\}\biggr]=0,
\end{align*}
which implies
\begin{align*}
&\dfrac{m-1}{2}\beta^2 A_{0} A_{\ell}
+A\left\{\beta^2\left(  A_{0 \ell}-A_{x^\ell}\right)  -m\beta\left(  A_{0} \beta_{\ell} +\beta_{0}A_\ell\right) \right\} +2mA^2\left\{ \dfrac{\beta}{m+1} \left( \beta_{x^\ell}-\beta_{0 \ell}\right) +\beta_\ell \beta_{0}\right\}\\  
&+\dfrac{2}{m+1}\beta^{m+2}A^{(-m+3)/2} \left( \beta_{0 \ell}-\beta_{x^\ell}\right)=0.
\end{align*}
From the above equation, we conclude that $ \bar{F}$ is projectively flat if and only if following four equations are satisfied.
\begin{equation}{\label{PF2.1}}
\dfrac{m-1}{2}\beta^2 A_{0} A_{\ell}=0
\end{equation}
\begin{equation}{\label{PF2.2}}
\beta^2\left(  A_{0 \ell}-A_{x^\ell}\right)  -m\beta\left(  A_{0} \beta_{\ell} +\beta_{0}A_\ell\right)=0
\end{equation}
\begin{equation}{\label{PF2.3}}
\dfrac{\beta}{m+1} \left( \beta_{x^\ell}-\beta_{0 \ell}\right) +\beta_\ell \beta_{0}=0
\end{equation}
\begin{equation}{\label{PF2.4}}
\beta_{0 \ell}=\beta_{x^\ell}.
\end{equation}
Further, from the equation (\ref{PF2.1}), we see that\\
either  $m=1$ or $A_{0} A_{\ell}=0.$\\
Now, if $m=1,$ then\\ 
(\ref{PF2.3}) reduces to (\ref{PF1.1}) and (\ref{PF2.2}), (\ref{PF2.4}) reduce to (\ref{PF1.2}).\\
But the equations (\ref{PF1.1}) and (\ref{PF1.2}) are necessary and sufficient conditions for Kropina-Randers changed  Finsler metric to be projectively flat. Therefore, we exclude the case $m=1$ for general case.\\
Then
\begin{equation}{\label{PF2.5}}
A_{0} A_{\ell}=0
\end{equation}
Also from the equations (\ref{PF2.3}) and (\ref{PF2.4}), we get
 \begin{equation}{\label{PF2.6}}
 \beta_{0} \beta_{\ell}=0
 \end{equation}
Above discussion  leads to the following theorem.
\begin{theorem}
Let $(M, \bar{F})$ be an $n-$dimensional Finsler space with $ \bar{F}= \dfrac{F^{m+1}}{\beta^m} + \beta \  \left( m\neq  -1, 0, 1\right) $ as generalized Kropina-Randers changed metric. Then $\bar{F}$ is projectively flat if and only if equations (\ref{PF2.2}), (\ref{PF2.4}), (\ref{PF2.5}) and   (\ref{PF2.6}) are satisfied.
\end{theorem}
Next, we find necessary and sufficient conditions for square-Randers changed Finsler metric 
$$ \bar{F}= \dfrac{(F+\beta)^{2}}{F} + \beta$$ to be projectively flat.\\
Let us put $ F^2=A$ in $\bar{F},$ then
\begin{align}{\label{53.1}}
\bar{F}=A^{1/2}+\dfrac{\beta^2}{A^{1/2}}+3\beta.
\end{align}
Differentiating (\ref{53.1}) w.r.t. $x^k,$ we get
\begin{equation}{\label{53.2}}
\bar{F}_{x^k}=\dfrac{1}{2}A^{-1/2}A_{x^k}+2\beta A^{-1/2}\beta_{x^k}-\dfrac{1}{2}\beta^2 A^{-3/2}A_{x^k}+3\beta_{x^k}.
\end{equation}
Differentiation of (\ref{53.2}) w.r.t. $y^\ell$ gives
\begin{equation}{\label{53.3}}
\begin{split}
\bar{F}_{x^ky^\ell}
&=\dfrac{1}{2}A^{-1/2}A_{x^ky^\ell}-\dfrac{1}{4}A^{-3/2}A_{\ell}A_{x^k}+2\beta A^{-1/2}\beta_{x^ky^\ell}+2\beta_\ell A^{-1/2}\beta_{x^k}-\beta A^{-3/2}A_\ell \beta_{x^k}\\
&-\dfrac{1}{2}\beta^2 A^{-3/2}A_{x^ky^\ell}-\beta \beta_\ell A^{-3/2}A_{x^k} +\dfrac{3}{4}\beta^2 A^{-5/2}A_\ell A_{x^k} +3\beta_{x^ky^\ell}.
\end{split}
\end{equation}
Contracting (\ref{53.3}) with $y^k$, we get
\begin{align*}
\bar{F}_{x^k y^\ell}y^k
&=\dfrac{1}{2}A^{-1/2}A_{0 \ell} -\dfrac{1}{4}A^{-3/2}A_{\ell}A_{0}+2\beta A^{-1/2}\beta_{0 \ell} +2\beta_\ell A^{-1/2}\beta_{0}-\beta A^{-3/2}A_\ell \beta_{0}\\ 
&\ \ -\dfrac{1}{2}\beta^2 A^{-3/2}A_{0 \ell}-\beta \beta_\ell A^{-3/2}A_{0} +\dfrac{3}{4}\beta^2 A^{-5/2}A_\ell A_{0} +3\beta_{0 \ell}\\
&=\dfrac{1}{4}A^{-5/2}\biggl[ 2A^{2}A_{0 \ell} -A A_{\ell}A_{0}+8\beta A^{2} \beta_{0 \ell} +8\beta_\ell A^{2} \beta_{0}-4\beta A A_\ell \beta_{0}\\ 
&\ \ \ \ - 2 \beta^2 A A_{0 \ell}-4\beta \beta_\ell A A_{0} +3\beta^2 A_\ell A_{0} +12 A^{5/2} \beta_{0 \ell} \biggr].
\end{align*}
From the equation (\ref{53.2}), we get
\begin{align*}
\bar{F}_{x^\ell}
&=\dfrac{1}{2}A^{-1/2}A_{x^\ell}+2\beta A^{-1/2}\beta_{x^\ell}-\dfrac{1}{2}\beta^2 A^{-3/2}A_{x^\ell}+3\beta_{x^\ell}\\
&=\dfrac{1}{4}A^{-5/2}\left[ 2 A^{2}A_{x^\ell}+8\beta A^{2}\beta_{x^\ell}-2\beta^2 A A_{x^\ell}+12 A^{5/2}\beta_{x^\ell}\right].
\end{align*}
We know that $ \bar{F}$ is projectively flat if and only if $ \bar{F}_{x^ky^\ell}y^k-\bar{F}_{x^\ell}=0, $\\
i.e., 
\begin{align*}
&\dfrac{1}{4}A^{-5/2}\biggl[ 2A^{2}A_{0 \ell} -A A_{\ell}A_{0}+8\beta A^{2} \beta_{0 \ell} +8\beta_\ell A^{2} \beta_{0}-4\beta A A_\ell \beta_{0}
- 2 \beta^2 A A_{0 \ell}-4\beta \beta_\ell A A_{0}\\& +3\beta^2 A_\ell A_{0} +12 A^{5/2} \beta_{0 \ell} \biggr]-\dfrac{1}{4}A^{-5/2}\left[ 2 A^{2}A_{x^\ell}+8\beta A^{2}\beta_{x^\ell}-2\beta^2 A A_{x^\ell}+12 A^{5/2}\beta_{x^\ell}\right] =0,
\end{align*}
i.e., 
\begin{align*}
&3\beta^2 A_\ell A_{0}
+2A\left\{\beta^2 \left( A_{x^\ell}- A_{0 \ell}\right)-2\beta \left( A_\ell \beta_{0}+\beta_\ell  A_{0}\right)-A_{\ell}A_{0} \right\}\\
&+2A^{2}\left\{ A_{0 \ell}- A_{x^\ell}+4\beta_\ell \beta_{0}+4\beta \left( \beta_{0 \ell}- \beta_{x^\ell}\right)   \right\}
 +12 A^{5/2}\left\{ \beta_{0 \ell} -\beta_{x^\ell} \right\} =0.
\end{align*}
From the above equation, we conclude that $ \bar{F}$ is projectively flat if and only if following four equations are satisfied.
\begin{equation}{\label{PF3.1}}
A_\ell A_{0}=0
\end{equation}
\begin{equation}{\label{PF3.2}}
\beta^2 \left( A_{x^\ell}- A_{0 \ell}\right)-2\beta \left( A_\ell \beta_{0}+\beta_\ell  A_{0}\right)-A_{\ell}A_{0}=0
\end{equation}
\begin{equation}{\label{PF3.3}}
4\beta \left( \beta_{0 \ell}- \beta_{x^\ell}\right)+A_{0 \ell}- A_{x^\ell}+4\beta_\ell \beta_{0}=0
\end{equation}
\begin{equation}{\label{PF3.4}}
\beta_{0 \ell} =\beta_{x^\ell}.
\end{equation}
Above discussion  leads to the following theorem.
\begin{theorem}
Let $(M, \bar{F})$ be an $n-$dimensional Finsler space with $ \bar{F}= \dfrac{(F+\beta)^{2}}{F} + \beta$ as a square-Randers changed metric. Then $\bar{F}$ is projectively flat if and only if the equations (\ref{PF3.1}), (\ref{PF3.2}), (\ref{PF3.3}), and (\ref{PF3.4}) are satisfied.
\end{theorem}
Next, we find necessary and sufficient conditions for Matsumoto-Randers changed Finsler metric 
$$ \bar{F}= \dfrac{F^{2}}{F-\beta} + \beta$$ 
to be projectively flat.\\
Let us put $ F^2=A$ in $\bar{F},$ then
\begin{align}{\label{54.1}}
\bar{F}= \dfrac{A}{\sqrt{A}-\beta} + \beta.
\end{align}
Differentiating (\ref{54.1}) w.r.t. $x^k,$ we get
\begin{equation}{\label{54.2}}
 \bar{F}_{x^k}=\dfrac{1}{\sqrt{A}-\beta} A_{x^k}-\dfrac{\sqrt{A}}{2(\sqrt{A}-\beta)^2} A_{x^k} +\dfrac{A}{(\sqrt{A}-\beta)^2} \beta_{x^k}+\beta_{x^k}.
\end{equation}
Differentiation of  (\ref{54.2}) further w.r.t. $y^\ell$ gives
\begin{equation}{\label{54.3}}
\begin{split}
\bar{F}_{x^k y^\ell}
&=\dfrac{1}{\sqrt{A}-\beta} A_{x^k y^\ell}-\dfrac{1}{(\sqrt{A}-\beta)^2} A_{x^k}\left(\dfrac{1}{2\sqrt{A}}A_\ell-\beta_\ell \right) 
-\dfrac{\sqrt{A}}{2(\sqrt{A}-\beta)^2} A_{x^k y^\ell}\\&-\dfrac{1}{4\sqrt{A}(\sqrt{A}-\beta)^2} A_\ell A_{x^k}+\dfrac{\sqrt{A}}{(\sqrt{A}-\beta)^3} \left(\dfrac{1}{2\sqrt{A}}A_\ell-\beta_\ell \right)  A_{x^k} +\dfrac{A}{(\sqrt{A}-\beta)^2} \beta_{x^k y^\ell}\\
&+\dfrac{1}{(\sqrt{A}-\beta)^2} A_\ell \beta_{x^k}-\dfrac{2A}{(\sqrt{A}-\beta)^3} \left(\dfrac{1}{2\sqrt{A}}A_\ell-\beta_\ell \right) \beta_{x^k}
+\beta_{x^k y^\ell}.
\end{split}
\end{equation}
Contracting (\ref{54.3}) with $y^k$, we get
\begin{align*}
\bar{F}_{x^k y^\ell}y^k
&=\dfrac{1}{\sqrt{A}-\beta} A_{0 \ell}-\dfrac{1}{(\sqrt{A}-\beta)^2} A_{0}\left(\dfrac{1}{2\sqrt{A}}A_\ell-\beta_\ell \right) 
-\dfrac{\sqrt{A}}{2(\sqrt{A}-\beta)^2} A_{0 \ell}\\
&\ \ -\dfrac{1}{4\sqrt{A}(\sqrt{A}-\beta)^2} A_\ell A_{0}+\dfrac{\sqrt{A}}{(\sqrt{A}-\beta)^3} \left(\dfrac{1}{2\sqrt{A}}A_\ell-\beta_\ell \right)  A_{0} +\dfrac{A}{(\sqrt{A}-\beta)^2} \beta_{0 \ell}\\
&\ \ +\dfrac{1}{(\sqrt{A}-\beta)^2} A_\ell \beta_{0}-\dfrac{2A}{(\sqrt{A}-\beta)^3} \left(\dfrac{1}{2\sqrt{A}}A_\ell-\beta_\ell \right) \beta_{0}
+\beta_{0 \ell}.
\end{align*}
Simplifying, we get
\begin{align*}
\bar{F}_{x^k y^\ell}y^k
&=\dfrac{1}{4\sqrt{A}(\sqrt{A}-\beta)^3}
\bigg\{8A^2 \beta_{0 \ell}+2A^{3/2}\left( A_{0\ell}+4\beta_0 \beta_\ell-8 \beta \beta_{0 \ell}  \right)+6A\left( 2\beta^2\beta_{0 \ell} -\beta A_{0 \ell}\right)\\
&\ \ +\sqrt{A}\left(-4\beta^3 \beta_{0 \ell}+4\beta^2 A_{0 \ell}-4\beta\left( A_0\beta_{\ell}+A_\ell \beta_0 \right)-A_0 A_\ell \right)+3\beta A_0 A_\ell\bigg\}.
\end{align*}
From the equation (\ref{54.2}), we get
\begin{align*}
\bar{F}_{x^\ell}
=&\dfrac{1}{\sqrt{A}-\beta} A_{x^\ell}-\dfrac{\sqrt{A}}{2(\sqrt{A}-\beta)^2} A_{x^\ell} +\dfrac{A}{(\sqrt{A}-\beta)^2} \beta_{x^\ell}+\beta_{x^\ell}\\
=&\dfrac{8A^2\beta_{x^\ell}+2A^{3/2}\left(A_{x^\ell}-8\beta\beta_{x^\ell}\right)-6A\left(\beta A_{x^\ell}-2\beta^2 \beta_{x^\ell}\right)+\sqrt{A}\left(4\beta^2A_{x^\ell}-4\beta^3 \beta_{x^\ell}\right)}{4\sqrt{A}(\sqrt{A}-\beta)^3}.
\end{align*}
We know that $ \bar{F}$ is projectively flat if and only if $ \bar{F}_{x^ky^\ell}y^k-\bar{F}_{x^\ell}=0, $\\
i.e., 
\begin{align*}
&\dfrac{1}{4\sqrt{A}(\sqrt{A}-\beta)^3}\bigg[
\bigg\{8A^2 \beta_{0 \ell}+2A^{3/2}\left( A_{0\ell}+4\beta_0 \beta_\ell-8 \beta \beta_{0 \ell}  \right)+6A\left( 2\beta^2\beta_{0 \ell} -\beta A_{0 \ell}\right)\\
&\ \ +\sqrt{A}\left(-4\beta^3 \beta_{0 \ell}+4\beta^2 A_{0 \ell}-4\beta\left( A_0\beta_{\ell}+A_\ell \beta_0 \right)-A_0 A_\ell \right)+3\beta A_0 A_\ell\bigg\}\\
&-\bigg\{8A^2\beta_{x^\ell}+2A^{3/2}\left(A_{x^\ell}-8\beta\beta_{x^\ell}\right)-6A\left(\beta A_{x^\ell}-2\beta^2 \beta_{x^\ell}\right)+\sqrt{A}\left(4\beta^2A_{x^\ell}-4\beta^3 \beta_{x^\ell}\right)\bigg\}\bigg]=0,
\end{align*}
i.e.,
\begin{align*}
&8A^2\left(\beta_{0 \ell}-\beta_{x^\ell}\right)+2A^{3/2}\bigg\{ A_{0\ell}-A_{x^\ell}+4\beta_0 \beta_\ell+8\beta\left(\beta_{x^\ell}-\beta_{0 \ell}\right) \bigg\}
+6A\bigg\{2\beta^2\left(\beta_{0 \ell}-\beta_{x^\ell}\right)+\beta\left(A_{x^\ell}-A_{0 \ell}\right) \bigg\}\\&
+\sqrt{A}\bigg\{4\beta^3\left(\beta_{x^\ell}-\beta_{0 \ell}\right) +4\beta^2\left(A_{0 \ell}-A_{x^\ell}\right)-4\beta\left( A_0\beta_{\ell}+A_\ell \beta_0 \right)-A_0 A_\ell\bigg\}+3\beta A_0 A_\ell=0.
\end{align*}
From the above equation, we conclude that $ \bar{F}$ is projectively flat if and only if following five equations are satisfied.
\begin{equation}{\label{PF4.1}}
\beta_{0 \ell}=\beta_{x^\ell}
\end{equation}
\begin{equation}{\label{PF4.2}}
 A_{0\ell}-A_{x^\ell}+4\beta_0 \beta_\ell+8\beta\left(\beta_{x^\ell}-\beta_{0 \ell}\right) =0
\end{equation}
\begin{equation}{\label{PF4.3}}
2\beta^2\left(\beta_{0 \ell}-\beta_{x^\ell}\right)+\beta\left(A_{x^\ell}-A_{0 \ell}\right)=0
\end{equation}
\begin{equation}{\label{PF4.4}}
4\beta^3\left(\beta_{x^\ell}-\beta_{0 \ell}\right) +4\beta^2\left(A_{0 \ell}-A_{x^\ell}\right)-4\beta\left( A_0\beta_{\ell}+A_\ell \beta_0 \right)-A_0 A_\ell=0
\end{equation}
\begin{equation}{\label{PF4.5}}
A_0 A_\ell=0
\end{equation}
Further, from the equations (\ref{PF4.1}) and (\ref{PF4.3}), we get 
\begin{equation}{\label{PF4.6}}
A_{x^\ell}=A_{0 \ell}.
\end{equation}
Again from the equations (\ref{PF4.1}), (\ref{PF4.6}) and (\ref{PF4.2}), we get 
\begin{equation}{\label{PF4.7}}
\beta_0 \beta_\ell=0,
\end{equation}
and from the equations (\ref{PF4.1}), (\ref{PF4.6}) , (\ref{PF4.5}) and (\ref{PF4.4}), we get 
\begin{equation}{\label{PF4.8}}
A_0\beta_{\ell}+A_\ell \beta_0=0.
\end{equation}
Above discussion  leads to the following theorem.
\begin{theorem}
Let $(M, \bar{F})$ be an $n-$dimensional Finsler space with $ \bar{F}= \dfrac{F^{2}}{F-\beta} + \beta$ as a Matsumoto-Randers changed metric. Then $ \bar{F}$  is projectively flat if and only if the following equations are satisfied:
	$$ A_{0}A_\ell=0,\ A_{0 \ell}=A_{x^\ell},\ \beta_{0}\beta_\ell=0,\ \beta_{0 \ell}=\beta_{x^\ell},\ A_\ell\beta_{0}+A_{0}\beta_\ell=0.  $$
\end{theorem}
Next, we find necessary and sufficient conditions for exponential-Randers changed Finsler metric
\begin{align*}
\bar{F}
= F e^{ \beta/F} + \beta
\end{align*} 
to be projectively flat.\\
Let us put $ F^2=A$ in $\bar{F},$ then
\begin{equation}{\label{55.1}}
\bar{F}=\sqrt{A}e^{\beta/\sqrt{A}}+\beta.
\end{equation}
Differentiating (\ref{55.1}) w.r.t.$x^k,$ we get
\begin{equation}{\label{55.2}}
\bar{F}_{x^k}
=e^{\beta/\sqrt{A}}\beta_{x^k}-\dfrac{\beta}{2A}e^{\beta/\sqrt{A}}A_{x^k}
+\dfrac{1}{2\sqrt{A}}e^{\beta/\sqrt{A}}A_{x^k} +\beta_{x^k}.
\end{equation}
Differentiation of  (\ref{55.2}) further w.r.t. $y^\ell$ gives
\begin{equation}{\label{55.3}}
\begin{split}
\bar{F}_{x^k y^\ell}
=&e^{\beta/\sqrt{A}}\beta_{x^ky^\ell}+e^{\beta/\sqrt{A}}\beta_{x^k}\left(\dfrac{1}{\sqrt{A}}\beta_\ell -\dfrac{\beta}{2A^{3/2}}A_\ell \right)
-\dfrac{\beta}{2A}e^{\beta/\sqrt{A}}A_{x^ky^\ell}-\dfrac{1}{2A}e^{\beta/\sqrt{A}}A_{x^k}\beta_\ell\\
&+\dfrac{\beta}{2A^2}e^{\beta/\sqrt{A}}A_{x^k}A_\ell -\dfrac{\beta}{2A}e^{\beta/\sqrt{A}}A_{x^k}\left(\dfrac{1}{\sqrt{A}}\beta_\ell -\dfrac{\beta}{2A^{3/2}}A_\ell \right)
+\dfrac{1}{2\sqrt{A}}e^{\beta/\sqrt{A}}A_{x^ky^\ell}\\
&-\dfrac{1}{4A^{3/2}}e^{\beta/\sqrt{A}}A_{x^k}A_\ell +\dfrac{1}{2\sqrt{A}}e^{\beta/\sqrt{A}}A_{x^k}\left(\dfrac{1}{\sqrt{A}}\beta_\ell -\dfrac{\beta}{2A^{3/2}}A_\ell \right)
+\beta_{x^ky^\ell}.
\end{split}
\end{equation}
Contracting (\ref{55.3}) with $y^k$, we get
\begin{align*}
L_{x^k y^\ell}y^k
=&e^{\beta/\sqrt{A}}\beta_{0\ell}+e^{\beta/\sqrt{A}}\beta_{0}\left(\dfrac{1}{\sqrt{A}}\beta_\ell -\dfrac{\beta}{2A^{3/2}}A_\ell \right)
-\dfrac{\beta}{2A}e^{\beta/\sqrt{A}}A_{0\ell}-\dfrac{1}{2A}e^{\beta/\sqrt{A}}A_{0}\beta_\ell\\
&+\dfrac{\beta}{2A^2}e^{\beta/\sqrt{A}}A_{0}A_\ell -\dfrac{\beta}{2A}e^{\beta/\sqrt{A}}A_{0}\left(\dfrac{1}{\sqrt{A}}\beta_\ell -\dfrac{\beta}{2A^{3/2}}A_\ell \right)
+\dfrac{1}{2\sqrt{A}}e^{\beta/\sqrt{A}}A_{0\ell}\\
&-\dfrac{1}{4A^{3/2}}e^{\beta/\sqrt{A}}A_{0}A_\ell +\dfrac{1}{2\sqrt{A}}e^{\beta/\sqrt{A}}A_{0}\left(\dfrac{1}{\sqrt{A}}\beta_\ell -\dfrac{\beta}{2A^{3/2}}A_\ell \right)
+\beta_{0\ell}\\
=&\dfrac{1}{4A^{5/2}}
\bigg[4A^{5/2}\left( e^{\beta/\sqrt{A}}+1\right) \beta_{0\ell}
+2A^2e^{\beta/\sqrt{A}}\left(2\beta_0 \beta_\ell+A_{0\ell} \right)
-2\beta A^{3/2} e^{\beta/\sqrt{A}}A_{0\ell}\\
&-Ae^{\beta/\sqrt{A}}\bigg\{A_{0}A_{\ell}+2\beta\left( A_0\beta_{\ell}+A_\ell \beta_0 \right)\bigg\}
+A^{1/2}\beta e^{\beta/\sqrt{A}} A_0 A_\ell+\beta^2e^{\beta/\sqrt{A}} A_0 A_\ell
\bigg].
\end{align*}
From the equation(\ref{55.2}), we get
\begin{align*}
\bar{F}_{x^\ell}
=&e^{\beta/\sqrt{A}}\beta_{x^\ell}-\dfrac{\beta}{2A}e^{\beta/\sqrt{A}}A_{x^\ell}
+\dfrac{1}{2\sqrt{A}}e^{\beta/\sqrt{A}}A_{x^\ell} +\beta_{x^\ell}\\
=&\dfrac{1}{4A^{5/2}}
\bigg[4A^{5/2}\left( e^{\beta/\sqrt{A}}+1\right) \beta_{x^\ell}
+2A^2e^{\beta/\sqrt{A}}A_{x^\ell}-2\beta A^{3/2} e^{\beta/\sqrt{A}}A_{x^\ell}\bigg].
\end{align*}
We know that $ \bar{F}$ is projectively flat if and only if $ \bar{F}_{x^ky^\ell}y^k-\bar{F}_{x^\ell}=0, $\\
i.e., 
\begin{align*}
&\dfrac{1}{4A^{5/2}}
\bigg[4A^{5/2}\left( e^{\beta/\sqrt{A}}+1\right) \beta_{0\ell}
+2A^2e^{\beta/\sqrt{A}}\left(2\beta_0 \beta_\ell+A_{0\ell} \right)
-2\beta A^{3/2} e^{\beta/\sqrt{A}}A_{0\ell}\\
&-Ae^{\beta/\sqrt{A}}\bigg\{A_{0}A_{\ell}+2\beta\left( A_0\beta_{\ell}+A_\ell \beta_0 \right)\bigg\}
+A^{1/2}\beta e^{\beta/\sqrt{A}} A_0 A_\ell+\beta^2e^{\beta/\sqrt{A}} A_0 A_\ell
\bigg]\\
-
&\dfrac{1}{4A^{5/2}}
\bigg[4A^{5/2}\left( e^{\beta/\sqrt{A}}+1\right) \beta_{x^\ell}
+2A^2e^{\beta/\sqrt{A}}A_{x^\ell}-2\beta A^{3/2} e^{\beta/\sqrt{A}}A_{x^\ell}\bigg]=0,
\end{align*}
i.e., 
\begin{align*}
&4A^{5/2}\left( e^{\beta/\sqrt{A}}+1\right)\left(\beta_{0\ell}-\beta_{x^\ell}\right) 
+2A^2e^{\beta/\sqrt{A}}\left(2\beta_0 \beta_\ell+A_{0\ell}-A_{x^\ell}\right)
-2\beta A^{3/2} e^{\beta/\sqrt{A}}\left( A_{0\ell}-A_{x^\ell}\right)\\&  -Ae^{\beta/\sqrt{A}}\bigg\{A_{0}A_{\ell}+2\beta\left( A_0\beta_{\ell}+A_\ell \beta_0 \right)\bigg\}
+A^{1/2}\beta e^{\beta/\sqrt{A}} A_0 A_\ell+\beta^2e^{\beta/\sqrt{A}} A_0 A_\ell=0.
\end{align*}
From the above equation, we conclude that $ \bar{F}$ is projectively flat if and only if following five equations are satisfied.
\begin{equation}{\label{PF5.1}}
\beta_{0\ell}=\beta_{x^\ell}
\end{equation}
\begin{equation}{\label{PF5.2}}
2\beta_0 \beta_\ell+A_{0\ell}-A_{x^\ell}=0
\end{equation}
\begin{equation}{\label{PF5.3}}
 A_{0\ell}-A_{x^\ell}=0\ \implies\ A_{0\ell}=A_{x^\ell}
\end{equation}
\begin{equation}{\label{PF5.4}}
A_{0}A_{\ell}+2\beta\left( A_0\beta_{\ell}+A_\ell \beta_0 \right)=0
\end{equation}
\begin{equation}{\label{PF5.5}}
A_0 A_\ell=0
\end{equation}
Further, from the equations (\ref{PF5.2}) and (\ref{PF5.3}), we get $$ \beta_0 \beta_\ell=0, $$
and from the equations (\ref{PF5.4}) and (\ref{PF5.5}), we get $$ A_0\beta_{\ell}+A_\ell \beta_0 =0.$$
The above discussion  leads to the following theorem.
\begin{theorem}
Let $(M, \bar{F})$ be an $n-$dimensional Finsler space with $\bar{F}
= F e^{ \beta/F} + \beta$ as an exponential-Randers changed metric. Then $\bar{F}$ is projectively flat if and only if the following equations are satisfied:
$$ A_{0}A_\ell=0,\ A_{0 \ell}=A_{x^\ell},\ \beta_{0}\beta_\ell=0,\ \beta_{0 \ell}=\beta_{x^\ell},\ A_\ell\beta_{0}+A_{0}\beta_\ell=0.  $$
\end{theorem}
Next, we find necessary and sufficient conditions for infinite series-Randers changed Finsler metric 
$$ \bar{F}=\dfrac{\beta^{2}}{\beta - F } + \beta $$
to be projectively flat.\\
Let us put $F^2=A$ in $\bar{F},$ then
\begin{align}{\label{56.1}}
\bar{F}=\dfrac{\beta^{2}}{\beta -\sqrt{A}} + \beta.
\end{align}
Differentiating (\ref{56.1}) w.r.t. $x^k,$ we get
\begin{align}{\label{56.2}}
\bar{F}_{x^k}
=\dfrac{2\beta}{\beta -\sqrt{A}}\beta_{x^k} 
-\dfrac{\beta^{2}}{(\beta -\sqrt{A})^2}\beta_{x^k}+\dfrac{\beta^{2}}{2\sqrt{A}(\beta -\sqrt{A})^2}A_{x^k}+\beta_{x^k}.
\end{align}
Differentiation of (\ref{56.2}) further w.r.t. $y^\ell$ gives
\begin{equation}{\label{56.3}}
\begin{split}
\bar{F}_{x^k y^\ell}
=&\dfrac{2\beta}{\beta -\sqrt{A}}\beta_{x^ky^\ell} +\dfrac{2}{\beta -\sqrt{A}}\beta_{x^k}\beta_\ell -\dfrac{2\beta}{(\beta -\sqrt{A})^2}\beta_{x^k}\left(\beta_{\ell}-\dfrac{1}{2\sqrt{A}}A_{\ell} \right)
-\dfrac{\beta^{2}}{(\beta -\sqrt{A})^2}\beta_{x^ky^\ell}\\
&-\dfrac{2\beta}{(\beta -\sqrt{A})^2}\beta_{x^k}\beta_\ell +\dfrac{2\beta^{2}}{(\beta-\sqrt{A})^3}\beta_{x^k}\left(\beta_{\ell}-\dfrac{1}{2\sqrt{A}}A_{\ell} \right)
+\dfrac{\beta^{2}}{2\sqrt{A}(\beta -\sqrt{A})^2}A_{x^ky^\ell}\\
&+\dfrac{\beta}{\sqrt{A}(\beta -\sqrt{A})^2}A_{x^k}\beta_\ell -\dfrac{\beta^{2}}{4A^{3/2}(\beta -\sqrt{A})^2}A_{x^k}A_\ell -\dfrac{\beta^{2}}{\sqrt{A}(\beta-\sqrt{A})^3}A_{x^k}\left(\beta_{\ell}-\dfrac{1}{2\sqrt{A}}A_{\ell} \right)\\&
+\beta_{x^ky^\ell}.
\end{split}
\end{equation}
Contracting (\ref{56.3}) with $y^k$, we get
\begin{align*}
L_{x^k y^\ell}y^k
=&\dfrac{2\beta}{\beta -\sqrt{A}}\beta_{0\ell} +\dfrac{2}{\beta -\sqrt{A}}\beta_{0}\beta_\ell -\dfrac{2\beta}{(\beta -\sqrt{A})^2}\beta_{0}\left(\beta_{\ell}-\dfrac{1}{2\sqrt{A}}A_{\ell} \right)
-\dfrac{\beta^{2}}{(\beta -\sqrt{A})^2}\beta_{0\ell}\\
&-\dfrac{2\beta}{(\beta -\sqrt{A})^2}\beta_{0}\beta_\ell +\dfrac{2\beta^{2}}{(\beta-\sqrt{A})^3}\beta_{0}\left(\beta_{\ell}-\dfrac{1}{2\sqrt{A}}A_{\ell} \right)
+\dfrac{\beta^{2}}{2\sqrt{A}(\beta -\sqrt{A})^2}A_{0\ell}\\
&+\dfrac{\beta}{\sqrt{A}(\beta -\sqrt{A})^2}A_{0}\beta_\ell -\dfrac{\beta^{2}}{4A^{3/2}(\beta -\sqrt{A})^2}A_{0}A_\ell -\dfrac{\beta^{2}}{\sqrt{A}(\beta-\sqrt{A})^3}A_{0}\left(\beta_{\ell}-\dfrac{1}{2\sqrt{A}}A_{\ell} \right)
+\beta_{0\ell}.
\end{align*}
Simplifying, we get
\begin{align*}
L_{x^k y^\ell}y^k
=&\dfrac{1}{4A^{3/2}(\beta-\sqrt{A})^3}
\bigg[-4A^3\beta_{0\ell}
+A^{5/2}\bigg\{ 20\beta\beta_{0\ell}+8\beta_0 \beta_\ell\bigg\}
-24A^2 \beta^2 \beta_{0\ell}\\
&+A^{3/2}\bigg\{ 8\beta^3 \beta_{0\ell}-2\beta^2A_{0\ell}-4\beta\left( A_0\beta_{\ell}+A_\ell \beta_0 \right)\bigg\}
+2\beta^3A A_{0\ell}
+3A^{1/2}\beta^2A_0A_\ell-\beta^3A_0A_\ell
\bigg].
\end{align*}
From the equation (\ref{56.2}), we get
\begin{align*}
\bar{F}_{x^\ell}
=&\dfrac{2\beta}{\beta -\sqrt{A}}\beta_{x^\ell} 
-\dfrac{\beta^{2}}{(\beta -\sqrt{A})^2}\beta_{x^\ell}+\dfrac{\beta^{2}}{2\sqrt{A}(\beta -\sqrt{A})^2}A_{x^\ell}+\beta_{x^\ell}\\
=&\dfrac{2\beta^2-4\sqrt{A}\beta+A}{(\beta-\sqrt{A})^2}\beta_{x^\ell} +\dfrac{\beta^2}{2\sqrt{A}(\beta-\sqrt{A})^2}A_{x^\ell}\\
=&\dfrac{1}{4A^{3/2}(\beta-\sqrt{A})^3}
\bigg[-4A^3\beta_{x^\ell}+20A^{5/2}\beta\beta_{x^\ell}-24A^2 \beta^2\beta_{x^\ell} +A^{3/2}\bigg\{8\beta^3\beta_{x^\ell}-2\beta^2A_{x^\ell}\bigg\}
+2A\beta^3A_{x^\ell}\bigg].
\end{align*}
We know that $ \bar{F}$ is projectively flat if and only if $ \bar{F}_{x^ky^\ell}y^k-\bar{F}_{x^\ell}=0, $\\
i.e., 
\begin{align*}
&\dfrac{1}{4A^{3/2}(\beta-\sqrt{A})^3}
\bigg[-4A^3\beta_{0\ell}
+A^{5/2}\bigg\{ 20\beta\beta_{0\ell}+8\beta_0 \beta_\ell\bigg\}
-24A^2 \beta^2 \beta_{0\ell}
+A^{3/2}\bigg\{ 8\beta^3 \beta_{0\ell}-2\beta^2A_{0\ell}\\
&-4\beta\left( A_0\beta_{\ell}+A_\ell \beta_0 \right)\bigg\}
+2\beta^3A A_{0\ell}
+3A^{1/2}\beta^2A_0A_\ell-\beta^3A_0A_\ell
\bigg]
-
\dfrac{1}{4A^{3/2}(\beta-\sqrt{A})^3}
\bigg[-4A^3\beta_{x^\ell}\\
&+20A^{5/2}\beta\beta_{x^\ell}-24A^2 \beta^2\beta_{x^\ell} +A^{3/2}\bigg\{8\beta^3\beta_{x^\ell}-2\beta^2A_{x^\ell}\bigg\}
+2A\beta^3A_{x^\ell}\bigg]=0.
\end{align*}
Simplifying, we get
\begin{align*}
&4A^3\bigg\{\beta_{x^\ell}-\beta_{0\ell}\bigg\}
+4A^{5/2}\bigg\{5\beta\left( \beta_{0\ell}-\beta_{x^\ell}\right) +2\beta_0 \beta_\ell\bigg\}
+24A^2 \beta^2\bigg\{\beta_{x^\ell}-\beta_{0\ell}\bigg\}\\&
+2A^{3/2}\bigg\{ 4\beta^3\left(\beta_{0\ell}-\beta_{x^\ell}\right) +\beta^2\left(A_{x^\ell} -A_{0\ell}\right) -2\beta\left( A_0\beta_{\ell}+A_\ell \beta_0 \right)\bigg\}
+2A\beta^3\bigg\{A_{0\ell}-A_{x^\ell}\bigg\}\\&
+3A^{1/2}\beta^2A_0A_\ell-\beta^3A_0A_\ell
=0.
\end{align*}
From the above equation, we conclude that $ \bar{F}$ is projectively flat if and only if following five equations are satisfied.
\begin{equation}{\label{PF6.1}}
\beta_{x^\ell}-\beta_{0\ell}=0\ \implies \ \beta_{x^\ell}=\beta_{0\ell}
\end{equation}
\begin{equation}{\label{PF6.2}}
5\beta\left( \beta_{0\ell}-\beta_{x^\ell}\right) +2\beta_0 \beta_\ell=0
\end{equation}
\begin{equation}{\label{PF6.3}}
 4\beta^3\left(\beta_{0\ell}-\beta_{x^\ell}\right) +\beta^2\left(A_{x^\ell} -A_{0\ell}\right) -2\beta\left( A_0\beta_{\ell}+A_\ell \beta_0 \right)=0
\end{equation}
\begin{equation}{\label{PF6.4}}
A_{0\ell}-A_{x^\ell}=0\ \implies\ A_{0\ell}=A_{x^\ell}
\end{equation}
\begin{equation}{\label{PF6.5}}
A_0A_\ell=0.
\end{equation}
Further, from the equations (\ref{PF6.1}) and (\ref{PF6.2}), we get
\begin{equation}{\label{PF6.6}}
\beta_0 \beta_\ell=0
\end{equation}
and from the equations (\ref{PF6.1}), (\ref{PF6.4}) and (\ref{PF6.3}), we get
\begin{equation}{\label{PF6.7}}
 A_0\beta_{\ell}+A_\ell \beta_0=0.
\end{equation}
Above discussion  leads to the following theorem.
\begin{theorem}
Let $(M, \bar{F})$ be an $n-$dimensional Finsler space with $ \bar{F}=\dfrac{\beta^{2}}{\beta - F } + \beta$ as an infinite series-Randers changed  metric. Then $\bar{F}$  is projectively flat if and only if the following  equations are satisfied:
		$$ A_{0}A_\ell=0,\ A_{0 \ell}=A_{x^\ell},\ \beta_{0}\beta_\ell=0,\ \beta_{0 \ell}=\beta_{x^\ell},\ A_\ell\beta_{0}+A_{0}\beta_\ell=0.  $$
\end{theorem}
\section{Dually flatness of Finsler metrics}
Firstly, we recall \cite{ShenRFGAIG} the following definition:
\begin{definition}
A Finsler metric $F$ on a smooth $n-$dimensional manifold $M$ is called locally dually flat if, at any point, there is a standard co-ordinate system $ (x^i, y^i) $ in $TM,$ $(x^i)$ called adapted local co-ordinate system, such that 
$$ L_{x^ky^\ell}y^k-2L_{x^\ell}=0,\  \text{where}\ L=F^2.$$
\end{definition}
Next, we find the necessary and sufficient conditions for locally dually flatness of all the metrics constructed via Randers change in section two.\\
First, we find necessary and sufficient conditions for Kropina-Randers changed Finsler metric
\begin{equation*}
\bar{F}= \dfrac{F^{2}}{\beta} + \beta
\end{equation*}
to be locally dually flat.\\
Let us put $F^2=A$ in $\bar{F},$ then
\begin{equation*}
\bar{F}= \dfrac{A}{\beta} + \beta.
\end{equation*}
\begin{equation}{\label{61.1}}
L=\bar{F}^2=\dfrac{A^2}{\beta^2}+2A+\beta^2.
\end{equation}
Differentiating (\ref{61.1}) w.r.t. $x^k,$ we get
\begin{equation}{\label{61.2}}
L_{x^k}=\dfrac{2A}{\beta^2}A_{x^k}-\dfrac{2A^2}{\beta^3}\beta_{x^k}+2A_{x^k}+2\beta \beta_{x^k}.
\end{equation}
Differentiation of (\ref{61.2}) further w.r.t. $y^\ell$ gives
\begin{equation}{\label{61.3}}
\begin{split}
L_{x^k y^\ell}
&=\dfrac{2A}{\beta^2}A_{x^k y^\ell} +\dfrac{2}{\beta^2}A_{x^k}A_\ell-\dfrac{4A}{\beta^3}\beta_\ell A_{x^k}
-\dfrac{2A^2}{\beta^3}\beta_{x^k y^\ell}-\dfrac{4A}{\beta^3}A_\ell \beta_{x^k}+ \dfrac{6A^2}{\beta^4}\beta_\ell \beta_{x^k}
+2A_{x^k y^\ell}\\&
+2\beta \beta_{x^k y^\ell}+2\beta_\ell\beta_{x^k}.
\end{split}
\end{equation}
Contracting (\ref{61.3}) with $y^k$, we get
\begin{align*}
L_{x^k y^\ell}y^k
&=\dfrac{2A}{\beta^2}A_{0 \ell} +\dfrac{2}{\beta^2}A_{0}A_\ell-\dfrac{4A}{\beta^3}\beta_\ell A_{0}-\dfrac{2A^2}{\beta^3}\beta_{0 \ell}-\dfrac{4A}{\beta^3}A_\ell \beta_{0}+ \dfrac{6A^2}{\beta^4}\beta_\ell \beta_{0}+2A_{0 \ell}+2\beta \beta_{0 \ell}+2\beta_\ell\beta_{0}\\
=&\dfrac{2}{\beta^4}\bigg[A^2\left( 3 \beta_\ell \beta_{0}-\beta\beta_{0 \ell}\right) + A\left( \beta^2A_{0 \ell} -2\beta \left( \beta_\ell A_{0}-A_\ell \beta_{0}\right) \right) +\beta^5 \beta_{0 \ell}+\beta^4\left( A_{0 \ell}
+\beta_\ell\beta_{0}\right)  +\beta^2 A_{0}A_\ell \bigg].
\end{align*}
Further, equation (\ref{61.2}) can be rewritten as
\begin{align*}
2 L_{x^\ell}
=&\dfrac{4A}{\beta^2}A_{x^\ell}-\dfrac{4A^2}{\beta^3}\beta_{x^\ell}+4A_{x^\ell}+4\beta \beta_{x^\ell}\\
=&\dfrac{2}{\beta^4}\bigg[-2A^2 \beta \beta_{x^\ell} + 2A \beta^2 A_{x^\ell}+2 \beta^4 A_{x^\ell}+2\beta^5 \beta_{x^\ell} \bigg].
\end{align*}
We know that $ \bar{F}$ is locally dually flat if and only if $ L_{x^ky^\ell}y^k-2L_{x^\ell}=0, $\\
i.e., 
\begin{align*}
\dfrac{2}{\beta^4}&\bigg[A^2\left( 3 \beta_\ell \beta_{0}-\beta\beta_{0 \ell}\right) + A\left( \beta^2A_{0 \ell} -2\beta \left( \beta_\ell A_{0}-A_\ell \beta_{0}\right) \right) +\beta^5 \beta_{0 \ell}+\beta^4\left( A_{0 \ell}
+\beta_\ell\beta_{0}\right)  +\beta^2 A_{0}A_\ell \bigg]\\& - \dfrac{2}{\beta^4}\bigg[  2A \beta^2 A_{x^\ell}-2A^2 \beta \beta_{x^\ell}+2 \beta^4 A_{x^\ell}+2\beta^5 \beta_{x^\ell} \bigg]=0,
\end{align*}
i.e.,
\begin{align*}
&A^2\left( 3 \beta_\ell \beta_{0}-\beta\beta_{0 \ell}+2 \beta \beta_{x^\ell} \right) + A\left\{ \beta^2\left( A_{0 \ell}- 2 A_{x^\ell}\right)  -2\beta \left( \beta_\ell A_{0}-A_\ell \beta_{0}\right) \right\} +\beta^5\left(  \beta_{0 \ell}-2 \beta_{x^\ell}\right) \\
&+\beta^4\left( A_{0 \ell}
+\beta_\ell\beta_{0} -2 A_{x^\ell}\right)  +\beta^2 A_{0}A_\ell =0.
\end{align*}
From the above equation, we conclude that $ \bar{F}$ is locally dually flat if and only if following three equations are satisfied.
\begin{equation}{\label{LDF1.1}}
3\beta_\ell \beta_{0}-\beta\beta_{0 \ell}+2 \beta \beta_{x^\ell}=0
\end{equation}
\begin{equation}{\label{LDF1.2}}
\beta^2\left( A_{0 \ell}- 2 A_{x^\ell}\right)  -2\beta \left( \beta_\ell A_{0}-A_\ell \beta_{0}\right)=0
\end{equation}
\begin{equation}{\label{LDF1.3}}
\beta^5\left(  \beta_{0 \ell}-2 \beta_{x^\ell}\right) +\beta^4\left( A_{0 \ell}
+\beta_\ell\beta_{0} -2 A_{x^\ell}\right)  +\beta^2 A_{0}A_\ell =0
\end{equation}
Hence, we have the following theorem.
 \begin{theorem}
Let $(M, \bar{F})$ be an $n-$dimensional Finsler space with $\bar{F}= \dfrac{F^{2}}{\beta} + \beta$ as a Kropina-Randers changed metric. Then $\bar{F}$ is locally dually flat if and only if equations (\ref{LDF1.1}), (\ref{LDF1.2}) and (\ref{LDF1.3}) are satisfied.
\end{theorem}
Next, we find necessary and sufficient conditions for generalized Kropina-Randers changed  Finsler metric 
$$ \bar{F}= \dfrac{F^{m+1}}{\beta^m} + \beta \  \left( m\neq 0, -1\right) $$
to be locally dually flat.\\ 
Let us put $ F^2=A$ in $\bar{F},$ then
\begin{equation*}
\bar{F}= \dfrac{A^{(m+1)/2}}{\beta^m} + \beta. 
\end{equation*}
\begin{align}{\label{62.1}}
L=\bar{F}^2
&=\dfrac{A^{m+1}}{\beta^{2m}}+2\dfrac{A^{(m+1)/2}}{\beta^{m-1}}+\beta^2.
\end{align}
Differentiating (\ref{62.1}) w.r.t. $x^k,$ we get
\begin{equation}{\label{62.2}}
L_{x^k}
=(m+1)\dfrac{A^{m}}{\beta^{2m}}A_{x^k}-2m\dfrac{A^{m+1}}{\beta^{2m+1}}\beta_{x^k}
+(m+1)\dfrac{A^{(m-1)/2}}{\beta^{m-1}}A_{x^k}-2(m-1)\dfrac{A^{(m+1)/2}}{\beta^{m}}\beta_{x^k}
+2\beta \beta_{x^k}.
\end{equation}
Differentiation of (\ref{62.2}) further  w.r.t. $y^\ell$ gives
\begin{equation}{\label{62.3}}
\begin{split}
L_{x^k y^\ell}
=&(m+1)\dfrac{A^{m}}{\beta^{2m}}A_{x^k y^\ell}+m(m+1)\dfrac{A^{m-1}}{\beta^{2m}}A_\ell A_{x^k}-2m(m+1)\dfrac{A^{m}}{\beta^{2m+1}}A_{x^k}\beta_\ell\\
&-2m\dfrac{A^{m+1}}{\beta^{2m+1}}\beta_{x^k y^\ell} -2m(m+1)\dfrac{A^{m}}{\beta^{2m+1}}A_\ell\beta_{x^k} +2m(2m+1)\dfrac{A^{m+1}}{\beta^{2m+2}}\beta_{x^k}\beta_\ell\\
&+(m+1)\dfrac{A^{(m-1)/2}}{\beta^{m-1}}A_{x^k y^\ell} +\dfrac{m^2-1}{2}\dfrac{A^{(m-3)/2}}{\beta^{m-1}}A_{x^k}A_\ell -(m^2-1)\dfrac{A^{(m-1)/2}}{\beta^{m}}A_{x^k}\beta_\ell\\
&-2(m-1)\dfrac{A^{(m+1)/2}}{\beta^{m}}\beta_{x^k y^\ell} -(m^2-1)\dfrac{A^{(m-1)/2}}{\beta^{m}}\beta_{x^k}A_\ell +2m(m-1)\dfrac{A^{(m+1)/2}}{\beta^{m+1}}\beta_{x^k}\beta_\ell\\
&+2\beta \beta_{x^k y^\ell}+2\beta_\ell \beta_{x^k}.
\end{split}
\end{equation}
Contracting (\ref{62.3}) with $y^k$, we get
\begin{align*}
L_{x^k y^\ell}y^k
=&(m+1)\dfrac{A^{m}}{\beta^{2m}}A_{0 \ell}+m(m+1)\dfrac{A^{m-1}}{\beta^{2m}}A_\ell A_{0}-2m(m+1)\dfrac{A^{m}}{\beta^{2m+1}}A_{0}\beta_\ell\\
&-2m\dfrac{A^{m+1}}{\beta^{2m+1}}\beta_{0 \ell} -2m(m+1)\dfrac{A^{m}}{\beta^{2m+1}}A_\ell\beta_{0} +2m(2m+1)\dfrac{A^{m+1}}{\beta^{2m+2}}\beta_{0}\beta_\ell\\
&+(m+1)\dfrac{A^{(m-1)/2}}{\beta^{m-1}}A_{0 \ell} +\dfrac{m^2-1}{2}\dfrac{A^{(m-3)/2}}{\beta^{m-1}}A_{0}A_\ell -(m^2-1)\dfrac{A^{(m-1)/2}}{\beta^{m}}A_{0}\beta_\ell\\
&-2(m-1)\dfrac{A^{(m+1)/2}}{\beta^{m}}\beta_{0 \ell} -(m^2-1)\dfrac{A^{(m-1)/2}}{\beta^{m}}\beta_{0}A_\ell +2m(m-1)\dfrac{A^{(m+1)/2}}{\beta^{m+1}}\beta_{0}\beta_\ell\\
&+2\beta \beta_{0 \ell}+2\beta_\ell \beta_{0}.
\end{align*}
After simplification, we get
\begin{align*}
L_{x^k y^\ell}y^k
=&\dfrac{A^{m-1}}{\beta^{2m+2}}\bigg[ (m+1) A \beta^2 A_{0 \ell}+m(m+1) \beta^{2} A_\ell A_{0}-2m(m+1)A\beta A_{0}\beta_\ell
-2mA^{2}\beta \beta_{0 \ell}\\
&-2m(m+1) A \beta A_\ell\beta_{0} +2m(2m+1)A^{2} \beta_{0}\beta_\ell  \bigg]
+\dfrac{A^{(m-3)/2}}{\beta^{m+1}}   \bigg[(m+1)A \beta^{2} A_{0 \ell}\\
&+\dfrac{m^2-1}{2} \beta^{2}A_{0}A_\ell -(m^2-1)A\beta A_{0}\beta_\ell
-2(m-1)A^2\beta \beta_{0 \ell} -(m^2-1)A \beta \beta_{0}A_\ell\\
&+2m(m-1)A^2 \beta_{0}\beta_\ell \bigg]+2\beta \beta_{0 \ell}+2\beta_\ell \beta_{0}.
\end{align*}
Further, equation (\ref{62.2}) can be rewritten as 
\begin{align*}
2 L_{x^\ell}
=&\dfrac{A^{m-1}}{\beta^{2m+2}}\bigg[2(m+1)A \beta^{2}A_{x^\ell}-4m A^{2}\beta\beta_{x^\ell}\bigg]
+\dfrac{A^{(m-3)/2}}{\beta^{m+1}}\bigg[2(m+1)A\beta^{2}A_{x^\ell}-4(m-1)A^{2}\beta\beta_{x^\ell}\bigg]\\&
+4\beta \beta_{x^\ell}.
\end{align*}
We know that $ \bar{F}$ is locally dually flat if and only if $ L_{x^ky^\ell}y^k-2L_{x^\ell}=0, $\\
i.e., 
\begin{align*}
&\dfrac{A^{m-1}}{\beta^{2m+2}}\bigg[ (m+1) A \beta^2 A_{0 \ell}+m(m+1) \beta^{2} A_\ell A_{0}-2m(m+1)A\beta A_{0}\beta_\ell
-2mA^{2}\beta \beta_{0 \ell}-2m(m+1) A \beta A_\ell\beta_{0}\\
&+2m(2m+1)A^{2} \beta_{0}\beta_\ell  \bigg]
+\dfrac{A^{(m-3)/2}}{\beta^{m+1}}   \bigg[(m+1)A \beta^{2} A_{0 \ell}
+\dfrac{m^2-1}{2} \beta^{2}A_{0}A_\ell -(m^2-1)A\beta A_{0}\beta_\ell\\&
-2(m-1)A^2\beta \beta_{0 \ell} -(m^2-1)A \beta \beta_{0}A_\ell
+2m(m-1)A^2 \beta_{0}\beta_\ell \bigg]+2\beta \beta_{0 \ell}+2\beta_\ell \beta_{0}\\&
-
\dfrac{A^{m-1}}{\beta^{2m+2}}\bigg[2(m+1)A \beta^{2}A_{x^\ell}-4m A^{2}\beta\beta_{x^\ell}\bigg]
-\dfrac{A^{(m-3)/2}}{\beta^{m+1}}\bigg[2(m+1)A\beta^{2}A_{x^\ell}-4(m-1)A^{2}\beta\beta_{x^\ell}\bigg]\\&
-4\beta \beta_{x^\ell}=0,
\end{align*}
i.e., 
\begin{align*}
&\dfrac{A^{m-1}}{\beta^{2m+2}}\bigg[ 2m A^{2}\bigg\{ \beta\left(2\beta_{x^\ell}-\beta_{0 \ell}\right)+(2m+1)\beta_{0}\beta_\ell\bigg\}+(m+1)A\bigg\{\beta^2\left(A_{0 \ell}-2A_{x^\ell}\right) -2m \beta\left(A_\ell\beta_{0}+A_{0}\beta_\ell\right)\bigg\}\\
&+m(m+1) \beta^{2} A_\ell A_{0}\bigg]
+\dfrac{A^{(m-3)/2}}{\beta^{m+1}}   \bigg[2(m-1)A^{2}\bigg\{\beta\left( 2\beta_{x^\ell}- \beta_{0 \ell}\right) +m\beta_{0}\beta_\ell\bigg\}+(m+1)A \bigg\{\beta^{2}\left(A_{0 \ell}-2A_{x^\ell}\right) \\
&-(m-1)\beta\left(A_{0}\beta_\ell+\beta_{0}A_\ell\right)\bigg\}+\dfrac{m^2-1}{2} \beta^{2}A_{0}A_\ell\bigg] +2\beta \left(\beta_{0 \ell}-2\beta_{x^\ell}\right)+2\beta_\ell \beta_{0}
=0.
\end{align*}
From the above equation, we conclude that $ \bar{F}$ is locally dually flat if and only if following seven equations are satisfied.
\begin{equation}{\label{LDF2.1}}
\beta\left(2\beta_{x^\ell}-\beta_{0 \ell}\right)+(2m+1)\beta_{0}\beta_\ell=0
\end{equation}
\begin{equation}{\label{LDF2.2}}
\beta^2\left(A_{0 \ell}-2A_{x^\ell}\right) -2m \beta\left(A_\ell\beta_{0}+A_{0}\beta_\ell\right)=0
\end{equation}
\begin{equation}{\label{LDF2.3}}
m(m+1) \beta^{2} A_\ell A_{0}=0\  \ \implies A_\ell A_{0}=0
\end{equation}
\begin{equation}{\label{LDF2.4}}
(m-1)\bigg\{\beta\left( 2\beta_{x^\ell}- \beta_{0 \ell}\right) +m\beta_{0}\beta_\ell\bigg\}=0
\end{equation}
\begin{equation}{\label{LDF2.5}}
\beta^{2}\left(A_{0 \ell}-2A_{x^\ell}\right) -(m-1)\beta\left(A_{0}\beta_\ell+\beta_{0}A_\ell\right)=0
\end{equation}
\begin{equation}{\label{LDF2.6}}
(m^2-1) \beta^{2}A_{0}A_\ell=0 \ \implies (m-1) A_{0}A_\ell=0
\end{equation}
\begin{equation}{\label{LDF2.7}}
\beta \left(\beta_{0 \ell}-2\beta_{x^\ell}\right)+\beta_\ell \beta_{0}
=0.
\end{equation}
Further, from the equation (\ref{LDF2.4}), we see that \\
either $m=1$ or $ \beta\left( 2\beta_{x^\ell}- \beta_{0 \ell}\right) +m\beta_{0}\beta_\ell=0 $\\ 
Now, if $m=1,$ then  
(\ref{LDF2.1}) reduces to (\ref{LDF1.1}),  
(\ref{LDF2.2}) reduces to (\ref{LDF1.2}), 
and (\ref{LDF2.3}),(\ref{LDF2.5}),(\ref{LDF2.7}) reduce to (\ref{LDF1.3}).\\
But the equations (\ref{LDF1.1}), (\ref{LDF1.2}) and (\ref{LDF1.3}) are necessary and sufficient conditions for Kropina-Randers changed Finsler metric to be locally dually flat. Therefore, we assume that $m\neq1$ for the general case.\\
Then 
\begin{equation}{\label{LDF2.8}}
\beta\left( 2\beta_{x^\ell}- \beta_{0 \ell}\right) +m\beta_{0}\beta_\ell=0.
\end{equation}
From the equations (\ref{LDF2.1}) and (\ref{LDF2.8}), we get
$$ \beta_{0}\beta_\ell=0. $$
Again from the equation (\ref{LDF2.1}), we get $$ \beta_{0 \ell}=2\beta_{x^\ell}. $$
Also from the equations (\ref{LDF2.2}) and (\ref{LDF2.5}), we get
$$ A_\ell\beta_{0}+A_{0}\beta_\ell=0. $$
Again from the equation (\ref{LDF2.2}), we get $$ A_{0 \ell}=2A_{x^\ell}. $$
Above discussion  leads to the following theorem.
\begin{theorem}
Let $(M, \bar{F})$ be an $n-$dimensional Finsler space with $ \bar{F}= \dfrac{F^{m+1}}{\beta^m} + \beta \  \left( m\neq -1, 0, 1\right) $ as generalized Kropina-Randers changed  metric. Then  $ \bar{F}$ is locally dually flat if and only if the following equations are satisfied:
	$$ A_{0}A_\ell=0,\ A_{0 \ell}=2A_{x^\ell},\ \beta_{0}\beta_\ell=0,\ \beta_{0 \ell}=2\beta_{x^\ell},\ A_\ell\beta_{0}+A_{0}\beta_\ell=0.  $$
\end{theorem}
Next, we find necessary and sufficient conditions for square-Randers changed Finsler metric 
$$ \bar{F}= \dfrac{(F+\beta)^{2}}{F} + \beta$$
to be locally dually flat.\\
Let us put $ F^2=A$ in $\bar{F},$ then
\begin{align*}
\bar{F}=A^{1/2}+\dfrac{\beta^2}{A^{1/2}}+3\beta.
\end{align*}
\begin{equation}{\label{63.1}}
\begin{split}
L=\bar{F}^2
&=\dfrac{\beta^4}{A}+\dfrac{6 \beta^3}{A^{1/2}}+6\beta A^{1/2}+A+11\beta^2.
\end{split}
\end{equation}
Differentiating (\ref{63.1}) w.r.t. $x^k,$ we get
\begin{align}{\label{63.2}}
L_{x^k}
=\dfrac{4\beta^3}{A} \beta_{x^k}-\dfrac{\beta^4}{A^2}A_{x^k}
+\dfrac{18 \beta^2}{A^{1/2}}\beta_{x^k}-\dfrac{3 \beta^3}{A^{3/2}} A_{x^k}
+6A^{1/2}\beta_{x^k}+\dfrac{3\beta}{A^{1/2}}A_{x^k}
+A_{x^k}+22\beta \beta_{x^k}.
\end{align}
Differentiation of (\ref{63.2}) further w.r.t. $y^\ell$ gives
\begin{align*}
L_{x^k y^\ell}
=&\dfrac{4\beta^3}{A} \beta_{x^k y^\ell}+\dfrac{12\beta^2}{A} \beta_{x^k}\beta_\ell -\dfrac{4\beta^3}{A^2} \beta_{x^k}A_\ell
-\dfrac{\beta^4}{A^2}A_{x^k y^\ell}-\dfrac{4\beta^3}{A^2}A_{x^k}\beta_\ell +\dfrac{2\beta^4}{A^3}A_{x^k}A_\ell
+\dfrac{18 \beta^2}{A^{1/2}}\beta_{x^k y^\ell}\\
&+\dfrac{36 \beta}{A^{1/2}}\beta_{x^k}\beta_\ell -\dfrac{9 \beta^2}{A^{3/2}}\beta_{x^k}A_\ell
-\dfrac{3\beta^3}{A^{3/2}} A_{x^k y^\ell}-\dfrac{9\beta^2}{A^{3/2}} A_{x^k}\beta_\ell+\dfrac{9\beta^3}{2A^{5/2}} A_{x^k}A_\ell
+6A^{1/2}\beta_{x^k y^\ell}\\
&+\dfrac{3}{A^{1/2}}\beta_{x^k}A_\ell
+\dfrac{3\beta}{A^{1/2}}A_{x^k y^\ell}+\dfrac{3}{A^{1/2}}A_{x^k}\beta_\ell -\dfrac{3\beta}{2A^{3/2}}A_{x^k}A_\ell
+A_{x^k y^\ell}+22\beta \beta_{x^k y^\ell}+22\beta_{x^k}\beta_\ell.
\end{align*}
Simplifying, we get
\begin{equation}{\label{63.3}}
\begin{split}
L_{x^k y^\ell}=&\dfrac{1}{A^3}\bigg[4\beta^3A^2 \beta_{x^k y^\ell}+12\beta^2A^2 \beta_{x^k}\beta_\ell -4\beta^3 A \beta_{x^k}A_\ell-\beta^4 A A_{x^k y^\ell} -4\beta^3 A A_{x^k}\beta_\ell +2\beta^4 A_{x^k}A_\ell\bigg]\\
&+\dfrac{3}{2A^{5/2}}\bigg[12\beta^2 A^2\beta_{x^k y^\ell}
+24 \beta A^{2}\beta_{x^k}\beta_\ell -6 \beta^2 A \beta_{x^k}A_\ell
-2\beta^3 A A_{x^k y^\ell}-6\beta^2AA_{x^k}\beta_\ell +3\beta^3 A_{x^k}A_\ell\\
+&4A^{3}\beta_{x^k y^\ell}+2A^{2}\beta_{x^k}A_\ell
+2\beta A^{2}A_{x^k y^\ell}+2A^{2}A_{x^k}\beta_\ell -\beta A A_{x^k}A_\ell\bigg]
+A_{x^k y^\ell}+22\beta \beta_{x^k y^\ell}+22\beta_{x^k}\beta_\ell.
\end{split}
\end{equation}	
Contracting (\ref{63.3}) with $y^k$, we get
\begin{align*}
L_{x^k y^\ell} y^k
=&\dfrac{1}{A^3}\bigg[4\beta^3A^2 \beta_{0\ell}+12\beta^2A^2 \beta_{0}\beta_\ell -4\beta^3 A \beta_{0}A_\ell-\beta^4 A A_{0\ell} -4\beta^3 A A_{0}\beta_\ell +2\beta^4 A_{0}A_\ell\bigg]\\
&+\dfrac{3}{2A^{5/2}}\bigg[12\beta^2 A^2\beta_{0 \ell}
+24 \beta A^{2}\beta_{0}\beta_\ell -6 \beta^2 A \beta_{0}A_\ell
-2\beta^3 A A_{0\ell}-6\beta^2AA_{0}\beta_\ell+3\beta^3 A_{0}A_\ell\\
&+4A^{3}\beta_{0\ell}+2A^{2}\beta_{0}A_\ell
+2\beta A^{2}A_{0\ell}+2A^{2}A_{0}\beta_\ell -\beta A A_{0}A_\ell\bigg]
+A_{0\ell}+22\beta \beta_{0\ell}+22\beta_{0}\beta_\ell.
\end{align*}
Further,  equation (\ref{63.2}) can be rewritten as 
\begin{align*}
2L_{x^\ell}
=&\dfrac{1}{A^3}\bigg[8\beta^3 A^2\beta_{x^\ell}-2\beta^4AA_{x^\ell}\bigg]
+\dfrac{3}{2A^{5/2}}\bigg[24\beta^2A^2\beta_{x^\ell}-4\beta^3A A_{x^\ell}
+8 A^3\beta_{x^\ell}+4\beta A^{2}A_{x^\ell}\bigg]\\
&+2A_{x^\ell}+44\beta \beta_{x^\ell}.
\end{align*}
We know that $ \bar{F}$ is locally dually flat if and only if $ L_{x^ky^\ell}y^k-2L_{x^\ell}=0, $\\
i.e., 
\begin{align*}
&\dfrac{1}{A^3}\bigg[4\beta^3A^2 \beta_{0\ell}+12\beta^2A^2 \beta_{0}\beta_\ell -4\beta^3 A \beta_{0}A_\ell-\beta^4 A A_{0\ell} -4\beta^3 A A_{0}\beta_\ell +2\beta^4 A_{0}A_\ell\bigg]\\
&+\dfrac{3}{2A^{5/2}}\bigg[12\beta^2 A^2\beta_{0 \ell}
+24 \beta A^{2}\beta_{0}\beta_\ell -6 \beta^2 A \beta_{0}A_\ell
-2\beta^3 A A_{0\ell}-6\beta^2AA_{0}\beta_\ell+3\beta^3 A_{0}A_\ell\\
&+4A^{3}\beta_{0\ell}+2A^{2}\beta_{0}A_\ell
+2\beta A^{2}A_{0\ell}+2A^{2}A_{0}\beta_\ell -\beta A A_{0}A_\ell\bigg]
+A_{0\ell}+22\beta \beta_{0\ell}+22\beta_{0}\beta_\ell\\
-
&\dfrac{1}{A^3}\bigg[8\beta^3 A^2\beta_{x^\ell}-2\beta^4AA_{x^\ell}\bigg]
-\dfrac{3}{2A^{5/2}}\bigg[24\beta^2A^2\beta_{x^\ell}-4\beta^3A A_{x^\ell}
+8 A^3\beta_{x^\ell}+4\beta A^{2}A_{x^\ell}\bigg]\\
&-2A_{x^\ell}-44\beta \beta_{x^\ell}=0.
\end{align*}
After simplification, we get
\begin{align*}
&\dfrac{1}{A^3}\bigg[4A^2\bigg\{\beta^3\left(\beta_{0\ell}-2\beta_{x^\ell}\right)+3\beta^2 \beta_{0}\beta_\ell\bigg\}+A\bigg\{\beta^4\left( 2A_{x^\ell}-A_{0\ell}\right) -4\beta^3\left( \beta_{0}A_\ell+A_{0}\beta_\ell\right) \bigg\} +2\beta^4 A_{0}A_\ell\bigg]\\
&+\dfrac{3}{2A^{5/2}}
\bigg[4A^{3}\bigg\{\beta_{0\ell}-2\beta_{x^\ell}\bigg\}+2A^2\bigg\{6\beta^2\left(\beta_{0 \ell}-2\beta_{x^\ell}\right)
+\beta\left(A_{0\ell}-2A_{x^\ell}+12\beta_{0}\beta_\ell\right)+\left(\beta_{0}A_\ell
+A_{0}\beta_\ell\right)\bigg\} \\
&+A\bigg\{2\beta^3\left(2A_{x^\ell}-A_{0\ell}\right) -6 \beta^2\left(\beta_{0}A_\ell
+A_{0}\beta_\ell\right)-\beta A_{0}A_\ell\bigg\}
+3\beta^3 A_{0}A_\ell\bigg]
+A_{0\ell}-2A_{x^\ell}+22\beta_{0}\beta_\ell\\
&+22\beta\left(\beta_{0\ell}
-2\beta_{x^\ell}\right) =0.
\end{align*}
From the above equation, we conclude that $ \bar{F}$ is locally dually flat if and only if following eight equations are satisfied.
\begin{equation}{\label{LDF3.1}}
\beta^3\left(\beta_{0\ell}-2\beta_{x^\ell}\right)+3\beta^2 \beta_{0}\beta_\ell=0
\end{equation} 
\begin{equation}{\label{LDF3.2}}
\beta^4\left( 2A_{x^\ell}-A_{0\ell}\right) -4\beta^3\left( \beta_{0}A_\ell+A_{0}\beta_\ell\right)=0
\end{equation} 
\begin{equation}{\label{LDF3.3}}
A_{0}A_\ell=0
\end{equation} 
\begin{equation}{\label{LDF3.4}}
\beta_{0\ell}-2\beta_{x^\ell}=0\ \implies\ \beta_{0\ell}=2\beta_{x^\ell}
\end{equation} 
\begin{equation}{\label{LDF3.5}}
6\beta^2\left(\beta_{0 \ell}-2\beta_{x^\ell}\right)
+\beta\left(A_{0\ell}-2A_{x^\ell}+12\beta_{0}\beta_\ell\right)+\left(\beta_{0}A_\ell
+A_{0}\beta_\ell\right)=0
\end{equation} 
\begin{equation}{\label{LDF3.6}}
2\beta^3\left(2A_{x^\ell}-A_{0\ell}\right) -6 \beta^2\left(\beta_{0}A_\ell
+A_{0}\beta_\ell\right)-\beta A_{0}A_\ell=0
\end{equation} 
\begin{equation}{\label{LDF3.7}}
A_{0}A_\ell=0
\end{equation} 
\begin{equation}{\label{LDF3.8}}
A_{0\ell}-2A_{x^\ell}+22\beta_{0}\beta_\ell+22\beta\left(\beta_{0\ell}
-2\beta_{x^\ell}\right)=0.
\end{equation}  
Further, from the equations (\ref{LDF3.1}) and  (\ref{LDF3.4}), we get
\begin{equation}{\label{LDF3.9}}
\beta_{0}\beta_\ell=0.
\end{equation} 
Again from the equations (\ref{LDF3.9}), (\ref{LDF3.4}) and  (\ref{LDF3.8}), we get
\begin{equation}{\label{LDF3.10}}
A_{0\ell}=2A_{x^\ell},
\end{equation} 
and from the equations (\ref{LDF3.10}), (\ref{LDF3.3}) and  (\ref{LDF3.6}), we get
\begin{equation}{\label{LDF3.11}}
\beta_{0}A_\ell+A_{0}\beta_\ell=0.
\end{equation} 
Above discussion  leads to the following theorem.
\begin{theorem}
Let $(M, \bar{F})$ be an $n-$dimensional Finsler space with $ \bar{F}= \dfrac{(F+\beta)^{2}}{F} + \beta$ as a square-Randers changed  metric. Then $ \bar{F}$ is locally dually flat if and only if the following equations are satisfied:
	$$ A_{0}A_\ell=0,\ A_{0 \ell}=2A_{x^\ell},\ \beta_{0}\beta_\ell=0,\ \beta_{0 \ell}=2\beta_{x^\ell},\ A_\ell\beta_{0}+A_{0}\beta_\ell=0.  $$
\end{theorem}
Next, we find necessary and sufficient conditions for Matsumoto-Randers changed Finsler metric 
$$ \bar{F}= \dfrac{F^{2}}{F-\beta} + \beta$$
to be locally dually flat.\\ 
Let us put $F^2=A$ in $\bar{F},$ then
$$ \bar{F}= \dfrac{A}{\sqrt{A}-\beta} + \beta.$$
\begin{equation}{\label{64.1}}
L=\bar{F}^2= \dfrac{A^2}{(\sqrt{A}-\beta)^2}+ \dfrac{2A\beta}{\sqrt{A}-\beta}+ \beta^2.
\end{equation}
Differentiating (\ref{64.1}) w.r.t. $x^k,$ we get
\begin{equation}{\label{64.2}}
\begin{split}
 L_{x^k}
=&\dfrac{2A}{(\sqrt{A}-\beta)^2} A_{x^k} -\dfrac{A^{3/2}}{(\sqrt{A}-\beta)^3}A_{x^k}+\dfrac{2A^2}{(\sqrt{A}-\beta)^3}\beta_{x^k}
+\dfrac{2\beta}{\sqrt{A}-\beta}A_{x^k}+\dfrac{2A}{\sqrt{A}-\beta} \beta_{x^k}\\& -\dfrac{\sqrt{A}\beta}{(\sqrt{A}-\beta)^2}A_{x^k}+\dfrac{2A\beta}{(\sqrt{A}-\beta)^2}\beta_{x^k}+ 2\beta \beta_{x^k}.
\end{split}
\end{equation}
Differentiation of  (\ref{64.2}) further w.r.t. $y^\ell$ gives
\begin{equation}{\label{64.3}}
\begin{split}
L_{x^k y^\ell}
=&\dfrac{2A}{(\sqrt{A}-\beta)^2} A_{x^k y^\ell}+\dfrac{2}{(\sqrt{A}-\beta)^2} A_{x^k}A_\ell -\dfrac{4A}{(\sqrt{A}-\beta)^3} A_{x^k}\left(\dfrac{1}{2\sqrt{A}}A_\ell-\beta_\ell \right) \\&
-
\dfrac{A^{3/2}}{(\sqrt{A}-\beta)^3}A_{x^k y^\ell} -\dfrac{3A^{1/2}}{2(\sqrt{A}-\beta)^3}A_{x^k}A_\ell +\dfrac{3A^{3/2}}{(\sqrt{A}-\beta)^4}A_{x^k}\left(\dfrac{1}{2\sqrt{A}}A_\ell-\beta_\ell\right)\\& 
+
\dfrac{2A^2}{(\sqrt{A}-\beta)^3}\beta_{x^k y^\ell}+\dfrac{4A}{(\sqrt{A}-\beta)^3}\beta_{x^k}A_\ell -\dfrac{6A^2}{(\sqrt{A}-\beta)^4}\beta_{x^k}\left(\dfrac{1}{2\sqrt{A}}A_\ell-\beta_\ell \right)\\& 
+
\dfrac{2\beta}{\sqrt{A}-\beta}A_{x^k y^\ell}+\dfrac{2}{\sqrt{A}-\beta}A_{x^k}\beta_\ell -\dfrac{2\beta}{(\sqrt{A}-\beta)^2}A_{x^k}\left(\dfrac{1}{2\sqrt{A}}A_\ell-\beta_\ell \right)\\& 
+
\dfrac{2A}{\sqrt{A}-\beta}\beta_{x^k y^\ell}+\dfrac{2}{\sqrt{A}-\beta}\beta_{x^k}A_\ell
-\dfrac{2A}{(\sqrt{A}-\beta)^2}\beta_{x^k}\left(\dfrac{1}{2\sqrt{A}}A_\ell-\beta_\ell \right)\\&
-
\dfrac{\sqrt{A}\beta}{(\sqrt{A}-\beta)^2}A_{x^ky^\ell} -\dfrac{\sqrt{A}}{(\sqrt{A}-\beta)^2}A_{x^k}\beta_\ell -\dfrac{\beta}{2\sqrt{A}(\sqrt{A}-\beta)^2}A_{x^k}A_\ell\\& +\dfrac{2\sqrt{A}\beta}{(\sqrt{A}-\beta)^3}A_{x^k}\left(\dfrac{1}{2\sqrt{A}}A_\ell-\beta_\ell \right)
+
\dfrac{2A\beta}{(\sqrt{A}-\beta)^2}\beta_{x^ky^\ell}+\dfrac{2A}{(\sqrt{A}-\beta)^2}\beta_{x^k}\beta_\ell+\dfrac{2\beta}{(\sqrt{A}-\beta)^2}\beta_{x^k}A_\ell\\& -\dfrac{4A\beta}{(\sqrt{A}-\beta)^3}\beta_{x^k}\left(\dfrac{1}{2\sqrt{A}}A_\ell-\beta_\ell \right)
+2\beta \beta_{x^ky^\ell}+2\beta_{x^k}\beta_\ell.
\end{split}
\end{equation}
Contracting (\ref{64.3}) with $y^k$, we get
\begin{align*}
L_{x^k y^\ell}y^k
=&
\dfrac{2A}{(\sqrt{A}-\beta)^2} A_{0\ell}+\dfrac{2}{(\sqrt{A}-\beta)^2} A_{0}A_\ell -\dfrac{2\sqrt{A}}{(\sqrt{A}-\beta)^3} A_{0}A_\ell+\dfrac{4A}{(\sqrt{A}-\beta)^3} A_{0}\beta_\ell \\&
-
\dfrac{A^{3/2}}{(\sqrt{A}-\beta)^3}A_{0\ell} -\dfrac{3A^{1/2}}{2(\sqrt{A}-\beta)^3}A_{0}A_\ell +\dfrac{3A}{2(\sqrt{A}-\beta)^4}A_{0}A_\ell-\dfrac{3A^{3/2}}{(\sqrt{A}-\beta)^4}A_{0}\beta_\ell\\& 
+
\dfrac{2A^2}{(\sqrt{A}-\beta)^3}\beta_{0\ell}+\dfrac{4A}{(\sqrt{A}-\beta)^3}\beta_{0}A_\ell -\dfrac{3A^{3/2}}{(\sqrt{A}-\beta)^4}\beta_{0}A_\ell+\dfrac{6A^2}{(\sqrt{A}-\beta)^4}\beta_{0}\beta_\ell\\& 
+
\dfrac{2\beta}{\sqrt{A}-\beta}A_{0\ell}+\dfrac{2}{\sqrt{A}-\beta}A_{0}\beta_\ell -\dfrac{\beta}{\sqrt{A}(\sqrt{A}-\beta)^2}A_{0}A_\ell+\dfrac{2\beta}{(\sqrt{A}-\beta)^2}A_{0}\beta_\ell
\end{align*}
\begin{align*}
&\ \ \ \ \ 
+
\dfrac{2A}{\sqrt{A}-\beta}\beta_{0\ell}+\dfrac{2}{\sqrt{A}-\beta}\beta_{0}A_\ell
-\dfrac{\sqrt{A}}{(\sqrt{A}-\beta)^2}\beta_{0}A_\ell+\dfrac{2A}{(\sqrt{A}-\beta)^2}\beta_{0}\beta_\ell\\
&\ \ \ \ \ -
\dfrac{\sqrt{A}\beta}{(\sqrt{A}-\beta)^2}A_{0\ell} -\dfrac{\sqrt{A}}{(\sqrt{A}-\beta)^2}A_{0}\beta_\ell -\dfrac{\beta}{2\sqrt{A}(\sqrt{A}-\beta)^2}A_{0}A_\ell\\& \ \ \ \ \ +\dfrac{\beta}{(\sqrt{A}-\beta)^3}A_{0}A_\ell-\dfrac{2\sqrt{A}\beta}{(\sqrt{A}-\beta)^3}A_{0}\beta_\ell
+
\dfrac{2A\beta}{(\sqrt{A}-\beta)^2}\beta_{0\ell}+\dfrac{2A}{(\sqrt{A}-\beta)^2}\beta_{0}\beta_\ell+\dfrac{2\beta}{(\sqrt{A}-\beta)^2}\beta_{0}A_\ell\\&\ \ \ \ \ -\dfrac{2\sqrt{A}\beta}{(\sqrt{A}-\beta)^3}\beta_{0}A_\ell+\dfrac{4A\beta}{(\sqrt{A}-\beta)^3}\beta_{0}\beta_\ell
+
2\beta \beta_{0\ell}+2\beta_{0}\beta_\ell.
\end{align*}
After simplification,we get
\begin{align*}
L_{x^k y^\ell} y^k
&=\dfrac{1}{2\sqrt{A}(\sqrt{A}-\beta)^4}
\bigg[8A^3\beta_{0\ell}
+2A^{5/2}\bigg\{\left(A_{0\ell}+10\beta_0 \beta_\ell+2\beta_{0\ell}\right)-4\beta\beta_{0\ell} \bigg\}\\&
+4A^2\bigg\{\left(A_0\beta_\ell+\beta_0 A_\ell\right)-\beta\left(A_{0\ell}+4\beta_0 \beta_\ell+2\beta_{0\ell}\right)-3\beta^2\beta_{0\ell}\bigg\}\\&
-4A^{3/2}\bigg\{4\beta\left(A_0\beta_\ell+\beta_0 A_\ell\right)+\beta^2\left(A_{0\ell}-5\beta_0 \beta_\ell-\beta_{0\ell}\right)-6\beta^3\beta_{0\ell}\bigg\}\\&
-2A\bigg\{\beta A_0A_\ell-3\beta^2\left(A_0\beta_\ell+\beta_0 A_\ell\right) -\beta^3\left(5A_{0\ell}-8\beta_0 \beta_\ell\right)+8\beta^4\beta_{0\ell}\bigg\}\\&
+4A^{1/2}\bigg\{2\beta^2A_0A_\ell+\beta^4 \left(\beta_0\beta_\ell-A_{0\ell}\right)+\beta^5\beta_{0\ell}\bigg\}
-3\beta^3 A_0A_\ell\bigg].
\end{align*}
 Further, equation (\ref{64.2}) can be rewritten as 
 \begin{align*}
 2L_{x^\ell}
 =&\dfrac{4A}{(\sqrt{A}-\beta)^2} A_{x^\ell} -\dfrac{2A^{3/2}}{(\sqrt{A}-\beta)^3}A_{x^\ell}+\dfrac{4A^2}{(\sqrt{A}-\beta)^3}\beta_{x^\ell}
 +\dfrac{4\beta}{\sqrt{A}-\beta}A_{x^\ell}+\dfrac{4A}{\sqrt{A}-\beta} \beta_{x^\ell}\\& -\dfrac{2\sqrt{A}\beta}{(\sqrt{A}-\beta)^2}A_{x^\ell}+\dfrac{4A\beta}{(\sqrt{A}-\beta)^2}\beta_{x^\ell}+ 4\beta \beta_{x^\ell}\\
 =&\dfrac{1}{2\sqrt{A}(\sqrt{A}-\beta)^4}
 \bigg[16A^3\beta_{x^\ell}
 +2A^{5/2}\bigg\{2A_{x^\ell}-8\beta\beta_{x^\ell}\bigg\}
 -4A^2\bigg\{2\beta A_{x^\ell}+6\beta^2\beta_{x^\ell}\bigg\}\\&
  -4A^{3/2}\bigg\{2\beta^2 A_{x^\ell}-12\beta^3\beta_{x^\ell}\bigg\}
  +2A\bigg\{10\beta^3 A_{x^\ell}-16\beta^4\beta_{x^\ell}\bigg\}
  -4A^{1/2}\bigg\{2\beta^4 A_{x^\ell}-2\beta^5\beta_{x^\ell}\bigg\}
 \bigg].
 \end{align*}
 We know that $ \bar{F}$ is locally dually flat if and only if $ L_{x^ky^\ell}y^k-2L_{x^\ell}=0, $\\
 i.e., 
 \begin{align*}
& \dfrac{1}{2\sqrt{A}(\sqrt{A}-\beta)^4}
 \bigg[8A^3\beta_{0\ell}
 +2A^{5/2}\bigg\{\left(A_{0\ell}+10\beta_0 \beta_\ell+2\beta_{0\ell}\right)-4\beta\beta_{0\ell} \bigg\}
 +4A^2\bigg\{\left(A_0\beta_\ell+\beta_0 A_\ell\right)\\&-\beta\left(A_{0\ell}+4\beta_0 \beta_\ell+2\beta_{0\ell}\right)-3\beta^2\beta_{0\ell}\bigg\}
 -4A^{3/2}\bigg\{4\beta\left(A_0\beta_\ell+\beta_0 A_\ell\right)+\beta^2\left(A_{0\ell}-5\beta_0 \beta_\ell-\beta_{0\ell}\right)-6\beta^3\beta_{0\ell}\bigg\}\\&
 -2A\bigg\{\beta A_0A_\ell-3\beta^2\left(A_0\beta_\ell+\beta_0 A_\ell\right) -\beta^3\left(5A_{0\ell}-8\beta_0 \beta_\ell\right)+8\beta^4\beta_{0\ell}\bigg\}
 +4A^{1/2}\bigg\{2\beta^2A_0A_\ell\\&+\beta^4 \left(\beta_0\beta_\ell-A_{0\ell}\right)+\beta^5\beta_{0\ell}\bigg\}
 -3\beta^3 A_0A_\ell\bigg]
 -
\dfrac{1}{2\sqrt{A}(\sqrt{A}-\beta)^4}
 \bigg[16A^3\beta_{x^\ell}
 +2A^{5/2}\bigg\{2A_{x^\ell}-8\beta\beta_{x^\ell}\bigg\}\\&
 -4A^2\bigg\{2\beta A_{x^\ell}+6\beta^2\beta_{x^\ell}\bigg\}
 -4A^{3/2}\bigg\{2\beta^2 A_{x^\ell}-12\beta^3\beta_{x^\ell}\bigg\}
 +2A\bigg\{10\beta^3 A_{x^\ell}-16\beta^4\beta_{x^\ell}\bigg\}\\&
 -4A^{1/2}\bigg\{2\beta^4 A_{x^\ell}-2\beta^5\beta_{x^\ell}\bigg\}
 \bigg]=0.
 \end{align*}
Simplifying, we get
\begin{align*}
&8A^3\bigg\{\beta_{0\ell}-2\beta_{x^\ell}\bigg\}
+2A^{5/2}\bigg\{\left(A_{0\ell}-2A_{x^\ell}+10\beta_0 \beta_\ell+2\beta_{0\ell}\right)-4\beta\left( \beta_{0\ell}-2\beta_{x^\ell}\right) \bigg\}\\&
+4A^2\bigg\{\left(A_0\beta_\ell+\beta_0 A_\ell\right)-\beta\left(A_{0\ell}-2A_{x^\ell}+4\beta_0 \beta_\ell+2\beta_{0\ell}\right)-3\beta^2\left( \beta_{0\ell}-2\beta_{x^\ell}\right) \bigg\}\\&
-4A^{3/2}\bigg\{4\beta\left(A_0\beta_\ell+\beta_0 A_\ell\right) +\beta^2\left(A_{0\ell}-2A_{x^\ell}-5\beta_0\beta_\ell-\beta_{0\ell}\right)-6\beta^3\left( \beta_{0\ell}-2\beta_{x^\ell}\right) \bigg\}\\&
-2A\bigg\{\beta A_0A_\ell-3\beta^2\left(A_0\beta_\ell+\beta_0 A_\ell\right) -\beta^3\left(5A_{0\ell}-10A_{x^\ell}-8\beta_0 \beta_\ell\right) +8\beta^4\left( \beta_{0\ell}-2\beta_{x^\ell}\right) \bigg\}\\&
+4A^{1/2}\bigg\{2\beta^2A_0A_\ell +\beta^4 \left(\beta_0\beta_\ell-A_{0\ell}+2A_{x^\ell}\right) +\beta^5\left( \beta_{0\ell} -2\beta_{x^\ell}\right) \bigg\}
-3\beta^3 A_0A_\ell=0.
\end{align*}
From the above equation, we conclude that $ \bar{F}$ is locally dually flat if and only if following seven equations are satisfied.
\begin{equation}{\label{LDF4.1}}
\beta_{0\ell}=2\beta_{x^\ell}
\end{equation} 
\begin{equation}{\label{LDF4.2}}
\left(A_{0\ell}-2A_{x^\ell}+10\beta_0 \beta_\ell+2\beta_{0\ell}\right)-4\beta\left( \beta_{0\ell}-2\beta_{x^\ell}\right)=0
\end{equation}  
\begin{equation}{\label{LDF4.3}}
\left(A_0\beta_\ell+\beta_0 A_\ell\right)-\beta\left(A_{0\ell}-2A_{x^\ell}+4\beta_0 \beta_\ell+2\beta_{0\ell}\right)-3\beta^2\left( \beta_{0\ell}-2\beta_{x^\ell}\right)=0
\end{equation}  
\begin{equation}{\label{LDF4.4}}
4\beta\left(A_0\beta_\ell+\beta_0 A_\ell\right) +\beta^2\left(A_{0\ell}-2A_{x^\ell}-5\beta_0\beta_\ell-\beta_{0\ell}\right)-6\beta^3\left( \beta_{0\ell}-2\beta_{x^\ell}\right)=0
\end{equation}  
\begin{equation}{\label{LDF4.5}}
\beta A_0A_\ell-3\beta^2\left(A_0\beta_\ell+\beta_0 A_\ell\right) -\beta^3\left(5A_{0\ell}-10A_{x^\ell}-8\beta_0 \beta_\ell\right) +8\beta^4\left( \beta_{0\ell}-2\beta_{x^\ell}\right)=0
\end{equation}  
\begin{equation}{\label{LDF4.6}}
2\beta^2A_0A_\ell +\beta^4 \left(\beta_0\beta_\ell-A_{0\ell}+2A_{x^\ell}\right) +\beta^5\left( \beta_{0\ell} -2\beta_{x^\ell}\right)=0
\end{equation}  
\begin{equation}{\label{LDF4.7}}
A_0A_\ell=0.
\end{equation}
Solving above seven equations, we get
	$$ A_{0}A_\ell=0,\ A_{0 \ell}=2A_{x^\ell},\ \beta_{0}\beta_\ell=0,\ \beta_{0 \ell}=0=2\beta_{x^\ell},\ A_\ell\beta_{0}+A_{0}\beta_\ell=0.  $$
Above discussion  leads to the following theorem.
\begin{theorem}
Let $(M, \bar{F})$ be an $n-$dimensional Finsler space with $ \bar{F}= \dfrac{F^{2}}{F-\beta} + \beta$ as a Matsumoto-Randers changed metric. Then $\bar{F}$ is locally dually flat if and only if the following equations are satisfied:
	$$ A_{0}A_\ell=0,\ A_{0 \ell}=2A_{x^\ell},\ \beta_{0}\beta_\ell=0,\ \beta_{0 \ell}=0=2\beta_{x^\ell},\ A_\ell\beta_{0}+A_{0}\beta_\ell=0.  $$
\end{theorem}
Next, we find necessary and sufficient conditions for exponential-Randers changed Finsler metric
\begin{align*}
\bar{F}= F e^{ \beta/F} + \beta
\end{align*}
to be locally dually flat.\\
Let us put $ F^2=A$ in $\bar{F},$ then
\begin{equation*}
\bar{F}=\sqrt{A}e^{\beta/\sqrt{A}}+\beta.
\end{equation*}
\begin{equation}{\label{65.1}}
L=\bar{F}^2=Ae^{2\beta/\sqrt{A}}+2\sqrt{A}\beta e^{\beta/\sqrt{A}}+\beta^2.
\end{equation}
Differentiating (\ref{65.1}) w.r.t. $x^k,$ we get
\begin{equation}{\label{65.2}}
\begin{split}
L_{x^k}
=&Ae^{2\beta/\sqrt{A}}\left(\dfrac{2}{\sqrt{A}}\beta_{x^k}-\dfrac{\beta}{A^{3/2}}A_{x^k}\right) +e^{2\beta/\sqrt{A}}A_{x^k}
+\dfrac{\beta}{\sqrt{A}} e^{\beta/\sqrt{A}}A_{x^k} +2\sqrt{A} e^{\beta/\sqrt{A}}\beta_{x^k}\\& +2\sqrt{A}\beta e^{\beta/\sqrt{A}}\left(\dfrac{1}{\sqrt{A}}\beta_{x^k}-\dfrac{\beta}{2A^{3/2}}A_{x^k}\right)
+2\beta \beta_{x^k}\\
=&2\sqrt{A}e^{2\beta/\sqrt{A}}\beta_{x^k} -\dfrac{\beta}{\sqrt{A}} e^{2\beta/\sqrt{A}}A_{x^k} +e^{2\beta/\sqrt{A}}A_{x^k}
+\dfrac{\beta}{\sqrt{A}} e^{\beta/\sqrt{A}}A_{x^k} +2\sqrt{A} e^{\beta/\sqrt{A}}\beta_{x^k}\\
&+2\beta e^{\beta/\sqrt{A}}\beta_{x^k} -\dfrac{\beta^2}{A}e^{\beta/\sqrt{A}}A_{x^k}
+2\beta \beta_{x^k}.
\end{split}
\end{equation}
Differentiation of (\ref{65.2}) further w.r.t. $y^\ell$ gives
\begin{equation}{\label{65.3}}
\begin{split}
L_{x^k y^\ell}
=&2\sqrt{A}e^{2\beta/\sqrt{A}}\beta_{x^ky^\ell} +\dfrac{1}{\sqrt{A}}e^{2\beta/\sqrt{A}}\beta_{x^k}A_\ell +2\sqrt{A}e^{2\beta/\sqrt{A}}\beta_{x^k}\left(\dfrac{2}{\sqrt{A}}\beta_{\ell}-\dfrac{\beta}{A^{3/2}}A_{\ell}\right)
-\dfrac{\beta}{\sqrt{A}} e^{2\beta/\sqrt{A}}A_{x^ky^\ell}\\& -\dfrac{1}{\sqrt{A}} e^{2\beta/\sqrt{A}}A_{x^k}\beta_\ell +\dfrac{\beta}{2A^{3/2}} e^{2\beta/\sqrt{A}}A_{x^k}A_\ell -\dfrac{\beta}{\sqrt{A}} e^{2\beta/\sqrt{A}}A_{x^k}\left(\dfrac{2}{\sqrt{A}}\beta_{\ell} -\dfrac{\beta}{A^{3/2}}A_{\ell}\right)\\
&+e^{2\beta/\sqrt{A}}A_{x^ky^\ell} +e^{2\beta/\sqrt{A}}A_{x^k}\left(\dfrac{2}{\sqrt{A}}\beta_{\ell}-\dfrac{\beta}{A^{3/2}}A_{\ell}\right)
+\dfrac{\beta}{\sqrt{A}} e^{\beta/\sqrt{A}}A_{x^ky^\ell} +\dfrac{1}{\sqrt{A}} e^{\beta/\sqrt{A}}A_{x^k}\beta_\ell \\&-\dfrac{\beta}{2A^{3/2}} e^{\beta/\sqrt{A}}A_{x^k}A_\ell +\dfrac{\beta}{\sqrt{A}} e^{\beta/\sqrt{A}}A_{x^k} \left(\dfrac{1}{\sqrt{A}}\beta_{\ell} -\dfrac{\beta}{2A^{3/2}}A_{\ell}\right)
+2\sqrt{A} e^{\beta/\sqrt{A}}\beta_{x^ky^\ell}\\& +\dfrac{1}{\sqrt{A}} e^{\beta/\sqrt{A}}\beta_{x^k}A_\ell +2\sqrt{A} e^{\beta/\sqrt{A}}\beta_{x^k}\left(\dfrac{1}{\sqrt{A}}\beta_{\ell} -\dfrac{\beta}{2A^{3/2}}A_{\ell}\right)
+2\beta e^{\beta/\sqrt{A}}\beta_{x^ky^\ell} +2 e^{\beta/\sqrt{A}}\beta_{x^k}\beta_\ell\\& +2\beta e^{\beta/\sqrt{A}}\beta_{x^k}\left(\dfrac{1}{\sqrt{A}}\beta_{\ell} -\dfrac{\beta}{2A^{3/2}}A_{\ell}\right)
-\dfrac{\beta^2}{A}e^{\beta/\sqrt{A}}A_{x^ky^\ell} -\dfrac{2\beta}{A}e^{\beta/\sqrt{A}}A_{x^k}\beta_\ell +\dfrac{\beta^2}{A^2}e^{\beta/\sqrt{A}}A_{x^k}A_\ell\\& -\dfrac{\beta^2}{A}e^{\beta/\sqrt{A}}A_{x^k}\left(\dfrac{1}{\sqrt{A}}\beta_{\ell} -\dfrac{\beta}{2A^{3/2}}A_{\ell}\right)
+2\beta \beta_{x^ky^\ell}+2 \beta_{x^k}\beta_\ell.
\end{split}
\end{equation}
Contracting (\ref{65.3}) with $y^k$, we get
\begin{align*}
L_{x^k y^\ell}y^k
=&2\sqrt{A}e^{2\beta/\sqrt{A}}\beta_{0\ell} +\dfrac{1}{\sqrt{A}}e^{2\beta/\sqrt{A}}\beta_{0}A_\ell +4e^{2\beta/\sqrt{A}}\beta_{0}\beta_{\ell} -\dfrac{2\beta}{A}e^{2\beta/\sqrt{A}}\beta_{0}A_{\ell}
-\dfrac{\beta}{\sqrt{A}} e^{2\beta/\sqrt{A}}A_{0\ell}\\
&-\dfrac{1}{\sqrt{A}} e^{2\beta/\sqrt{A}}A_{0}\beta_\ell +\dfrac{\beta}{2A^{3/2}} e^{2\beta/\sqrt{A}}A_{0}A_\ell
-\dfrac{2\beta}{A} e^{2\beta/\sqrt{A}}A_{0} \beta_{\ell} +\dfrac{\beta^2}{A^2} e^{2\beta/\sqrt{A}}A_{0}A_{\ell}\\
&+e^{2\beta/\sqrt{A}}A_{0\ell} +\dfrac{2}{\sqrt{A}}e^{2\beta/\sqrt{A}}A_{0}\beta_{\ell} -\dfrac{\beta}{A^{3/2}}e^{2\beta/\sqrt{A}}A_{0}A_{\ell}
+\dfrac{\beta}{\sqrt{A}} e^{\beta/\sqrt{A}}A_{0\ell} +\dfrac{1}{\sqrt{A}} e^{\beta/\sqrt{A}}A_{0}\beta_\ell 
\end{align*}
\begin{align*}
& \ \ \ \ \ -\dfrac{\beta}{2A^{3/2}} e^{\beta/\sqrt{A}}A_{0}A_\ell +\dfrac{\beta}{A} e^{\beta/\sqrt{A}}A_{0}\beta_{\ell}  -\dfrac{\beta^2}{2A^2} e^{\beta/\sqrt{A}}A_{0}A_{\ell}
+2\sqrt{A} e^{\beta/\sqrt{A}}\beta_{0\ell}\\
& \ \ \ \ \ +\dfrac{1}{\sqrt{A}} e^{\beta/\sqrt{A}}\beta_{0}A_\ell +2 e^{\beta/\sqrt{A}}\beta_{0}\beta_{\ell}  -\dfrac{\beta}{A}e^{\beta/\sqrt{A}}\beta_{0}A_{\ell}
+2\beta e^{\beta/\sqrt{A}}\beta_{0\ell} +2 e^{\beta/\sqrt{A}}\beta_{0}\beta_\ell\\
& \ \ \ \ \ +\dfrac{2\beta}{\sqrt{A}} e^{\beta/\sqrt{A}}\beta_{0}\beta_{\ell} -\dfrac{\beta^2}{A^{3/2}}e^{\beta/\sqrt{A}}\beta_{0}A_{\ell}
-\dfrac{\beta^2}{A}e^{\beta/\sqrt{A}}A_{0\ell} -\dfrac{2\beta}{A}e^{\beta/\sqrt{A}}A_{0}\beta_\ell +\dfrac{\beta^2}{A^2}e^{\beta/\sqrt{A}}A_{0}A_\ell\\
& \ \ \ \ \ -\dfrac{\beta^2}{A^{3/2}}e^{\beta/\sqrt{A}}A_{0}\beta_{\ell} +\dfrac{\beta^3}{2A^{5/2}}e^{\beta/\sqrt{A}}A_{0}A_{\ell}
+2\beta \beta_{0\ell}+2 \beta_{0}\beta_\ell.
\end{align*}
Simplifying, we get
\begin{align*}
L_{x^k y^\ell}y^k
=&\dfrac{1}{2A^{5/2}}
\bigg[4A^3e^{\beta/\sqrt{A}}\left(e^{\beta/\sqrt{A}}+1\right)\beta_{0\ell}
+2A^{5/2}\bigg\{2\beta\left( e^{\beta/\sqrt{A}}\beta_{0\ell}+\beta_{0\ell}\right) +e^{2\beta/\sqrt{A}}A_{0\ell}\\&+4e^{\beta/\sqrt{A}}\left( e^{\beta/\sqrt{A}}+1\right) \beta_0 \beta_\ell+2\beta_0 \beta_\ell \bigg\}
+2e^{\beta/\sqrt{A}}A^2\bigg\{ \beta\left( \left(1-e^{\beta/\sqrt{A}}\right)A_{0\ell}+2\beta_0\beta_\ell \right)\\&+\left( e^{\beta/\sqrt{A}}+1\right) \left( A_0\beta_\ell+\beta_0A_\ell\right)  \bigg\}
-2e^{\beta/\sqrt{A}}A^{3/2}\bigg\{\beta^2A_{0\ell}+\beta\left(2e^{\beta/\sqrt{A}}+1\right)\left( A_0\beta_\ell+\beta_0A_\ell\right) \bigg\}\\&
-e^{\beta/\sqrt{A}}A
\bigg\{2\left( A_0\beta_\ell+\beta_0A_\ell\right)
+\beta\left(e^{\beta/\sqrt{A}}+1\right)A_0A_\ell\bigg\}
+\beta^2e^{\beta/\sqrt{A}}\left( 2e^{\beta/\sqrt{A}}+1\right)\sqrt{A} A_0A_\ell\\&
+\beta^3e^{\beta/\sqrt{A}}A_0A_\ell
\bigg].
\end{align*}  
Further, equation (\ref{65.2}) can be rewritten as  
\begin{align*}
2L_{x^\ell}
=&4\sqrt{A}e^{2\beta/\sqrt{A}}\beta_{x^\ell} -\dfrac{2\beta}{\sqrt{A}} e^{2\beta/\sqrt{A}}A_{x^\ell} +2e^{2\beta/\sqrt{A}}A_{x^\ell}
+\dfrac{2\beta}{\sqrt{A}} e^{\beta/\sqrt{A}}A_{x^\ell} +4\sqrt{A} e^{\beta/\sqrt{A}}\beta_{x^\ell}\\
&+4\beta e^{\beta/\sqrt{A}}\beta_{x^\ell} -\dfrac{2\beta^2}{A}e^{\beta/\sqrt{A}}A_{x^\ell}
+4\beta \beta_{x^\ell}\\
=&4\sqrt{A}e^{\beta/\sqrt{A}}\left(e^{\beta/\sqrt{A}}+1\right)\beta_{x^\ell} +2e^{2\beta/\sqrt{A}}A_{x^\ell} +4\beta\left(e^{\beta/\sqrt{A}}+1\right)\beta_{x^\ell}\\& -\dfrac{2\beta}{\sqrt{A}}e^{\beta/\sqrt{A}}\left(e^{\beta/\sqrt{A}}-1\right)A_{x^\ell}-\dfrac{2\beta^2}{A}e^{\beta/\sqrt{A}}A_{x^\ell} \\
=&\dfrac{1}{2A^{5/2}}
\bigg[8A^3e^{\beta/\sqrt{A}}\left(e^{\beta/\sqrt{A}}+1\right)\beta_{x^\ell}
+4A^{5/2}\bigg\{e^{2\beta/\sqrt{A}}A_{x^\ell}+2\beta\left( e^{\beta/\sqrt{A}}+1\right) \beta_{x^\ell} \bigg\}\\&
-4\beta A^2 e^{\beta/\sqrt{A}}\left( e^{\beta/\sqrt{A}}-1\right) A_{x^\ell}
-4\beta^2A^{3/2}e^{\beta/\sqrt{A}} A_{x^\ell}
\bigg].
\end{align*}
We know that $ \bar{F}$ is locally dually flat if and only if $ L_{x^ky^\ell}y^k-2L_{x^\ell}=0, $\\
i.e., 
\begin{align*}
&\dfrac{1}{2A^{5/2}}
\bigg[4A^3e^{\beta/\sqrt{A}}\left(e^{\beta/\sqrt{A}}+1\right)\beta_{0\ell}
+2A^{5/2}\bigg\{2\beta\left( e^{\beta/\sqrt{A}}\beta_{0\ell}+\beta_{0\ell}\right) +e^{2\beta/\sqrt{A}}A_{0\ell}\\&+4e^{\beta/\sqrt{A}}\left( e^{\beta/\sqrt{A}}+1\right) \beta_0 \beta_\ell+2\beta_0 \beta_\ell \bigg\}
+2e^{\beta/\sqrt{A}}A^2\bigg\{ \beta\left( \left(1-e^{\beta/\sqrt{A}}\right)A_{0\ell}+2\beta_0\beta_\ell \right)\\&+\left( e^{\beta/\sqrt{A}}+1\right) \left( A_0\beta_\ell+\beta_0A_\ell\right)  \bigg\}
-2e^{\beta/\sqrt{A}}A^{3/2}\bigg\{\beta^2A_{0\ell}+\beta\left(2e^{\beta/\sqrt{A}}+1\right)\left( A_0\beta_\ell+\beta_0A_\ell\right) \bigg\}\\&
-e^{\beta/\sqrt{A}}A
\bigg\{2\left( A_0\beta_\ell+\beta_0A_\ell\right)
+\beta\left(e^{\beta/\sqrt{A}}+1\right)A_0A_\ell\bigg\}
+\beta^2e^{\beta/\sqrt{A}}\left( 2e^{\beta/\sqrt{A}}+1\right)\sqrt{A} A_0A_\ell\\&
+\beta^3e^{\beta/\sqrt{A}}A_0A_\ell
\bigg]
-
\dfrac{1}{2A^{5/2}}
\bigg[8A^3e^{\beta/\sqrt{A}}\left(e^{\beta/\sqrt{A}}+1\right)\beta_{x^\ell}
+4A^{5/2}\bigg\{e^{2\beta/\sqrt{A}}A_{x^\ell}+2\beta\left( e^{\beta/\sqrt{A}}+1\right) \beta_{x^\ell} \bigg\}\\&
-4\beta A^2 e^{\beta/\sqrt{A}}\left( e^{\beta/\sqrt{A}}-1\right) A_{x^\ell}
-4\beta^2A^{3/2}e^{\beta/\sqrt{A}} A_{x^\ell}
\bigg]=0.
\end{align*} 
Simplifying, we get
\begin{align*}
&4A^3e^{\beta/\sqrt{A}}\left(e^{\beta/\sqrt{A}}+1\right)\bigg\{\beta_{0\ell} -2\beta_{x^\ell}\bigg\}\\&
+2A^{5/2}\bigg\{2\beta\left( e^{\beta/\sqrt{A}}+1\right) \left(\beta_{0\ell}-2\beta_{x^\ell}\right)  +e^{2\beta/\sqrt{A}}\left(A_{0\ell}-2A_{x^\ell}\right)+4e^{\beta/\sqrt{A}}\left( e^{\beta/\sqrt{A}}+1\right) \beta_0 \beta_\ell+2\beta_0 \beta_\ell \bigg\}\\&
+2e^{\beta/\sqrt{A}}A^2\bigg\{\beta\left( \left(1-e^{\beta/\sqrt{A}}\right)\left(A_{0\ell}-2A_{x^\ell}\right) +2\beta_0\beta_\ell \right)+\left( e^{\beta/\sqrt{A}}+1\right) \left( A_0\beta_\ell+\beta_0A_\ell\right)\bigg\}\\&
-2e^{\beta/\sqrt{A}}A^{3/2}\bigg\{\beta^2\left(A_{0\ell}-2A_{x^\ell}\right) +\beta\left(2e^{\beta/\sqrt{A}}+1\right)\left( A_0\beta_\ell+\beta_0A_\ell\right) \bigg\}\\&
-e^{\beta/\sqrt{A}}A
\bigg\{2\left( A_0\beta_\ell+\beta_0A_\ell\right)
+\beta\left(e^{\beta/\sqrt{A}}+1\right)A_0A_\ell\bigg\}\\&
+\beta^2e^{\beta/\sqrt{A}}\left( 2e^{\beta/\sqrt{A}}+1\right)\sqrt{A} A_0A_\ell
+\beta^3e^{\beta/\sqrt{A}}A_0A_\ell=0.
\end{align*} 
From the above equation, we conclude that $ \bar{F}$ is locally dually flat if and only if following seven equations are satisfied.
\begin{equation}{\label{LDF5.1}}
\left(e^{\beta/\sqrt{A}}+1\right)\bigg\{\beta_{0\ell} -2\beta_{x^\ell}\bigg\}=0\ \implies\ \beta_{0\ell}=2\beta_{x^\ell}
\end{equation} 
\begin{equation}{\label{LDF5.2}}
2\beta\left( e^{\beta/\sqrt{A}}+1\right) \left(\beta_{0\ell}-2\beta_{x^\ell}\right)  +e^{2\beta/\sqrt{A}}\left(A_{0\ell}-2A_{x^\ell}\right)+4e^{\beta/\sqrt{A}}\left( e^{\beta/\sqrt{A}}+1\right) \beta_0 \beta_\ell+2\beta_0 \beta_\ell=0
\end{equation} 
\begin{equation}{\label{LDF5.3}}
\beta\left( \left(1-e^{\beta/\sqrt{A}}\right)\left(A_{0\ell}-2A_{x^\ell}\right) +2\beta_0\beta_\ell \right)+\left( e^{\beta/\sqrt{A}}+1\right) \left( A_0\beta_\ell+\beta_0A_\ell\right)=0
\end{equation} 
\begin{equation}{\label{LDF5.4}}
\beta^2\left(A_{0\ell}-2A_{x^\ell}\right) +\beta\left(2e^{\beta/\sqrt{A}}+1\right)\left( A_0\beta_\ell+\beta_0A_\ell\right)=0
\end{equation} 
\begin{equation}{\label{LDF5.5}}
2\left( A_0\beta_\ell+\beta_0A_\ell\right)
+\beta\left(e^{\beta/\sqrt{A}}+1\right)A_0A_\ell=0
\end{equation} 
\begin{equation}{\label{LDF5.6}}
\left( 2e^{\beta/\sqrt{A}}+1\right) A_0A_\ell=0\ \implies\ A_0A_\ell=0.
\end{equation} 
Further, from the equations (\ref{LDF5.5}) and (\ref{LDF5.6}), we get
\begin{equation}{\label{LDF5.8}}
 A_0\beta_\ell+\beta_0A_\ell=0
\end{equation} 
Again, from the equations (\ref{LDF5.4}) and (\ref{LDF5.8}), we get
\begin{equation}{\label{LDF5.9}}
A_{0\ell}=2A_{x^\ell},
\end{equation} 
and from the equations (\ref{LDF5.3}), (\ref{LDF5.8}) and (\ref{LDF5.9}), we get
\begin{equation}{\label{LDF5.10}}
\beta_0\beta_\ell=0.
\end{equation} 
Above discussion  leads to the following theorem.
\begin{theorem}
Let $(M, \bar{F})$ be an $n-$dimensional Finsler space with $\bar{F}
= F e^{ \beta/F} + \beta$ as an exponential-Randers changed metric. Then $\bar{F}$ is locally dually flat if and only if the following equations are satisfied:
	$$ A_{0}A_\ell=0,\ A_{0 \ell}=2A_{x^\ell},\ \beta_{0}\beta_\ell=0,\ \beta_{0 \ell}=2\beta_{x^\ell},\ A_\ell\beta_{0}+A_{0}\beta_\ell=0.  $$
\end{theorem}
Next, we find necessary and sufficient conditions for infinite series-Randers changed  Finsler metric 
$$ \bar{F}=\dfrac{\beta^{2}}{\beta - F } + \beta $$ 
to be locally dually flat.\\
Let us put $F^2=A$ in $\bar{F},$ then
\begin{align*}
\bar{F}=\dfrac{\beta^{2}}{\beta -\sqrt{A}} + \beta.
\end{align*}
\begin{equation}{\label{66.1}}
 L=\bar{F}^2= \dfrac{\beta^{4}}{(\beta -\sqrt{A})^2} +\dfrac{2\beta^{3}}{\beta -\sqrt{A}} + \beta^2.
\end{equation}
Differentiating (\ref{66.1}) w.r.t. $x^k,$ we get
\begin{equation}{\label{66.2}}
\begin{split}
L_{x^k}
=&\dfrac{4\beta^{3}}{(\beta -\sqrt{A})^2}\beta_{x^k}-\dfrac{2\beta^{4}}{(\beta -\sqrt{A})^3}\left(\beta_{x^k}-\dfrac{1}{2\sqrt{A}}A_{x^k} \right)   
+\dfrac{6\beta^{2}}{\beta -\sqrt{A}}\beta_{x^k} -\dfrac{2\beta^{3}}{(\beta -\sqrt{A})^2}\left(\beta_{x^k}-\dfrac{1}{2\sqrt{A}}A_{x^k} \right)  \\&  
+ 2\beta\beta_{x^k}\\
=&-\dfrac{2\beta^{3}\sqrt{A}}{(\beta -\sqrt{A})^3}\beta_{x^k}
+\dfrac{6\beta^{2}}{\beta -\sqrt{A}}\beta_{x^k}+2\beta\beta_{x^k}
+\dfrac{\beta^{4}}{\sqrt{A}(\beta -\sqrt{A})^3}A_{x^k}
+\dfrac{\beta^{3}}{\sqrt{A}(\beta -\sqrt{A})^2}A_{x^k}.
\end{split}
\end{equation}
Differentiation of (\ref{66.2}) further w.r.t. $y^\ell$  gives
\begin{equation}{\label{66.3}}
\begin{split}
L_{x^ky^\ell}
=&-\dfrac{2\beta^{3}\sqrt{A}}{(\beta -\sqrt{A})^3}\beta_{x^ky^\ell}
-\dfrac{6\beta^{2}\sqrt{A}}{(\beta -\sqrt{A})^3}\beta_{x^k}\beta_\ell
-\dfrac{\beta^{3}}{\sqrt{A}(\beta -\sqrt{A})^3}\beta_{x^k}A_\ell
+\dfrac{6\beta^{3}\sqrt{A}}{(\beta -\sqrt{A})^4}\beta_{x^k}\left(\beta_\ell-\dfrac{A_\ell}{2\sqrt{A}} \right)\\&
+
\dfrac{6\beta^{2}}{\beta -\sqrt{A}}\beta_{x^ky^\ell}
+\dfrac{12\beta}{\beta -\sqrt{A}}\beta_{x^k}\beta_\ell
-\dfrac{6\beta^{2}}{(\beta -\sqrt{A})^2}\beta_{x^k}
\left(\beta_\ell-\dfrac{1}{2\sqrt{A}}A_\ell \right)
+
2\beta\beta_{x^ky^\ell}+2\beta_{x^k}\beta_\ell\\&
+
\dfrac{\beta^{4}}{\sqrt{A}(\beta -\sqrt{A})^3}A_{x^ky^\ell}
+\dfrac{4\beta^{3}}{\sqrt{A}(\beta -\sqrt{A})^3}A_{x^k}\beta_\ell
-\dfrac{\beta^{4}}{2{A}^{3/2}(\beta -\sqrt{A})^3}A_{x^k}A_\ell\\&
-\dfrac{3\beta^{4}}{\sqrt{A}(\beta -\sqrt{A})^4}A_{x^k}
\left(\beta_\ell-\dfrac{1}{2\sqrt{A}}A_\ell \right)
+
\dfrac{\beta^{3}}{\sqrt{A}(\beta -\sqrt{A})^2}A_{x^ky^\ell}
+\dfrac{3\beta^{2}}{\sqrt{A}(\beta -\sqrt{A})^2}A_{x^k}\beta_\ell\\&
-\dfrac{\beta^{3}}{2A^{3/2}(\beta -\sqrt{A})^2}A_{x^k}A_\ell
-\dfrac{2\beta^{3}}{\sqrt{A}(\beta -\sqrt{A})^3}A_{x^k}
\left(\beta_\ell-\dfrac{1}{2\sqrt{A}}A_\ell \right).
\end{split}
\end{equation}
Contracting (\ref{66.3}) with $y^k$, we get
\begin{align*}
L_{x^k y^\ell}y^k
=&-\dfrac{2\beta^{3}\sqrt{A}}{(\beta -\sqrt{A})^3}\beta_{0\ell}
-\dfrac{6\beta^{2}\sqrt{A}}{(\beta -\sqrt{A})^3}\beta_{0}\beta_\ell
-\dfrac{\beta^{3}}{\sqrt{A}(\beta -\sqrt{A})^3}\beta_{0}A_\ell
+\dfrac{6\beta^{3}\sqrt{A}}{(\beta -\sqrt{A})^4}\beta_{0} \left(\beta_\ell-\dfrac{1}{2\sqrt{A}}A_\ell \right)\\&
+
\dfrac{6\beta^{2}}{\beta -\sqrt{A}}\beta_{0\ell}
+\dfrac{12\beta}{\beta -\sqrt{A}}\beta_{0}\beta_\ell
-\dfrac{6\beta^{2}}{(\beta -\sqrt{A})^2}\beta_{0}
\left(\beta_\ell-\dfrac{1}{2\sqrt{A}}A_\ell \right)
+
2\beta\beta_{0\ell}+2\beta_0\beta_{\ell}\\&
+
\dfrac{\beta^{4}}{\sqrt{A}(\beta -\sqrt{A})^3}A_{0\ell}
+\dfrac{4\beta^{3}}{\sqrt{A}(\beta -\sqrt{A})^3}A_{0}\beta_\ell
-\dfrac{\beta^{4}}{2{A}^{3/2}(\beta -\sqrt{A})^3}A_{0}A_\ell\\&
-\dfrac{3\beta^{4}}{\sqrt{A}(\beta -\sqrt{A})^4}A_{0}
\left(\beta_\ell-\dfrac{1}{2\sqrt{A}}A_\ell \right)
+
\dfrac{\beta^{3}}{\sqrt{A}(\beta -\sqrt{A})^2}A_{0\ell}
+\dfrac{3\beta^{2}}{\sqrt{A}(\beta -\sqrt{A})^2}A_{0}\beta_\ell\\&
-\dfrac{\beta^{3}}{2A^{3/2}(\beta -\sqrt{A})^2}A_{0}A_\ell
-\dfrac{2\beta^{3}}{\sqrt{A}(\beta -\sqrt{A})^3}A_{0}
\left(\beta_\ell-\dfrac{1}{2\sqrt{A}}A_\ell \right).
\end{align*}
Simplifying, we get
\begin{align*}
L_{x^k y^\ell}y^k
=&\dfrac{1}{2{A}^{3/2}(\beta -\sqrt{A})^4}
\bigg[4A^{7/2}\bigg\{ \beta\beta_{0\ell}+2\beta_{0}\beta_\ell\bigg\}
-28A^3\bigg\{\beta^2\beta_{0\ell}+2\beta\beta_{0}\beta_\ell\bigg\}\\
+&8A^{5/2}\bigg\{8\beta^3\beta_{0\ell}+15\beta^2\beta_{0}\beta_\ell\bigg\}
+2A^2\bigg\{-28\beta^4\beta_{0\ell}-40\beta^3\beta_{0}\beta_\ell+\beta^3A_{0\ell}+3\beta^2\left( A_0\beta_\ell+\beta_{0}A_\ell\right)\bigg\}\\
+&2A^{3/2}\bigg\{8\beta^5\beta_{0\ell}+10\beta^4\beta_{0}\beta_\ell-3\beta^4A_{0\ell}-8\beta^3\left( A_0\beta_\ell+\beta_{0}A_\ell\right)\bigg\}\\
+&A\bigg\{4\beta^5A_{0\ell}+4\beta^4\left( A_0\beta_\ell+\beta_{0}A_\ell\right) -3\beta^3A_0A_\ell\bigg\}
+8\beta^4\sqrt{A}A_0A_\ell
-2\beta^5A_0A_\ell
\bigg].
\end{align*}
Further,  equation (\ref{66.2}) can be rewritten as  
\begin{align*}
L_{x^\ell}
=&-\dfrac{4\beta^{3}\sqrt{A}}{(\beta -\sqrt{A})^3}\beta_{x^\ell}
+\dfrac{12\beta^{2}}{\beta -\sqrt{A}}\beta_{x^\ell}+4\beta\beta_{x^\ell}
+\dfrac{2\beta^{4}}{\sqrt{A}(\beta -\sqrt{A})^3}A_{x^\ell}
+\dfrac{2\beta^{3}}{\sqrt{A}(\beta -\sqrt{A})^2}A_{x^\ell}\\
=&\dfrac{4\beta\left(4\beta^{3}-10\beta^2\sqrt{A}+6\beta A-A^{3/2}\right) }{(\beta -\sqrt{A})^3}\beta_{x^\ell}
\dfrac{4\beta^{3}\sqrt{A}}{(\beta -\sqrt{A})^3}A_{x^\ell}\\
=&\dfrac{1}{2{A}^{3/2}(\beta -\sqrt{A})^4}
\bigg[8\beta A^{7/2}\beta_{x^\ell}
-56A^3\beta^2\beta_{x^\ell}
+128A^{5/2}\beta^3\beta_{x^\ell}
-2A^2\bigg\{56\beta^4\beta_{x^\ell}-2\beta^3A_{x^\ell}\bigg\}\\
&+2A^{3/2}\bigg\{16\beta^5\beta_{x^\ell}-6\beta^4A_{x^\ell}\bigg\}
+8\beta^5AA_{x^\ell}
\bigg].
\end{align*}
We know that $ \bar{F}$ is locally dually flat if and only if $ L_{x^ky^\ell}y^k-2L_{x^\ell}=0, $\\
i.e., 
\begin{align*}
&\dfrac{1}{2{A}^{3/2}(\beta -\sqrt{A})^4}
\bigg[4A^{7/2}\bigg\{ \beta\beta_{0\ell}+2\beta_{0}\beta_\ell\bigg\}
-28A^3\bigg\{\beta^2\beta_{0\ell}+2\beta\beta_{0}\beta_\ell\bigg\}
+8A^{5/2}\bigg\{8\beta^3\beta_{0\ell}+15\beta^2\beta_{0}\beta_\ell\bigg\}\\&
+2A^2\bigg\{-28\beta^4\beta_{0\ell}-40\beta^3\beta_{0}\beta_\ell+\beta^3A_{0\ell}+3\beta^2\left( A_0\beta_\ell+\beta_{0}A_\ell\right)\bigg\}
+2A^{3/2}\bigg\{8\beta^5\beta_{0\ell}+10\beta^4\beta_{0}\beta_\ell\\&-3\beta^4A_{0\ell}-8\beta^3\left( A_0\beta_\ell+\beta_{0}A_\ell\right)\bigg\}
+A\bigg\{4\beta^5A_{0\ell}+4\beta^4\left( A_0\beta_\ell+\beta_{0}A_\ell\right) -3\beta^3A_0A_\ell\bigg\}
+8\beta^4\sqrt{A}A_0A_\ell\\&
-2\beta^5A_0A_\ell
\bigg]
-
\dfrac{1}{2{A}^{3/2}(\beta -\sqrt{A})^4}
\bigg[8\beta A^{7/2}\beta_{x^\ell}
-56A^3\beta^2\beta_{x^\ell}
+128A^{5/2}\beta^3\beta_{x^\ell}
-2A^2\bigg\{56\beta^4\beta_{x^\ell}-2\beta^3A_{x^\ell}\bigg\}\\
&+2A^{3/2}\bigg\{16\beta^5\beta_{x^\ell}-6\beta^4A_{x^\ell}\bigg\}
+8\beta^5AA_{x^\ell}
\bigg]=0.
\end{align*}
After simplification, we get
\begin{align*}
&4A^{7/2}\bigg\{\beta\left(\beta_{0\ell}-2\beta_{x^\ell}\right) +2\beta_{0}\beta_\ell\bigg\}
-28A^3\beta\bigg\{\beta\left(\beta_{0\ell}-2\beta_{x^\ell}\right)+2\beta_{0}\beta_\ell\bigg\}\\
&+2A^2\beta^2\bigg\{-28\beta^2\left(\beta_{0\ell}-2\beta_{x^\ell}\right)-40\beta\beta_{0}\beta_\ell+\beta\left( A_{0\ell}-2A_{x^\ell}\right)+3\left( A_0\beta_\ell+\beta_{0}A_\ell\right)\bigg\}\\
&+8A^{5/2}\beta^2\bigg\{8\beta\left(\beta_{0\ell}-2\beta_{x^\ell}\right)+15\beta_{0}\beta_\ell\bigg\}\\
&+2A^{3/2}\beta^3\bigg\{8\beta^2\left(\beta_{0\ell}-2\beta_{x^\ell}\right)+10\beta\beta_{0}\beta_\ell-3\beta\left( A_{0\ell}-2A_{x^\ell}\right)-8\left( A_0\beta_\ell+\beta_{0}A_\ell\right)\bigg\}\\&
+A\beta^3\bigg\{4\beta^2\left( A_{0\ell}-2A_{x^\ell}\right)+4\beta\left( A_0\beta_\ell+\beta_{0}A_\ell\right) -3A_0A_\ell\bigg\}
+8\beta^4\sqrt{A}A_0A_\ell-2\beta^5A_0A_\ell=0.
\end{align*}
From the above equation, we conclude that $ \bar{F}$ is locally dually flat if and only if following six equations are satisfied.
\begin{equation}{\label{LDF6.1}}
\beta\left(\beta_{0\ell}-2\beta_{x^\ell}\right) +2\beta_{0}\beta_\ell=0
\end{equation} 
\begin{equation}{\label{LDF6.2}}
-28\beta^2\left(\beta_{0\ell}-2\beta_{x^\ell}\right)-40\beta\beta_{0}\beta_\ell+\beta\left( A_{0\ell}-2A_{x^\ell}\right)+3\left( A_0\beta_\ell+\beta_{0}A_\ell\right)=0
\end{equation} 
\begin{equation}{\label{LDF6.3}}
8\beta\left(\beta_{0\ell}-2\beta_{x^\ell}\right)+15\beta_{0}\beta_\ell=0
\end{equation} 
\begin{equation}{\label{LDF6.4}}
8\beta^2\left(\beta_{0\ell}-2\beta_{x^\ell}\right)+10\beta\beta_{0}\beta_\ell-3\beta\left( A_{0\ell}-2A_{x^\ell}\right)-8\left( A_0\beta_\ell+\beta_{0}A_\ell\right)=0
\end{equation} 
\begin{equation}{\label{LDF6.5}}
4\beta^2\left( A_{0\ell}-2A_{x^\ell}\right)+4\beta\left( A_0\beta_\ell+\beta_{0}A_\ell\right) -3A_0A_\ell=0
\end{equation} 
\begin{equation}{\label{LDF6.6}}
A_0A_\ell=0
\end{equation} 
Further, from the equations (\ref{LDF6.1}) and (\ref{LDF6.3}), we get  
\begin{equation}{\label{LDF6.7}}
\beta_{0\ell}=2\beta_{x^\ell}
\end{equation} 
and 
\begin{equation}{\label{LDF6.8}}
\beta_{0}\beta_\ell=0.
\end{equation} 
Again, from the equations (\ref{LDF6.7}), (\ref{LDF6.8}), (\ref{LDF6.2}) and (\ref{LDF6.4}), we get  
\begin{equation}{\label{LDF6.9}}
A_{0\ell}=2A_{x^\ell},
\end{equation} 
and 
\begin{equation}{\label{LDF6.10}}
A_0\beta_\ell+\beta_{0}A_\ell=0.
\end{equation}
Above discussion  leads to the following theorem. 
\begin{theorem}
Let $(M, \bar{F})$ be an $n-$dimensional Finsler space with $ \bar{F}=\dfrac{\beta^{2}}{\beta - F } + \beta$
	ass an infinite series-Randers changed metric. Then $ \bar{F}$  is locally dually flat if and only if the following equations are satisfied:
	$$ A_{0}A_\ell=0,\ A_{0 \ell}=2A_{x^\ell},\ \beta_{0}\beta_\ell=0,\ \beta_{0 \ell}=2\beta_{x^\ell},\ A_\ell\beta_{0}+A_{0}\beta_\ell=0.  $$
\end{theorem}


\begin{thebibliography}{99}
	
	\bibitem{S.I.Ama H.Nag} S. I. Amari and H. Nagaoka, \textit{Methods of information geometry}, Translations of Mathematical Monographs, AMS, 191, Oxford Univ. Press, 2000.
	
	\bibitem{Antonelli}	P. L. Antonelli, R. S. Ingarden, and M. Matsumoto, \textit{The Theory of Sprays and Finsler spaces with Applications in Physics and Biology}, Vol. 58, Springer Science $\&$ Business Media, 2013.
			
	\bibitem{SAB.GS.2016PFDF} S. A. Baby and G. Shanker, On Randers-Conformal Change of Finsler Space with Special  $(\alpha, \beta)-$metrics,  \textit{International Journal of Contemporary Mathematical Sciences}, \textbf{11} (9)(2016), 415-423.
	
	\bibitem{BCS} D. Bao, S. S. Chern and Z. Shen, \textit{An Introduction to Riemann-Finsler Geometry}, GTM 200, Springer-Verlag, New York, 2000.
	
	\bibitem{Chern1996} S. S. Chern, Finsler geometry is just Riemannian geometry without quadratic restriction, \textit{Notices of AMS}, \textbf{43} (9)(1996), 959-963.
	
	\bibitem{ChernShenRFG} S. S. Chern and Z. Shen, \textit{Riemann-Finsler Geometry}, Nankai Tracts in Mathematics, Vol. 6, World Scientific, Singapore, 2005.
	
	
	\bibitem{GHAM1903} G. Hamel, $\ddot{U}$ber die Geometrieen, in denen die Geraden die K$\ddot{u}$rzesten sind, \textit{Mathematische Annalen}, \textbf{57}  (2)(1903), 231–264.
	
	
	\bibitem{Hilbert4prob} D. Hilbert, Gesammelte Abhandlungen, Band III, 303–305, or, \textit{Archiv f. Math. u. Phys.3. Reihe, Bd. I}, (1901), 44–65, 213–232.
	
	
	\bibitem{M.Mat1972} M. Matsumoto, On C-reducible Finsler-spaces, \textit{Tensor, N. S.}, \textbf{24} (1972),  29-37.
	
	\bibitem{GS.SAB.2015PF}	G. Shanker and S. A. Baby, On the projective flatness of a Finsler space with infinite series $(\alpha, \beta)-$metric $\dfrac{\beta^2}{\beta-\alpha}$,  \textit{South East Asian J. of Math. and MAth. Sci.}, \textbf{11} (1)(2015), 17-24.
	
	\bibitem{GS.SAB.2016PF}	G. Shanker and S. A. Baby, On the Riemann and Ricci curvature of  Finsler spaces with special  $(\alpha, \beta)-$metric,  \textit{International J.  Math. Combin.}, \textbf{4} (2016), 61-69.
	
	
	\bibitem{GS.SAB.2017PFDF}	G. Shanker and S. A. Baby, On Kropina-Randers change of $m$th-root Finsler metric,  \textit{JP Journal of Geometry and Topology}, \textbf{20} (2)(2017), 115-128.
	
	\bibitem{GS.PG.2013PF}	G. Shanker and P. Gupta, Projectively Flat Finsler Spaces with
	Special $(\alpha, \beta)-$metric,  \textit{Journal of International Academy of Physical Sciences}, \textbf{17} (4)(2013), 369-376.
	
	\bibitem{GS.R.2013PF}	G. Shanker and Ravindra, On Randers change of exponential metric,  \textit{Applied Sciences}, \textbf{15} (2013), 94-103.
	
	\bibitem{GS.RKS.RDSK.2017PFDF}	G. Shanker, R. K. Sharma and R. D. S. Kushwaha, On the Matsumoto change  of a Finsler space with $m-$th root metric,  \textit{to appear in Acta Matematica Academie Paedagogiace Nyiregyh$\acute{a}$ziensis}, \textbf{32} (2)(2017).
	
	\bibitem{ShenRFGAIG} Z. Shen, Riemann-Finsler geometry with applications to information geometry, \textit{Chin. Ann. Math.}, \textbf{27B} (1)(2006), 73-94.
	
	\bibitem{C.Sh1984}	C. Shibata, On invariant tensors of $\beta-$change of Finsler metrics,  \textit{J. Math. Kyoto Univ.}, \textbf{24} (1984), 163-188.
	
	
	\bibitem{RY.GS.2013PF} R. Yadav and G. Shanker, On some projectively flat $(\alpha, \beta)-$metrics,  \textit{Gulf Journal of Mathematics}, \textbf{1} (2013), 72-77.	
\end{thebibliography}
\end{document}